\documentclass[11pt]{article}
\usepackage[utf8]{inputenc}

\usepackage[a4paper,margin=0.9in]{geometry}

\usepackage{graphicx} 
\usepackage[inline]{enumitem} 
\newcommand{\inlineitem}[1][]{%
\ifnum\enit@type=\tw@
    {\descriptionlabel{#1}}
  \hspace{\labelsep}%
\else
  \ifnum\enit@type=\z@
       \refstepcounter{\@listctr}\fi
    \quad\@itemlabel\hspace{\labelsep}%
\fi}
\makeatother

\usepackage[utf8]{inputenc}
\usepackage{amsmath}
\usepackage{amsfonts}
\usepackage{amssymb}
\usepackage{textcomp}
\usepackage{amssymb}
\usepackage{todonotes}
\usepackage{tikz-cd}
\usepackage{lipsum}
\usepackage{adjustbox}

\usepackage{float}
\usepackage{mathptmx}
\usepackage{stmaryrd}

\usepackage[utf8]{inputenc}
\usepackage{graphicx}
\usepackage{subfig}
\usepackage{mathtools}
\usepackage{tikz}
\usepackage[all,cmtip]{xy}
\usepackage{tikz-cd}
\usetikzlibrary{arrows}
\usetikzlibrary{positioning}
\usetikzlibrary{matrix}
\usepackage{authblk}

\usepackage[noadjust]{cite}
\usepackage{graphicx}
\usepackage{caption}
\usepackage{mathtools}
\usepackage{subcaption}
\graphicspath{ {./image/} }
\usepackage{amsthm}
\PassOptionsToPackage{hyphens}{url}\usepackage{hyperref}
\usepackage{cleveref}
\usepackage{amsthm}
\usepackage[dvipsnames]{xcolor}
\usepackage{tikz-cd}
\tikzcdset{scale cd/.style={every label/.append style={scale=#1},
    cells={nodes={scale=#1}}}}

\theoremstyle{remark}
\newtheorem*{remark}{Remark}

\theoremstyle{example}
\newtheorem{theorem}{Theorem}[section]
\newtheorem*{theorem*}{Theorem}

\newtheorem*{corollary*}{Corollary}

\newtheorem{Question}[theorem]{Question}
\newtheorem{lemma}[theorem]{Lemma}
\newtheorem{proposition}[theorem]{Proposition}
    
\theoremstyle{definition}
\newtheorem{definition}[theorem]{Definition}
\newtheorem{example}[theorem]{Example}

\newtheorem{theoremL}{Theorem}

\newtheorem{corollaryL}[theoremL]{Corollary}

\usepackage{graphicx}  
\usepackage{cite} 

\usepackage{lipsum}

\usepackage{graphicx}

\title{On the mapping class group of 4-dimensional 1-handlebodies via Budney--Gabai invariants}
\author{Weizhe Niu}
\date{}

\begin{document}

\maketitle

\begin{abstract}
We define an invariant $(W_3)_m$ for $\pi_0\mathrm{Diff}(\natural_m S^1\times D^3,\partial)$ for $m\geq 1$ that generalizes Budney--Gabai's $W_3$ invariant. We give a computational framework inspired by Budney--Gabai and use it to calculate the invariant for all unknotted barbell difeomorphisms of $\natural_m S^1\times D^3$ for $m=1,2$. This allows us to detect more linearly independent elements in $\pi_0\mathrm{Diff}(S^1\times D^3,\partial)$, and to prove that $\pi_0\mathrm{Diff}( \natural_2 S^1\times D^3,\partial)/ \left( \pi_0 \mathrm{Diff}(S^1\times D^3,\partial)\right)^2$ admits infinitely generated subgroups generated by unknotted barbell diffeomorphisms, leading to infinitely many properly embedded separating 3-balls that are non-isotopic relative to the boundary.
\end{abstract}

\section{Introduction}
This paper studies the mapping class group of $\natural_m S^1\times D^3$, $m\geq 1$, i.e.~4-dimensional 1-handlebodies. For $m=1$, Budney--Gabai \cite{Budney-gabai,Budney-gabai2} proved that $\pi_0\mathrm{Diff}(S^1\times D^3,\partial)$ is infinitely generated. They proposed a method for creating diffeomorphisms of a 4-manifold by embedding a ``barbell manifold'' inside of it and performing a construction similar to the ``point-pushing'' map on surfaces, and exhibited infinitely many isotopically non-trivial diffeomorphisms of $S^1\times D^3$. The main results of this paper are as follows.

\begin{theoremL}[cf.  Chapter 5]
For $m\geq 1$, there is an invariant
    $$(W_3)_m\colon \pi_0\mathrm{Diff}(\natural_m S^1\times D^3,\partial)\to \pi_5 C_3'[\natural_m S^1\times D^3]\otimes \mathbb{Q}/R$$
that depends on the choice of a properly embedded 3-ball $\Delta$ in $\natural_m S^1\times D^3$. For $m=2$, it detects an infinite family of linearly independent elements of $\pi_0\mathrm{Diff}(S^1\times D^3 \natural S^1\times D^3,\partial)/ \left(\pi_0 \mathrm{Diff}(S^1\times D^3,\partial)\right)^2$ represented by unknotted barbell diffeomorphisms.
\end{theoremL}
For $m=1,2$, we present formulae for $(W_3)_m$ for all unknotted barbells in $\natural_m S^1\times D^3$ (cf. Theorem \ref{barbellfactor} and Theorem \ref{barbellformulaem111} for $m=1$, and Lemma \ref{reallinearfactorisation}, Lemma \ref{type123formula} and Lemma \ref{type46formula} for $m=2$).
\begin{corollaryL}[cf. Theorem \ref{chapter5main1}]
    There exist infinitely many properly embedded non-isotopic separating 3-balls with common boundaries in $S^1\times D^3\natural S^1\times D^3$.
\end{corollaryL}

\begin{remark}
    Tatsuoka \cite{tatsuoka} produces an infinite sequence of non-isotopic 3-spheres in $S^1\times D^3\# S^1\times D^3$ which can lead to a sequence of non-isotopic separating 3-balls in $S^1\times D^3\natural S^1\times D^3$. However, the non-triviality of the separating 3-balls we produce in this paper do not follow from her argument as we will discuss in Section 5.3.
\end{remark}

For $m=1$, if we fix $\Delta=\{\mathrm{pt}\}\times D^3\subset S^1\times D^3$, then Budney--Gabai's $W_3$ invariant \cite{Budney-gabai,Budney-gabai2} is recovered. We generalize \cite{Budney-gabai,Budney-gabai2} by calculating $W_3$ for all unknotted barbell diffeomorphisms of $S^1\times D^3$. We detect more linearly independent elements of $\pi_0\mathrm{Diff}(S^1\times D^3,\partial)$ in addition to those in \cite{Budney-gabai}. In the following theorem, $\theta_k(e_i,e_j)$ denotes certain barbell diffeomorphisms of $S^1\times D^3$ (cf. Example \ref{firstexample} and Figure \ref{barbelltheta9}).

 \begin{theoremL}[cf. Theorem \ref{chapter4main}]
\label{chapter4main}
        The elements $\theta_k(e_{k-1},e_{k-3})$ for $k\geq 6$ of $\pi_0\mathrm{Diff}(S^1\times D^3,\partial)$ are linearly independent. Furthermore, these elements are linearly independent to $\delta_k\coloneqq\theta_k(e_{k-1},e_{k-2})$. More generally, there exist linearly independent elements $\theta_k(e_{k-1},e_{k-m})$, for $m\in \{3,4,\dots,[(k-1)/2]-1\}$ with $k\geq 2m-1$ where$[(k-1)/2]$ is the integer part of $(k-1)/2$.
    \end{theoremL}

The paper is structured as follows. Chapter 2 recalls the definition of barbell diffeomorphism. The invariant $(W_3)_m$ is defined in Chapter 3. Chapter 4 computes $(W_3)_1$ for unknotted barbell diffeomorphisms of $S^1\times D^3$ and Chapter 5 computes $(W_3)_2$ for unknotted barbell diffeomorphisms of $\natural_2S^1\times D^3$.
\vspace{3mm}

\par\noindent\textbf{Acknowledgement.} Most of this work is adapted from the author’s PhD thesis \cite{weizheniuthesis}, completed at the University of Glasgow under the supervision of Brendan Owens and Mark Powell. The author thanks them for their extensive comments and many helpful discussions. The author also thanks Ryan Budney for valuable conversations, and the author’s postdoctoral mentor Jianfeng Lin for helpful comments.
\section{Barbell diffeomorphism}

In this chapter, we recall the construction of barbell diffeomorphism of 4-manifolds. We keep the discussion concise and simple, and refer the reader to \cite{Budney-gabai,Budney-gabai2} or the author's PhD thesis \cite{weizheniuthesis} for details.

\begin{definition}
    The \textit{model barbell} in $\mathbb{R}^3\subset \mathbb{R}^4$ is the union of two 2-spheres of radius 1 centered at $(-2,0,0)$ and $(2,0,0)$ in $\mathbb{R}^3$ with the interval $[-1,1]\times (0,0)\subset \mathbb{R}\times \mathbb{R}^2$ connecting the two 2-spheres. The two 2-spheres are called the \textit{cuff spheres}, and the interval is called the \textit{bar}. The \textit{thickened model barbell} $\mathcal{B}$ is a closed neighbourhood of the model barbell in $\mathbb{R}^4$. It can be identified with the boundary connected sum of two trivial $D^2$-bundles over $S^2$ which is naturally isomorphic to $(S^2\times B^1 \natural S^2\times B^1)\times I$, and we denote the $I$-coordinate by $t$. In other words, this can be obtained by taking a tubular neighbourhood of the model barbell in $\mathbb{R}^3$ and taking the product with $I$. 
\end{definition}

    \begin{figure}
  \centering
  \begin{minipage}[b]{0.4\textwidth}
\begin{tikzpicture}[scale=0.90]
\filldraw[color=black!100, fill=black!0,thick](-3,0) circle (1.5);
\filldraw[color=black!100,fill=black!0, thick](3,0) circle (1.5);
\draw[black, thick] (-1.5,0) -- (1.5,0);
\draw[black,thick] (-4.5,0) .. controls (-3,-0.5) .. (-1.5,0);
\draw[black,dashed] (-4.5,0) .. controls (-3,0.5) .. (-1.5,0);
\draw[black,thick] (4.5,0) .. controls (3,-0.5) .. (1.5,0);
\draw[black,dashed] (4.5,0) .. controls (3,0.5) .. (1.5,0);
    \end{tikzpicture}
    \centering
    \caption{The model barbell.}
    \label{standardmodelbarbell}
  \end{minipage}
  \hfill
  \begin{minipage}[b]{0.4\textwidth}
\begin{tikzpicture}[scale=0.90]
\draw[black,thick] (-3,1.5) .. controls (0,0.5) .. (3,1.5);
\draw[black,thick] (-3,1.5) .. controls (0,2.5) .. (3,1.5);
\draw[black,thick] (-3,-1.5) .. controls (0,-2.5) .. (3,-1.5);
\draw[black, thick] (-3,1.5) -- (-3,-1.5);
\draw[black, thick] (3,1.5) -- (3,-1.5);
\filldraw[color=black!100, fill=black!0,thick](-1.5,-0.5) circle (0.5);
\draw[black,thick] (-2,-0.5) .. controls (-1.5,-0.75) .. (-1,-0.5) ;
\filldraw[black] (-1.5,-1.7) circle (0pt) node[anchor=south]{$(\alpha_1)_0$};
\draw[dotted,black, thick] (-1.5,0) -- (-1.5,1.5) node[anchor=west]{$E_1$};
\filldraw[color=black!100, fill=black!0,thick](1.5,-0.5) circle (0.5);
\filldraw[black] (1.5,-1.7) circle (0pt) node[anchor=south]{$(\alpha_2)_0$};
\draw[black,thick] (2,-0.5) .. controls (1.5,-0.75) .. (1,-0.5);
\draw[dotted,black, thick] (1.5,0) -- (1.5,1.5) node[anchor=west]{$E_2$};
\filldraw[black] (5,0) circle (0pt) node[anchor=east]{$\times [-1,1]$};
  \draw [-stealth](0,-0.5) -- (0.5,-0.5) node[anchor=south]{$z$};
  \draw [-stealth](0,-0.5) -- (0,0) node[anchor=west]{$y$};
    \draw [-stealth](0,-0.5) -- (0.5,0) node[anchor=west]{$x$};
\end{tikzpicture}
\centering
\caption{The thickened model barbell.}
\label{barbellcylinder}
  \end{minipage}
\end{figure}

\begin{figure}

\includegraphics[width=0.4\textwidth]{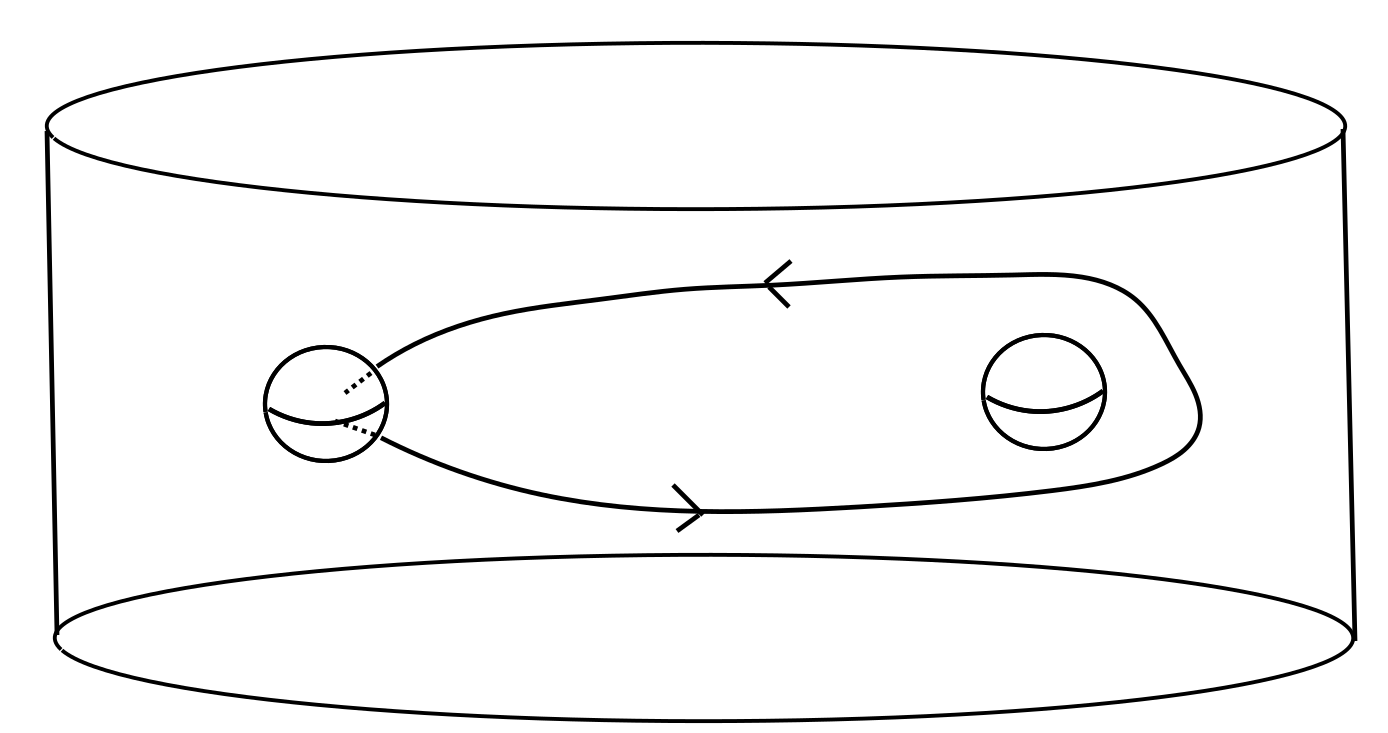}
\centering
\caption{The $t=0$ slice of the thickened model barbell.}
\label{barbell}

\end{figure}

\begin{figure}
    \centering
    \includegraphics[width=0.3\linewidth]{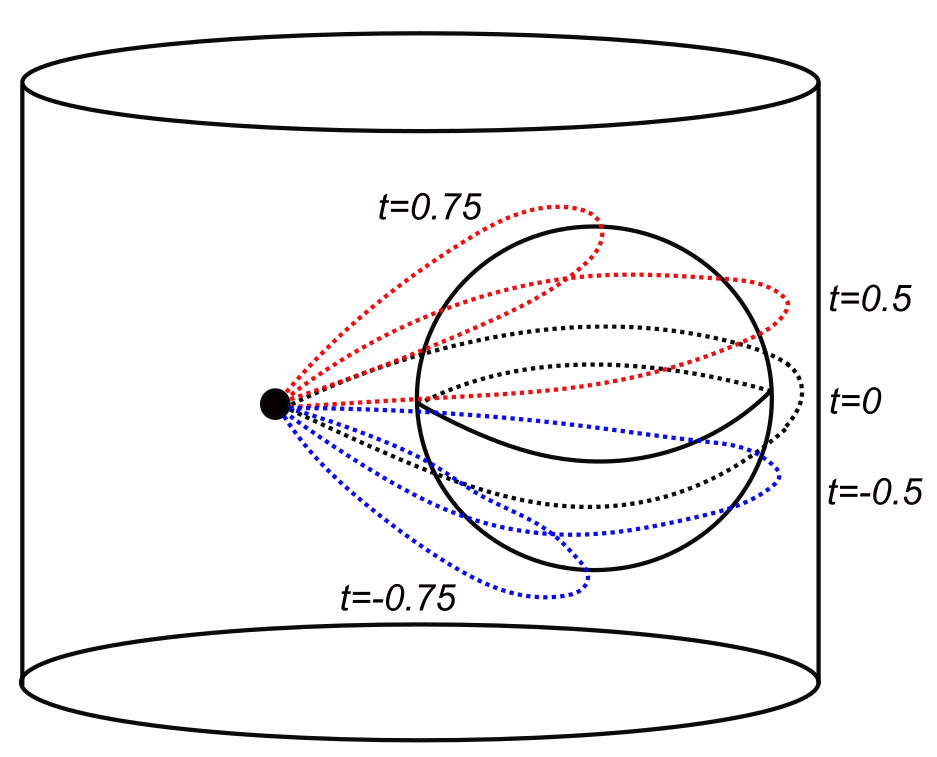}
    \caption{The time-dependent path in the barbell that defines the barbell diffeomorphism via isotopy extension.}
    \label{barbellpathtime}
\end{figure}

Denote the two complementary 4-balls of the thickened model barbell $\mathcal{B}$ by $\alpha_1$ and $\alpha_2$; both are parameterized by $I\times B^3$. For $t\in[-1,1]$, let $(\alpha_i)_t$ denote the intersection between $\alpha_i$ and the $t$-slice $\mathcal{B}_t$ of $\mathcal{B}$. The two slices $(\alpha_1)_0$ and $(\alpha_2)_0$ are shown in Figure \ref{barbellcylinder}. The two 2-disks $E_1$ and $E_2$ represent two properly embedded orthogonal 2-disks to the 2-spheres $\partial (\alpha_1)_0$ and $\partial (\alpha_2)_0$. More precisely, if we give $\mathcal{B}$ a standard handle structure with one 0-handle and two 2-handles, then they can be viewed as the cocores of the two 2-handles. At the $t=0$ slice of $\mathcal{B}$, push the 3-ball $(\alpha_1)_0$ along a closed loop around $(\alpha_2)_0$, but not touching  $(\alpha_2)_0$, as in Figure \ref{barbell}. This leads to a path of embedded 3-balls in $(\mathcal{B})_0\cup (\alpha_1)_0$. As $t$ approaches $\pm 1$, swing the path to the upward/downward direction, tracing out two hemispheres, until it reaches the trivial arc orthogonal to $\alpha_1$ (parallel to the normal 2-disk $E_1$ to the 2-sphere $\partial (\alpha_1)_0$). See Figure \ref{barbellpathtime} for the loops when $t=\pm 0.5$ and $t=\pm 0.75$. Note that we draw the left 2-sphere  (and the complementary 3-ball $(\alpha_1)_0)$ it bounds) as a black dot to make it easier to visualise. The above defines a loop of embedded 4-balls, or more precisely, a loop of string links (i.e.~a loop of embeddings of $I\times D^3$), in $\mathcal{B}\cup \alpha_1$ based at $\alpha_1$ that is $t$-level-preserving. We also arrange so that it is away from the boundary of $\mathcal{B}$. Now, applying the isotopy extension theorem to this loop in the manifold $\mathcal{B}\cup \alpha_1$ gives rise to a well-defined isotopy class of boundary-preserving diffeomorphisms of $\mathcal{B}$, which only depends on the isotopy class of this $t$-level preserving path by the uniqueness of the isotopy extension theorem. This is called the \textit{barbell diffeomorphism}. 

\begin{definition}
    Let $X$ be an orientable 4-manifold. Let $\mathcal{B}$ be an embedded barbell in $X$ whose cuffs have trivial normal bundles. The induced barbell diffeomorphism $\Phi_{\mathcal{B}}\in \mathrm{Diff}(X,\partial)$ is defined by extending the barbell diffeomorphism using the identity map to the rest of the ambient manifold $X$, and it gives rise to an element $[\Phi_{\mathcal{B}}]\in \pi_0\mathrm{Diff}(X,\partial)$. This is called a \textit{barbell implantation}.
\end{definition}

Since $X$ is orientable and the embedded cuff spheres admit unique framings that agree with the orientation of $X$, such an embedding is determined by embedding two disjoint 2-spheres together with a framed embedded bar connecting them. As argued in Remark 5.12 in \cite{Budney-gabai}, different framings of the bar give rise to different barbell implantations that are related by full right-hand twists of the $x$-$y$ plane as one travels along the bar. Such twists preserve isotopy classes of barbell implantations since they fix the union of $E_1$ and its image under the barbell diffeomorphism setwise (see Remark 5.4 of \cite{Budney-gabai}) for details). Therefore, when talking about barbell implantations, we usually do not distinguish between an embedded barbell and an embedded thickened model barbell. We sometimes omit the word embedded as well.

\begin{definition}
    A barbell in $X$ is called \textit{unknotted} if its two cuff spheres are both isotopically unknotted. By an isotopically unknotted 2-sphere in $X$ we mean an embedded 2-sphere that bounds a $D^3$.
\end{definition}

We consider the case when $X=S^1\times D^3$. An unknotted barbell in $S^1\times D^3$ is determined by the relative isotopy class of the bar. Namely, if we choose two embedded unknotted spheres $B$ and $R$ in $X$, together with two points $b_0$ on $B$ and $r_0$ on $R$, then the space of unknotted barbells is determined by isotopy classes of embeddings $\pi_0\mathrm{Emb}(I,X;b_0,r_0)$ which can be described by a word in the free group $F_3$ with 3 generators: the meridians $\nu_R$ and $\nu_B$  of the two spheres, and the circle factor $t$ of $S^1\times D^3$. We orient the bar by saying that it starts from $B$ and ends at $R$. In other words, the induced barbell diffeomorphism is obtained by looping $B$ around $R$. In fact, the opposite orientation gives the inverse of the barbell diffeomorphism, which can be seen by directly tracking the definitions.

If the bar links the sphere $B$ before looping around the $S^1$-factor at all, or links the sphere $R$ just before reaching $R$ in the end, then these linkings can be eliminated by an isotopy that directly drags the spheres $B$ and $R$ out of the parts of the bar near the two ends (or equivalently, through an isotopy that straightens the bar). Therefore, such barbells are completely determined by the double coset $$\langle \nu_B \rangle\backslash\langle \nu_B, \nu_R,t\rangle/\langle \nu_R \rangle.$$ 
In other words, we do not need to consider words that start with $\nu_B$ or end with $\nu_R$.

\begin{figure}
    \centering
    \includegraphics[width=0.45\textwidth]{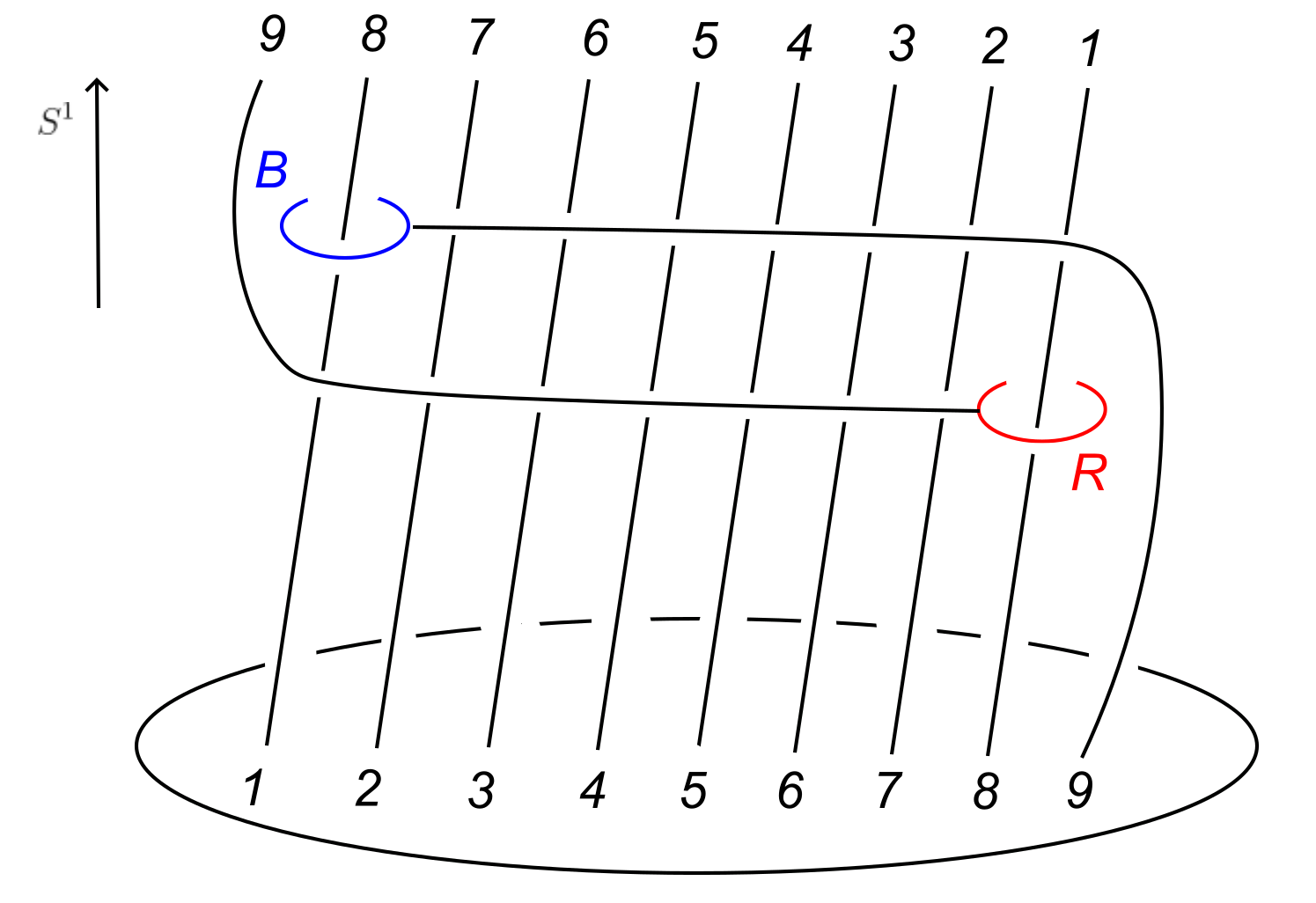}
    \caption{The embedded barbell $\theta_{10}(e_8,e_8)\coloneqq \theta_{10}((0,\dots,0,1,0),(0,\dots,0,1,0))$. Strand number $k$ on the top is glued to strand number $10-k$ on the bottom.}
    \label{barbelltheta9}
\end{figure}
\begin{example}
\label{firstexample}
 There is a subclass of such barbells studied in \cite{Budney-gabai,Budney-gabai2}, which will play an important role in this paper. For $k \in \mathbb{Z}^{+}$ and $(v)_i, (w)_i \in \mathbb{Z}^{k-1}$ where $i=1,2,\dots,k-1$, let $\theta_k(v, w)$ denote the following barbell in $S^1\times D^3$: choose two disjoint, parallel embedded 2-spheres $B\subset \{s_2\}\times D^3\subset S^1\times D^3$, $R$ in 
    $\{s_1\}\times D^3\subset S^1\times D^3$ with $s_1\neq s_2$. We specify a bar connecting them. The integer $k-1$ indicates the number of times the bar moves around the $S^1$ factor in total, and we require the bar to go in the negative $S^1$ direction (we give the product orientation to $S^1\times D^3$) only throughout. If we cut the bar with a 3-ball $\{s_0\}\times D^3\subset S^1\times D^3$ with $s_0\notin \{s_1,s_2\}$ such that the order of the triple $(s_0,s_1,s_2)$ agrees with the orientation of $S^1$, then there are $k-1$ intersection points between the bar and $\{s_0\}\times D^3$. We arrange these intersection points to lie on the same line in $\{s_0\}\times D^3$ and give them indices $\{1,2,\dots,k-1\}$ as shown in Figure \ref{barbelltheta9} for $\theta_{10}((0,\dots,0,1,0),(0,\dots,0,1,0))$. Note that for better compatibility with Budney--Gabai's work, we follow the convention used in \cite{Budney-gabai2} that we label the strands in opposite directions at the top and bottom. So, the journey of the bar from $B$ to $R$ leads to $k-1$ vertical (i.e.~parallel to the circle direction) strands with indices and $R$ can link each of these strands. We use each entry of the vector $(v)_i$ with $i\geq 1$ to indicate the signed number of times the $i$-th vertical strand of the bar wraps around the sphere $R$, counting from left to right at the bottom of Figure \ref{barbelltheta9}. We fix the convention such that a vertical strand going through a cuff sphere from top to bottom (i.e.~pointing to the negative $S^1$ direction) is positive. Similarly, for the blue cuff $B$, following the convention as in \cite{Budney-gabai2}, we give indices to the same strands, but with reversed order. These are shown in Figure \ref{barbelltheta9} on the top. Each vector entry $(w)_i$ with $i\geq 1$ indicates the signed number of times the $i$-th vertical strand of the bar wraps around the sphere $B$ counting from right to left on the top. When both $B$ and $R$ link the same strand, we only consider the barbell such that such a strand first links $B$, and then links $R$. Note that we always keep $B$ in $\{s_2\}\times D^3$ and $R$ in $\{s_1\}\times D^3$. In the language of free group words, $\theta_k(v, w)$ corresponds to the word 
$$\nu_R^{v_{k-1}}t^{-1}\nu_B^{w_{1}} \nu_R ^{v_{k-2}} t^{-1}\dots t^{-1}\nu_B^{w_{k-2}}\nu_R ^{v_1} t^{-1} \nu_B ^{w_{k-1}}.$$
Let $e_i\in \mathbb{Z}^{k-1}$ be the vector whose only non-zero entry is the $i$-th entry and is $1$. As we will see in Section 4.3, the barbells $\theta_k(e_i,e_j)$ will be important as building blocks of barbells in $S^1\times D^3$.

\end{example}

\section{$(W_3)_m$ invariant for $\pi_0\mathrm{Diff}(\natural_mS^1\times D^3,\partial)$}
Budney--Gabai \cite{Budney-gabai, Budney-gabai2} defined a rational homotopy invariant $W_3$ for 
$\pi_0\mathrm{Diff}(S^1\times D^3,\partial)$ and used it to detect isotopically non-trivial diffeomorphisms of $S^1\times D^3$. In this chapter, we define an invariant
$$(W_3)_m\colon \pi_0\mathrm{Diff}(\natural_m S^1\times D^3,\partial)\to \pi_2\mathrm{Emb}(I,\natural_m S^1\times D^3;I_0) \to \pi_5 C_3'[\natural_m S^1\times D^3]\otimes \mathbb{Q}/R.$$
for $m\geq 1$ as a composition of two maps, namely the scanning map $s$, which depends on a choice of a properly embedded 3-ball in $\natural_m S^1\times D^3$, and a map defined on $\pi_2\mathrm{Emb}(I,\natural_m S^1\times D^3;I_0)$ that relies on the mapping space model \cite{sinha}. Here $I_0$ is a chosen properly embedded interval, which we will elaborate on shortly, and $R$ represents some relations that will become clear soon. We will illustrate that when $m=1$ and for a natural choice of a properly embedded 3-ball in $S^1\times D^3$, we recover Budney-Gabai's original $W_3$ invariant.

\subsection{Fulton-MacPherson compactification}

In this section, we briefly introduce the Fulton-MacPherson compactification of configuration spaces of manifolds, mostly following \cite{sinha2} and \cite{sinha}. We will need this notion for the definition of $(W_3)_m$. We present the results and properties that we need and omit proofs. In the end of this section, we introduce the \textit{mapping space model} due to Dev Sinha \cite{sinha} for the space of properly embedded intervals with fixed endpoints in a 4-manifold $M$. For readers who are already familiar with these notions, this section serves the purpose of setting up our notation, and one can skip the details otherwise.

\begin{definition}
    
Let $M$ be a compact manifold embedded in some Euclidean space $\mathbb{R}^{N+1}$, and let $C_n(M)$ denote the $n$-th ordered configuration space of $M$, i.e.~the space of ordered, distinct $n$-tuples of points in $M$. We define the \textit{Fulton-MacPherson compactification} of $C_n(M)$ by $C_n[M]=\mathrm{Bl}_{\Delta } M^n$, i.e.~the blow-up along the fat diagonal (the subspace in which at least two points coincide). More precisely,
  \begin{itemize}
  \item For $m,n\geq 1$, let $C_m(n)$ be the space of ordered $m$-tuples of distinct points in the indexing set $\{1,2,\dots,n\}$. Note that $C_m(n)$ is empty if $m>n$.
      \item For $(i,j)\in C_2(n)$, let $\pi_{ij}\colon C_n(M)\to S^N$ be the map that sends $(x_i)\in M^n\subset (\mathbb{R}^{N+1})^n$ to the unit vector in the direction of $x_i-x_j$.
      \item  For $(i,j,k)\in C_3(n)$, let $s_{ijk}\colon C_n(M)\to [0,\infty]$ be the map that sends $(x_i)$ to $(|x_i-x_j|/|x_i-x_k|)$.
      \item Let $A_n[M]$ be the product $M^n\times (S^N)^{C_2(n)}\times [0,\infty]^{C_3(n)}$.
      \item Define $\alpha_n\colon C_n(M)\to A_n[M]$ to be the product $\mathrm{Inc}\times (\pi_{ij})\times (s_{ijk})$ where $\mathrm{Inc}\colon C_n(M)\to M^n$ is the natural inclusion.
  \end{itemize}
 Then $C_n[M]$ is defined as the closure of the image of $C_n(M)$ in $A_n[M]$. We summarize some of the important properties of $C_n[M]$ following the first few sections of \cite{sinha}.
\begin{itemize}
\item The homeomorphism type of $C_n[M]$ is independent of the embedding of $M$ in $\mathbb{R}^{N+1}$.

    \item If $M$ is compact, then $C_n[M]$ is compact.
    \item The inclusion $\mathrm{Inc}\colon C_n(M)\to M^n$ factors through a surjective projection $C_n[M]\to M^n$.
    \item An embedding $f\colon M\to N$ induces a map $ev_n(f)\colon C_n[M]\to C_n[N]$ extending the natural induced map on $C_n(M)$.
    \item $C_n[M]$ is a manifold with corners with $C_n(M)$ as its interior, equipped with a preferred stratification structure.
    \item The inclusion $C_n(M)\to C_n[M]$ is a homotopy equivalence.
\end{itemize}

\end{definition}

\begin{example}
    For $M=I$, the space $C_2[I]$ is a union of two triangles. More generally, $C_n[I]$ is a disjoint union of $n$-simplices, having one component for each ordering of $n$ points.
\end{example}

There is a preferred \textit{stratification} structure of $C_n[M]$ that can be defined combinatorially in the following way. For us, a stratification of a space $X$ is a collection of disjoint subspaces $\{X_c\}$ called strata, such that the intersection of the closures of any two strata is the closure of some stratum.

\begin{definition}
Let $\Phi_n$ denote the category of rooted, connected \textit{trees}, with $n$ leaves labelled by 1,\dots,$n$, and with no bivalent internal vertices. Each tree admits a natural orientation when a root vertex is fixed by defining the direction pointing away from the root as positive. For $T$ and $T'\in \Phi_n$, there is a unique morphism between $T$ and $T'$ if they are isomorphic up to contraction of some non-leaf edges. Two leaves are called \textit{root-joined} if the unique paths to the root vertex intersect only at the root vertex.
\end{definition}
\begin{definition}
         Given a set $S$, an \textit{exclusion relation} $R$ is a subset of $C_3(S)$ satisfying
\begin{itemize}
\item if $(x,y,z)\in R$, then $(y,x,z)\in R$ and $(x,z,y)\notin R$;
\item if $(x,y,z)\in R$ and $(w,x,y)\in R$, then $(w,x,z)\in R$.

\end{itemize}
The collection of all exclusion relations of $S$ is denoted by $\mathrm{Ex}(S)$.
\end{definition}

\begin{definition}
    A \textit{parenthesization} of a set $S$ is a collection of nested subsets of $S$ with each of them having cardinality greater than one.
\end{definition}
There exists a map $\mathcal{F}_n$ from $\mathrm{Ex}(n)$ of the indexing set $\{1,2,\dots,n\}$ to the collection of parenthesizations of the same set.  For $R\in \mathrm{Ex}(n)$, this is given by taking the following sets.
\begin{itemize}
    \item $A_{\sim i,\text{¬}k}$ that contain all $j$ such that $(i,j,k)$ is in $R$ and also $i$ if such a $j$ exists.
\end{itemize}

\begin{definition}
\label{tree}
    For a parenthesization $P$ of $\{1,2,\dots,n\}$ given by a collection $\{A_\alpha\}$ of nested subsets, there is a tree $T(P)\in \Phi_n$ defined in the following way:
 \begin{itemize}
    \item Each $A_\alpha$ gives an internal vertex $v_\alpha$.
    \item There is an edge between $v_\alpha$ and $v_\beta$ if $A_\alpha \subset A_\beta$ but there are no proper inclusions $A_\alpha \subset A_\gamma \subset A_\beta$.
    \item A root vertex with edges connecting it to all internal vertices corresponding to maximal $A_\alpha$.
    \item Leaves (i.e.~vertices without edges with positive orientation) with an edge connecting the $i$-th leaf to either the vertex $v_\alpha$ where $A_\alpha$ is the minimal set containing $i$, or to the root vertex if no such $A_\alpha$ exists.

    \end{itemize}
\end{definition}

For $x=(x_i,u_{ij},d_{ijk})\in C_n[M]$, we can define an exclusion relation $R(x)\in \mathrm{Ex}(n)$ by letting $(i,j,k)\in R(x)$ if $d_{ijk}=0$. Then Definition \ref{tree} gives rise to a tree denoted by $T(x)$. For $T\in \Phi_n$, let $C_T(M)$ denote the space of points $x\in C_n[M]$ such that $T(x)=T$ and let $C_T[M]$ be its closure. This defines a stratification of $C_n[M]$ (see Section 3 of \cite{sinha}). 

We will also make use of another space $C_k'[M]$, defined as the pullback of the $k$-fold product of the unit tangent bundle of $M$ to $C_k[M]$. Intuitively, it contains points decorated with tangent vectors. It is equipped with a natural stratification induced from that on $C_k[M]$. The following diagram describes their relationship:
\[
\begin{tikzcd}
C_k'[M] \arrow{r} \arrow[swap]{d} & (STM)^k \arrow{d} \\
C_k[M]  \arrow{r} & M^k.
\end{tikzcd}
\]
There is one more variant of  $C_n[M]$ as introduced in \cite{Budney-gabai, Budney-gabai2, sinha} that is involved in Budney--Gabai's invariants.

\begin{definition}
    
Let $C_n(M,\partial)$ be the subspace of $C_{n+2}(M)$ such that the first and last points are fixed (and distinct) points $x_0$ and $x_1$ in the boundary of $M$. The compactification of $C_n(M,\partial)$, denoted by $C_n[M,\partial]$, is defined as the closure of the image of the map $$\mathrm{Inc}\times (\pi_{ij}) \times (s_{ijk})\colon C_{n+2}(M)\to A_{n+2}[M]$$
restricting to the subdomain $C_n(M,\partial)$. Similarly, let $C_n'[M,\partial]$ be the subspace of $C_{n+2}'[M]$ that maps down to $C_n[M,\partial]$.

\end{definition}

There is a stratification on $C_n[M,\partial]$, hence on $C_n'[M,\partial]$ defined in a similar way as in Definition \ref{tree} using a (full) subcategory $\Phi_n'$ of $\Phi_n$. This subcategory contains trees whose roots have valence at least 2 and such that each set of leaves over the same vertex is \textit{consecutive} for all non-root vertices, meaning that if the indices $i$ and $j$ are leaves over the same vertex, then $i<k<j$ implies that $k$ is also a leaf over this vertex. Here, a leaf is said to be \textit{over} a vertex if this vertex is contained in the unique path from the leaf to the root.

From now on, we will work with the \textit{linearly ordered component} of $C_n'[I,\partial]$ for which the order of the points in the configuration agrees with the order they occur in the interval, and with all vectors pointing to the positive direction. We still use the same notation $C_n'[I,\partial]$, but the reader should keep in mind that equivalent constructions can be performed in other components as well.

\begin{figure}
      \centering
      \includegraphics[width=0.35\linewidth]{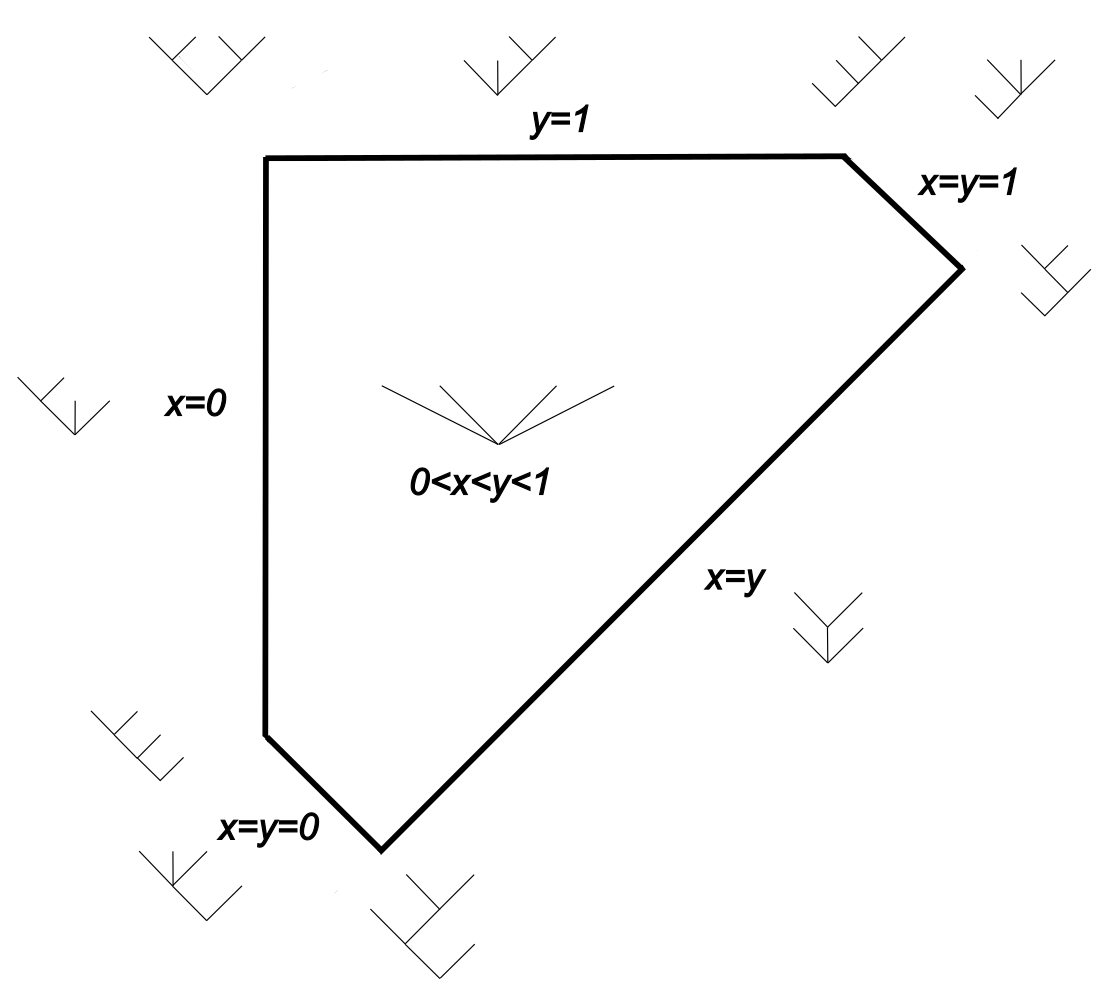}
      \caption{The space $C_2[I,\partial]$ as the fourth Stasheff polytope with trees annotated.}
      \label{polytope}
  \end{figure}

\begin{example}
    Consider $C_2[I,\partial]$ as the closure of the image of the map $$\mathrm{Inc}\times (\pi_{ij}) \times (S_{ijk})\colon C_4[I]\to A_4[I]$$
    restricting to $C_2(I,\partial)=\{(0,x,y,1)\colon 0<x<y<1\}$ where we restrict to one component only. This space can be parameterized in the following way: 
    \begin{itemize}
        \item The interior contains points with $0<x<y<1$
        \item Three boundary faces $0<x=y<1$, $x=0$ and $y=1$
        \item Two extra faces denoted by $x=y=0$ and $x=y=1$.
    \end{itemize}
 Following the bijection $\mathcal{F}_4$, we list the largest subset in each of the corresponding parenthesizations:
\begin{itemize}
    \item The interior corresponds to $\{1,2,3,4\}$
    \item The three (standard) boundary faces (with the vertices removed) correspond to $\{2,3\}$, $\{1,2\}$ and $\{3,4\}$
    \item The two extra faces correspond to $\{1,2,3\}$ and $\{2,3,4\}$
    \item The vertex $x=0, y=1$ corresponds to two sets $\{1,2\}$ and $\{3,4\}$.

\end{itemize}
          It is parameterized by the fourth Stasheff polytope as shown in Figure \ref{polytope}. The corresponding trees at different strata are also shown in Figure \ref{polytope}. In fact, the $(n+2)$-th Stasheff polytopes also parameterize higher-dimensional compactifications $C_n[I,\partial]$. The reader can refer to Theorem 4.19 in \cite{sinha} for more details.
\end{example}

\begin{definition}
Let $x=(x_i,v_i)\times (u_{ij})\times (s_{ijk})\in C_T'[M,\partial]$ where $v_i\in T_{x_i}M$. Then $x$ is called \textit{aligned} with respect to $T$ if for all $i$ and $j$ that are not root-joined (thus we have $x_i=x_j)$, we have $v_i=v_j$ and $u_{ij}$ is the image of $v_i$ under the Jacobian of the embedding of $M$ in $\mathbb{R}^N$. The subspace of aligned points is called the \textit{aligned sub-stratum} and is denoted by $C_T^{\alpha}[M,\partial]$.
\end{definition}

\begin{definition}
    A map $f\colon C_n'[I,\partial]\to C_n[M,\partial]$ that respects the substratification by aligned points, meaning that $C_T^{\alpha}[I,\partial]$ gets mapped to $C_T^{\alpha}[M,\partial]$, is called \textit{stratum-preserving and aligned}.
\end{definition}

The Goodwillie-Weiss embedding calculus tower (see \cite{10.2140/gt.1999.3.67}) provides an approximation of the space $\mathrm{Emb}(I,M)$ of smooth proper embeddings of the interval in $M$ as shown in Figure \ref{embedding calculus}. The maps $ev_k$ are called the $k$-th \textit{evaluation maps}. We state the \textit{mapping space model}.

\begin{theorem}[Mapping space model, Chapter 5 of \cite{sinha}]
    
The space of aligned stratum-preserving maps is denoted by $$\mathrm{Map}(C_n'[I,\partial],C_n'[M,\partial])$$ and is weakly homotopy equivalent to the $n$-th element $T_n \mathrm{Emb}(I,M)$ in the tower of embedding calculus, as long as the dimension of $M$ is 4 or greater. Furthermore, the composition of this (weak) homotopy equivalence with the $k$-th evaluation map (which we still denote by $ev_k$) is given by sending $f$ to the map $$(x_1,x_2,\dots,x_n)\mapsto (f(x_1),\dots,f(x_n),f'(x_1),\dots,f'(x_n))$$ with $f'(x_i)$ unit-normalized. Note that we omit the tangent vectors in $C_n'[I,\partial]$ as they are all positive.

\end{theorem}

\begin{figure}
    \centering

\centering
\[
\xymatrix{
\mathrm{Emb}(I,M) \ar[r]^{ev_k} \ar[dr]_{ev_{k-1}} & T_k \mathrm{Emb}(I,M)\ar[d] \\
 & T_{k-1} \mathrm{Emb}(I,M)}
\]

\caption{The embedding calculus tower.}
\label{embedding calculus}
\end{figure}

\subsection{Definitions and construction}
In the following, we set up preparations for the definition of an invariant
$$(W_3)_m\colon \pi_0\mathrm{Diff}(\natural_m S^1\times D^3,\partial)\to \pi_2\mathrm{Emb}(I,\natural_m S^1\times D^3;I_0) \to \pi_5 C_3'[\natural_m S^1\times D^3]\otimes \mathbb{Q}/R.$$
For $n\geq 1$, the projection map
$$C_n(\natural_m S^1\times D^3)\to C_{n-1}(\natural_m S^1\times D^3)$$ 
is a fibration that admits a section (see for example Section 4 of \cite{cohen2010introduction}), therefore the homotopy groups of $C_n(\natural_m S^1\times D^3)$ are isomorphic to the homotopy groups of the product 
$X_0\times X_1\times\cdots\times X_{n-1}$
where $X_i$ is $\natural_m S^1\times D^3$ with $i$ punctures. In other words, we can understand the homotopy groups of $C_n(\natural_m S^1\times D^3)$ by understanding the homotopy groups of $$(\vee_m S^1)\times (\vee_m S^1\vee S^3)\times (\vee_m S^1\vee S^3\vee S^3)\times \cdots\times (\vee_m S^1\vee_{n-1} S^3).$$ 
Therefore, $\pi_1 C_n(\natural_m S^1\times D^3)$ is generated by $nm$ generators $$(t_1)_i, (t_2)_i,\dots (t_n)_i$$ for $i=1,2,\dots, m$. The group 
$\pi_3 C_n(\natural_m S^1\times D^3)\otimes \mathbb{Q}$ is generated by elements of the form of $(t_l)_i^s . \omega_{jk}$ for $s\in \mathbb{Z}$, $i=1,2,\dots,m$ and $l$, $j$, $k\in \{1,2,\dots,n\}$ with the following relations satisfied:
\begin{itemize}
    \item $\omega_{ij}=\omega_{ji}$ and $\omega_{ii}=0$
    \item $(t_l)_i .\omega_{jk}=\omega_{jk}$ if $l\notin \{j,k\}$
    \item $(t_j)_i .\omega_{jk}=(t_k)_i^{-1}.\omega_{jk}.$
\end{itemize}
Here $\omega_{ij}$ denotes the element with all points fixed except the $j$-th point which orbits around the $i$-th point along a 3-sphere.

Recall that the \textit{Whitehead product} operation for a topological space $X$ is a graded quasi-Lie algebra structure on the homotopy groups of $X$. For $f\in \pi_k X$ and $ g\in \pi_l X$, we have
$$[f,g]\colon S^{k+l-1}\to S^k\vee S^l\to X,$$
where the first map is given by the attaching map of the top cell of $S^k\times S^l$ and the second map is the wedge $f\vee g$.  A basic example is $X=S^2$ and $f=g=\mathrm{Id}\in \pi_2 S^2$. The Whitehead product is twice the Hopf map $2 \nu\colon S^3\to S^2$. The group $\pi_5 C_n(\natural_m S^1\times D^3)\otimes \mathbb{Q}$ is generated by the Whitehead products of the above elements with the following extra relations satisfied:
\begin{itemize}
    \item $[\omega_{ij},\omega_{lm}]=0$ if $\{i,j\}\cap \{l,m\}=\emptyset$
    \item $[\omega_{ij},\omega_{jl}]=[\omega_{jl},\omega_{li}]=[\omega_{li},\omega_{ij}]$
    \item $(t_l)_i.[f,g]=[(t_l)_i.f,(t_l)_i. g].$
\end{itemize}

In the following, we define the \textit{scanning map} 
$$s\colon \pi_0\mathrm{Diff}(\natural_m S^1\times D^3,\partial)\to \pi_2\mathrm{Emb}(I,\natural_m S^1\times D^3;I_0)$$ 
which is only defined up to a choice of a properly embedded 3-ball in $\natural_m S^1\times D^3$. In particular, let $D_s\subset\natural_m S^1\times D^3$ be a properly embedded 3-disk called the \textit{scanning disk}. We choose our basepoint interval $I_0=0\times 0 \times [-1,1]\subset D_s\cong [-1,1]\times [-1,1] \times [-1,1]$. An element $[\Phi]\in \pi_0\mathrm{Diff}(\natural_m S^1\times D^3,\partial)$ represented by a diffeomorphism $\Phi$ maps the intervals $u\times v\times [-1,1]\subset D_s$, $u,v\in[-1,1]$ to a two-parameter family of properly embedded intervals, leading to an element in $\pi_2\mathrm{Emb}(I,\natural_m S^1\times D^3;I_0)$ after connecting the end points with the endpoints of $I_0$ by straight arcs in the boundary faces $I \times I\times 0$ and $I\times I\times 1$. In other words, a double loop of embedded intervals. An isotopy of $\Phi$ induces homotopies of the image $\Phi(u\times v\times [-1,1])\subset D_s\subset \natural_m S^1\times D^3$ relative to the boundary faces. Therefore, the map $s$ is well-defined up to isotopy.

\begin{figure}
    \centering
    \includegraphics[width=0.5\linewidth]{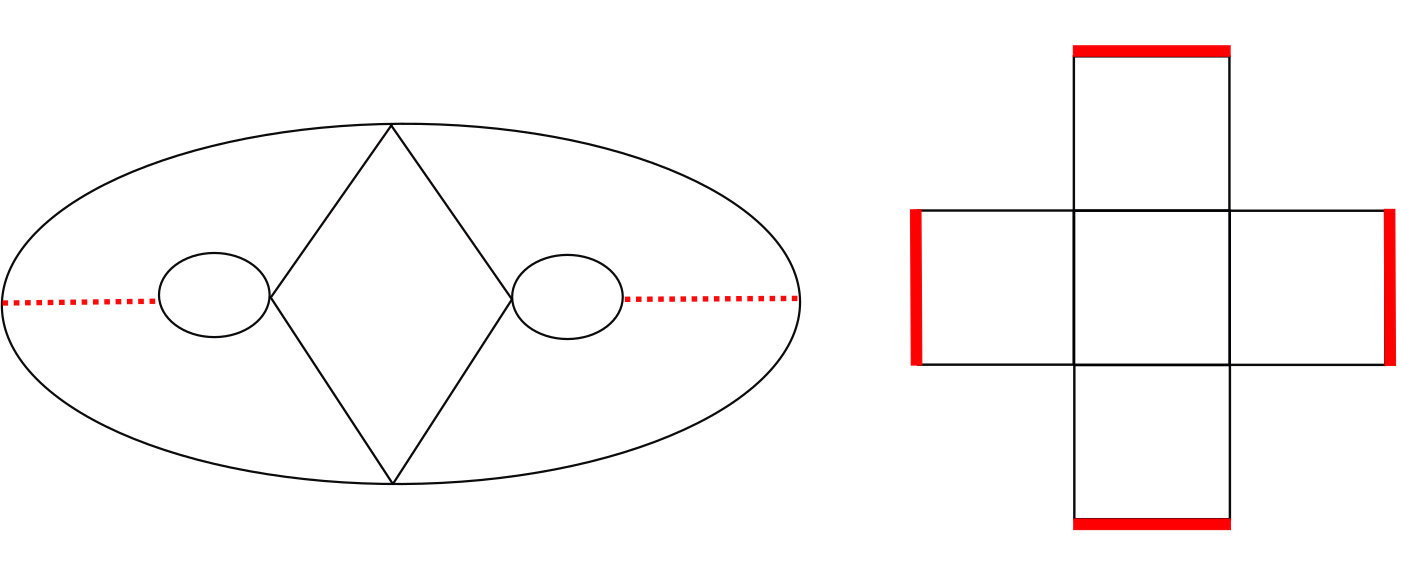}
    \caption{The fundamental domain of the universal cover of $S^1\times D^3\natural S^1\times D^3$.}
    \label{fundamental domain}
\end{figure}

Before defining the second map of $(W_3)_m$, we first take a look at the universal cover of $\natural_m S^1\times D^3$. Remove $m$ disjoint sub-disks from $D^2$ and cut along a disjoint union of homotopically trivial (relative to the boundary) arcs connecting each of the sub-disks to the boundary of $D^2$. See Figure \ref{fundamental domain} for the case $m=2$. This gives us a fundamental domain. We then build a space $T_m$ by gluing infinitely many fundamental domains together along their cut edges (drawn in red in Figure \ref{fundamental domain}). The (4-dimensional) universal cover $T_m\times I\times I$ is obtained by taking product with $I^2$. We call the space $U_m\coloneqq T_m\times I\subset T_m\times I^2$ the \textit{universal tree}. 

We now define the second map, which will be denoted by 
$$\overline{ev}_3\colon \pi_2 \mathrm{Emb}(I,\natural_mS^1\times D^3;I_0)\to \pi_5C_3'[S^1\times D^3]\otimes \mathbb{Q}.$$
This will make use of the mapping space model and involves the second and third stages of the embedding calculus tower. We start with analyzing the second stage. As discussed in the end of the previous section, the $k$-th element of the embedding calculus tower $T_k \mathrm{Emb}(I,\natural_mS^1\times D^3)$ can be modelled by stratum-preserving, aligned maps $C_{k+2}'[I,\partial]\to C_{k+2}'[\natural_mS^1\times D^3,\partial]$. For $[f]\in \pi_n\mathrm{Emb}(I,\natural_mS^1\times D^3;I_0)$ represented by a map 
$$f\colon S^n\to \mathrm{Emb}(I,\natural_mS^1\times D^3;I_0),$$
we have an induced aligned, stratum-preserving map:
$$ev_2(f)\colon S^n\times C_2'[I,\partial]\to C_4'[\natural_mS^1\times D^3,\partial] .$$
The space $C_2[I,\partial]$ is a truncated triangle (see Figure \ref{polytope}) with two vertices blown up into two extra edges denoted by $x=y=0$ and $x=y=1$. On these faces, the map $ev_k(f)$ gives a constant interval which is the image of $I_0$. Therefore, it is enough to ignore the fixed first and last coordinates and view this map as a stratum-preserving, aligned map $$ev_2(f)\colon S^{n}\times C_2[I]\to C_2'[\natural_mS^1\times D^3].$$ 
In practice, as Budney--Gabai indicated in \cite{Budney-gabai}, it is often helpful to just think of $C_{k+2}'[M,\partial]$ as $C_k'[M]$ but keep in mind that extra stratification information is encoded to make sure the spaces have the correct homotopy types.

For $n=1$, the restriction of $ev_2(f)$ to the three boundary facets of $C_2[I]$, which we denote by $x_1=0$, $x_2=1$ and $x_1=x_2$, by the aligned and stratum-preserving property, for each of these sub-strata, has codomain diffeomorphic to $C_1'[\natural_mS^1\times D^3]\cong(\natural_mS^1\times D^3)\times S^3$, and hence homotopy equivalent to $(\vee_m S^1)\times S^3$, and $$\pi_2 ((\vee_m S^1)\times S^3)=\pi_1\Omega((\vee_m S^1)\times S^3)=0.$$  
Therefore, these restrictions are null homotopic. By attaching null homotopies along these boundary facets, one obtains a map $$\overline{ev}_2(f)\colon S^3 \to C_2'[\natural_m S^1\times D^3],$$
i.e.~this is an element in $\pi_3 C_2'[S^1\times D^3]$. More precisely, these null homotopies can be constructed through the following Propositions \ref{constructnull} and \ref{constructnull2}, which are inspired by ideas in Ryan Budney's talk at Glasgow \cite{ryan}.

\begin{proposition}
\label{constructnull}
        Let $f\colon S^{n}\times C_2[I]\to C_2'[S^1\times D^3]$ be an aligned stratum-preserving map with $n\geq1$. The restriction of $f$ to the boundary facets of $C_2[I]$ is null homotopic.
\end{proposition}
\begin{proof}
    For $f\colon S^{n}\times C_2[I]\to C_2'[S^1\times D^3]$, we lift it to the universal cover $F\colon S^n\times C_2[I]\to \overline{C_2'[S^1\times D^3]}\subset C_2'[\overline{S^1\times D^3}]$. This map can be thought of as a map $(p,t_1,t_2)\to (q_1,q_2,v_1,v_2)$ where $p\in S^n$, and $v_1$ and $v_2$ are the corresponding velocity vectors. The points $q_1$ and $q_2$ are points in $\mathbb{R}\times D^3$. Along the $t_1=t_2$ facet, as $t_1$ approaches $t_2$, the velocity vectors agree with the direction of the collision, as given by $\mathrm{lim}_{t_1\to t_2}\frac{q_2-q_1}{t_2-t_1}$. The map $(q_1,q_2,\frac{q_2-q_1}{t_2-t_1},\frac{q_2-q_1}{t_2-t_1})$ defines an extension to the entire triangle, and the restrictions to the $t_1=0$ and $t_2=1$ facets can be straight-line homotoped to constant maps since they point to convex spaces (actually half spaces). Similarly, as $t_1$ approaches $0$, the corresponding vector is given by the formula $\mathrm{lim}_{t_1\to 0}\frac{f(t_2+t_1)-q_2}{t_1}$ which also defines an extension to the entire triangle. In addition, the restriction to the other two edges $t_1=t_2$ and $t_2=0$ can be homotoped to the constant map by straight-line homotopies. 

The null homotopy for the edges $t_1=0$ and $t_2=1$ can be constructed in a similar manner. We give details for the $t_1=0$ edge to illustrate this, and a null homotopy for the $t_2=1$ edge can be constructed using the same idea. The restriction of a map $f\colon S^n\times C_2[I]\to C_2'[S^1\times D^3]$ to the edge $t_1=0$ is given by
$$(p,0, t_2)\longmapsto (f_p(0),f_p(t_2),\mathrm{lim}_{t_1\to 0} \frac{f_p(t_1)-f_p(0)}{t_1},v_2).$$
As before, we define an extension of it to the remaining part of the triangle:
$$(p,t_1,t_2)\longmapsto (f_p(t_1),f_p(t_2),\frac{f_p(t_1)-f_p(0)}{t_1},v_2).$$
On the $t_1=t_2$ edge, this map is given by:
$$(p,t_2,t_2)\to  (f_p(t_2),f_p(r_2),\frac{f_p(t_1)-f_p(0)}{t_1},\mathrm{lim}_{t_1\to t_2}\frac{f_p(t_1)-f_p(t_2)}{t_1-t_2})$$
where the last coordinate is null homotopic by the null homotopy described in the last paragraph. Now, observe that the vector $\frac{f_p(t_1)-f_p(0)}{t_1}$ is confined to a hemisphere for all values of $t_1$, and hence it is null homotopic through the straight line homotopy.
\end{proof}

\begin{proposition}
\label{constructnull2}
For $f\colon S^n \times I\to \natural_m S^1\times D^3$ representing $[f]\in \pi_n\mathrm{Emb}(I,\natural_m S^1\times D^3;I_0)$ with $n\geq 1$, the restriction of the second evaluation
    $$ev_2(f)\colon S^n \times C_2[I]\to  C_2'[\natural_m S^1\times D^3]$$
    to the boundary facets of $C_2[I]$ is null homotopic.
\end{proposition}
\begin{proof}
   We use the same ideas from the proof of Proposition \ref{constructnull}. We describe a null-homotopy along the $x_1=x_2$ facet, and the other two facets are similar. As in the $S^1\times D^3$ case, along the $x_1=x_2$ facet, as $x_1$ approaches $x_2$, the velocity vectors agree with the direction of the collision. Along this facet of $C_2[I]$, we lift the restriction of $f$ to it to the universal cover $\overline{C_2'[\natural_m S^1\times D^3]}$ and observe that this map is given by the derivative on the facet $x_1=x_2$ (cf. discussion in the end of Section 4.1). There is an extension of this derivative to the entire triangle as follows. For $p\in S^2$ and $x_1,x_2\in C_2[I]$, we only need to specify one tangent vector. There is a unique geodesic connecting $q_1=f_p(x_1)$ and $q_2=f_p(x_2)$. This geodesic is given by a concatenation of intervals. Take the initial arc of the concatenation and define the direction vector (after scaling it a unit vector) determined by this arc to be our extension. In particular, such a map can be thought of as a map 
    \begin{align*}
        S^n\times C_2[I]\to C_2'[\natural_m S^1\times D^3]\\
        (p,x_1,x_2)\to (q_1,q_2,v)
    \end{align*}
    where $p\in S^n$, $q_i=f(x_i)$ for $i=1,2$, and $v$ is the tangent vector we just described. Now, the same argument applied in Proposition \ref{constructnull} works equally well here. Namely, if we choose some $I_0=\{pt\}\times I\times 0\subset T_m\times I^2$, then the restriction of this extension to the other two edges $x_1=0$ and $x_2=1$ can be homotoped to $I_0$ since these vectors point to one side of a 3-sphere (one of the ``vertical'' directions $I_0$).

    Similarly, one can construct null homotopies that can be attached along the other two facets $x_1=0$ and $x_2=1$ facets, just as in Proposition \ref{constructnull}. \end{proof}
The next lemma is a generalization of Proposition 3.2 of \cite{Budney-gabai}.
\begin{lemma}
    For $f\in \pi_1\mathrm{Emb}(I,\natural_m S^1\times D^3;I_0)$, $\overline{ev}_2(f)$ defines an element in $\pi_3 C_2'[\natural_m S^1\times D^3]\otimes \mathbb{Q}\cong\mathbb{Z}[t_1^{\pm},t_2^{\pm}\dots,t_m^{\pm}]\oplus \mathbb{Z}^2$ modulo the following relations:
    \begin{itemize}
    \item (0,1,0) for the $x_1=0$ facet
   \item $(t_i^0,1,1)$ for $i=1,\dots,m$ for the $x_1=x_2$ facet
    \item (0,0,1) for the $x_2=1$ facet,
\end{itemize}
where the first coordinate represents the Laurent polynomials part and the remaining two coordinates represent the $\mathbb{Z}^2$ part, which represents the degrees of the two velocity vectors.
\end{lemma}

\begin{proof}
These relations are obtained from the inclusions of the edges. The first relation comes from the inclusion map $C_1'[I]\cong \{x_1=0\}\times S^3 \to C_2'[I]$ of the $x_1=0$ edge sending $((0,x_2),1)\to (0,x_2,1,0)$. Therefore, the induced map 
$$\pi_3 C_1'[\natural_m S^1\times D^3]\cong \mathbb{Z}\to \pi_3 C_2'[\natural_m S^1\times D^3]\cong \mathbb{Z}[t_1^{\pm},t_2^{\pm}\dots,t_m^{\pm}]\oplus \mathbb{Z}^2$$ is given by $1\to (0,1,0)$. The third relation is obtained in exactly the same way, but with an induced map $1\to (0,0,1)$. The second relation comes from the inclusion of the edge $x_1=x_2$ inducing $((x_1,x_1),1)$ to $(x_1,t_i^0.x_1,1,1)$ for $i=1,\dots,m$.\end{proof}

Therefore, we have defined a homomorphism
$$\overline{ev}_2\colon \pi_1\mathrm{Emb}(I,\natural_m S^1\times D^3)\to \mathbb{Z}[t_1^{\pm},t_2^{\pm}\dots,t_m^{\pm}]/\langle 1 \rangle.$$
that will guide the way of our construction of $(W_3)_m$. Note that the generators $t_i$ generally do not commute with each other. We now move to consider the third evaluation map 
$$ev_3(f)\colon S^2\times C_3[I]\to C_3'[\natural_m S^1\times D^3]$$
for $[f]\in \pi_2\mathrm{Emb}(I,\natural_mS^1\times D^3;I_0)$. By Lemma \ref{constructnull2}, the restriction to the four facets of $C_3[I]$ can be homotoped to elements in $$\pi_4 C_3'[\natural_m S^1\times D^3]\cong \pi_4 (\vee_m S^1)\times (\vee_m S^1\vee S^3)\times (\vee_m S^1\vee S^3\vee S^3)\times (S^3)^3.$$
Since $\pi_4 S^3\cong\mathbb{Z}/2$, these restrictions are torsion. If the order of the restriction to a boundary facet is $o$, then the restriction of the map $ev_3(of)$ to the boundary facets is null homotopic, thus  we define 
$$\frac{1}{o}\overline{ev}_3(of)\coloneqq S^5 \to C_3'[\natural_m S^1\times D^3]. $$
by attaching null homotopies along the 4 boundary facets. It remains to argue that this is a rational homotopy invariant up to choices of null homotopies. This is done by exploring inclusions of the boundary facets
$$C_2[I]\to C_3[I]$$
and the induced maps 
$$\pi_5 C_2'[\natural_m S^1\times D^3]\to \pi_5 C_3'[\natural_m S^1\times D^3].$$
Before investigating these further, we set up some notations. We use the notation $\langle (t)_i\rangle$ to denote the free group with $m$ generators
    $$\langle (t)_1,(t)_2,\dots,(t)_m\rangle$$ and similarly the notation $\langle(t_1)_i\rangle$ and $\langle (t_2)_i\rangle$ for the corresponding free groups with $m$ generators. Also, $\nu((t)_1,\dots,(t)_m)$ denotes an element in $$\langle (t)_1,(t)_2,\dots,(t)_m\rangle$$ and similarly $\nu((t_1)_1,\dots,(t_1)_m)$ and $\nu((t_2)_1,\dots,(t_2)_m)$ denote elements in the groups $\langle(t_1)_i\rangle$ and $\langle (t_2)_i\rangle$ respectively.

The group $\pi_5 C_2'(\natural_m S^1\times D^3)$ is isomorphic to the following quotient group $$\mathbb{Z}[ (t_1)_i^{\pm},(t_2)_i^{\pm}]/\langle \nu((t_1)_i)=\nu^{-1}((t_2)_i),\forall \nu( (t)_1,(t)_2,\dots,(t)_m) \rangle$$mod torsion where $i,j=1,2,\dots,m$ and $(t_1)_i$ do not commute with $(t_2)_j$ for all $i,j$. A polynomial $(t_1)_i^{\alpha}(t_2)_j^{\beta}$ represents the Whitehead bracket $[(t_1)_i^{\alpha}(t_2)_j^{\beta}.\omega_{12},\omega_{12}]$. The following lemma is a generalization of Proposition 3.3 in \cite{Budney-gabai}.

\begin{lemma}
\label{verycomplicatedrelation}

    For each of the four facets $x_1=0$, $x_1=x_2$, $x_2=x_3$, $x_3=1$ of the tetrahedron $C_3[I]$, attaching a null homotopy (modulo torsion) $S^2\times C_2[I]\times I\to C_2[\natural_m S^1\times D^3]'$ gives rise to the following relations in 
    $$\pi_5 C_2'(\natural_m S^1\times D^3).$$
\begin{enumerate}
    \item For the $x_1=0$ facet, the generators $(t_1)_i^{\alpha}.\omega_{12}$ are mapped to $(t_2)_i^{\alpha}.\omega_{23}$ by the induced map on homotopy groups of the inclusion map, thus one obtains the relations $[(t_2)_i^{\alpha}.\omega_{23},(t_2)_j^{\beta}.\omega_{23}]$ for all $\alpha,\beta\in \mathbb{Z}$. More generally, the elements  $\nu((t_1)_i).\omega_{12}$ are mapped to $\nu((t_2)_i).\omega_{23}$, leading to the relations $$[\nu((t_2)_i).\omega_{23},\mu((t_2)_i).\omega_{23}]$$ for $\mu,\nu\in \langle (t_2)_i\rangle$.
    
    \item For the $x_1=x_2$ facet, the inclusion doubles the first coordinate $(x_1,x_2)\to (x_1,\delta x_1,x_2)$ where $\delta x_1$ is a small perturbation in the direction of the corresponding velocity vector. The induced map on homotopy groups maps $\omega_{12}$ to $\omega_{13}+\omega_{23}$, $(t_1)_i$ to $(t_1)_i (t_2)_i$ and fixes $(t_2)_i$, leading to the relations $[(t_1)_i^{\alpha}.\omega_{13}+(t_2)_i^{\alpha}.\omega_{23},(t_1)_j^{\beta}.\omega_{13}+(t_2)_j^{\beta}\omega_{23}]$. More generally, $\nu((t_1)_i)$ are mapped to $\nu((t_1 t_2)_i)$, leading to the relations $$[\nu((t_1t_2)_i)(\omega_{13}+\omega_{23}),\mu((t_1t_2)_i)(\omega_{13}+\omega_{23})]$$ for $\mu,\nu\in \langle (t)_i\rangle $.
    
    \item For the $x_2=x_3$ facet, the inclusion map doubles the second coordinate thus the induced map sends $\omega_{12}$ to $\omega_{12}+\omega_{13}$, $(t_2)_i$ to $(t_2)_i (t_3)_i$ and fixes $(t_1)_i$ leading to the relations $[(t_1)_i^{\alpha}.\omega_{12}+(t_1)_i^{\alpha}\omega_{13},(t_1)_j^{\beta}.\omega_{12}+(t_1)_j^{\beta}\omega_{13}]$. More generally, $\nu((t_1)_i)$ are mapped to $\nu((t_1)_i)$, leading to the relations $$[\nu((t_1)_i)(\omega_{13}+\omega_{12}),\mu((t_1)_i)(\omega_{13}+\omega_{12})]$$ for $\mu,\nu\in \langle (t_1)_i\rangle $.

    \item For the $x_3=1$ facet, the inclusion map fixes $x_1$ and $x_2$ and maps $x_3$ to $(1,0)$. Therefore, similarly to the $x_1=0$ facet, leading to the relations $[(t_1)_i^{\alpha}.\omega_{12},(t_1)_j^{\beta}.\omega_{12}]$. More generally, the relations $$[\nu((t_1)_i).\omega_{12},\mu((t_1)_i).\omega_{12}]$$ are satisfied for $\mu,\nu\in \langle (t_1)_i\rangle$.
\end{enumerate}
\end{lemma}
The relations 1 and 4 kill the relevant brackets that only involve $\omega_{23}$ and $\omega_{12}$ respectively. Thus the relations from 2 are simplified to \begin{align*}
   & [\nu((t_1t_2)_i).\omega_{13},\mu((t_1t_2)_i).\omega_{23}]+
    [\nu((t_1t_2)_i).\omega_{23},\mu((t_1t_2)_i).\omega_{13}]+
    [\nu((t_1t_2)_i).\omega_{13}, \mu((t_1t_2)_i).\omega_{13}]
\end{align*}and relations from 3 are simplified to \begin{align*}
    &[\nu((t_1)_i).\omega_{12},\mu((t_1)_i).\omega_{13}]+[\nu((t_1)_i).\omega_{13},\mu((t_1)_i).\omega_{12}]+[\nu((t_1)_i).\omega_{13}, \mu((t_1)_i).\omega_{13}].
\end{align*}
Now, observe that $$[\nu((t_1t_2)_i).\omega_{13}, \mu((t_1t_2)_i).\omega_{13}]=[\nu((t_1)_i).\omega_{13}, \mu((t_1)_i).\omega_{13}]$$ since $t_2$ acts trivially on $\omega_{13}$. Therefore, we can merge relations 2 and 3 into one relation in the following way. Relations 2 give rise to
\begin{align*}
   & [\nu((t_1t_2)_i).\omega_{13},\mu((t_1t_2)_i).\omega_{23}]+
    [\nu((t_1t_2)_i).\omega_{23},\mu((t_1t_2)_i).\omega_{13}]=\\
    &\nu ((t_1)_i)\mu^{-1}((t_1)_i)\mu^{-1}((t_3)_i).[\omega_{13},\omega_{23}]+
    \mu((t_1)_i)\nu^{-1}((t_1)_i)\nu ^{-1}((t_3)_i).[\omega_{23},\omega_{13}]=\\
     &(-\nu ((t_1)_i)\mu^{-1}((t_1)_i)\mu^{-1}((t_3)_i)+ \mu((t_1)_i)\nu^{-1}((t_1)_i)\nu^{-1} ((t_3)_i)).[\omega_{12},\omega_{23}],
\end{align*}
and relations 3 give rise to
\begin{align*}
      &[\nu((t_1)_i).\omega_{12},\mu((t_1)_i).\omega_{13}]+[\nu((t_1)_i).\omega_{13},\mu((t_1)_i).\omega_{12}]=\\
    &\nu((t_1)_i) \nu ((t_3)_i) \mu^{-1} ((t_3)_i).[\omega_{12},\omega_{13}]+
    \mu((t_1)_i)\mu((t_3)_i)\nu^{-1}((t_3)_i)[\omega_{13},\omega_{12}]=\\
    &(-\nu((t_1)_i) \nu ((t_3)_i) \mu^{-1} ((t_3)_i)+\mu((t_1)_i)\mu((t_3)_i)\nu^{-1}((t_3)_i))[\omega_{12},\omega_{23}].
\end{align*}
Therefore, we get the following relation as an analogue of Budney--Gabai's \textit{Hexagon relation} as discussed in the end of Section 4.2.
\begin{theorem}
 \label{newhexa}
The rational homotopy group $\pi_5 C_3'(\natural_m S^1\times D^3)\otimes \mathbb{Q}$ is generated by $(t_1)_i,(t_3)_i$, $i=1,2,\dots,m$ with the following relation satisfied.
    \begin{align*}
    &-\nu ((t_1)_i)\mu^{-1}((t_1)_i)\mu^{-1}((t_3)_i)+ \mu((t_1)_i)\nu^{-1}((t_1)_i)\nu^{-1} ((t_3)_i)=\\
    &-\nu((t_1)_i) \nu ((t_3)_i) \mu^{-1} ((t_3)_i)+\mu((t_1)_i)\mu((t_3)_i)\nu^{-1}((t_3)_i).
\end{align*}
After a change of variable and rearrangements of the equation, we can write the relation in the following way:
\begin{align*}
   & \nu((t_1)_i)\mu((t_3)_i)+\mu^{-1}((t_1)_i)\nu^{-1}((t_3)_i)=\\
   &\nu^{-1}((t_1)_i)\mu((t_3)_i)\nu^{-1}((t_3)_i)+\nu((t_1)_i)\mu^{-1}((t_1)_i)\nu((t_3)_i).
\end{align*}
\end{theorem}
We will denote this codomain by $\Lambda_m$. Therefore, we have now defined a map
$$(\overline{ev}_3)_m\colon \pi_2 \mathrm{Emb}(I,\natural_m S^1\times D^3) \to \Lambda_m = \pi_5 C_3'[\natural_m S^1\times D^3]\otimes \mathbb{Q}/R$$
where $R$ stands for the relations described in Lemmas \ref{verycomplicatedrelation} and \ref{newhexa}. We remark that if we set $m=1$ and choose $D_s$ to be a disk $\{pt\}\times D^3\subset S^1\times D^3$, then the discussion so far in this section coincide with the $W_3$ defined in \cite{Budney-gabai}, and recovers the Hexagon relation defined for $S^1\times D^3$ in \cite{Budney-gabai}:
 $$t_1^{\alpha-\beta}t_3^{-\beta}+t_1^{\beta}t_3^{\beta-\alpha}=t_1^\alpha t_3^{\alpha-\beta}+t_1^{\beta-\alpha}t_3^{-\alpha} $$
for $\alpha,\beta\in \mathbb{Z}$. Therefore, our definition is indeed a generalisation of Budney-Gabai's $W_3$ invariant.

\section{Computing $(W_3)_1$}

In this chapter, we compute $(W_3)_1$ for barbell diffeomorphisms of $ S^1\times D^3$, following the method of \cite{Budney-gabai2}. Budney--Gabai \cite{Budney-gabai,Budney-gabai2} proposed two different approaches in the calculation of $W_3$ of unknotted barbell diffeomorphisms of $S^1\times D^3$. The first method is to write the image of a barbell diffeomorphism under the scanning map as a linear combination of certain \textit{fundamental classes} they denote by $G(p,q)$ in $\pi_2\mathrm{Emb}(I,S^1\times D^3)$. Budney--Gabai (see Theorem 8.1 in \cite{Budney-gabai}) provided a formula to calculate $W_3$ of $\theta_k(v,w)$ (cf. Example \ref{firstexample}) using this approach. The second method is via intersection theory and the Pontryagin-Thom construction. This approach is only completed in \cite{Budney-gabai2} for the calculations for a special class of barbells denoted by $\delta_k$ (see Theorem 3.1 in \cite{Budney-gabai2}), but not for other unknotted barbells in $S^1\times D^3$. \begin{theorem}[Section 3 and Figure 10 of \cite{Budney-gabai2}]
\label{deltakkkk}
    The $W_3$ invariants of barbell diffeomorphisms induced from $$\delta_k\coloneqq\theta_k((0,\dots,1,0),(0,\dots,0,1))$$ 
    for $k\geq 4$ are nonzero and linearly independent. Therefore, the induced barbell diffeomorphisms $\Phi_{\delta_k}$ from $\delta_k$ are isotopically non-trivial and linearly independent for $k\geq 4$.
\end{theorem}

We aim to complete the calculations for unknotted barbells in $S^1\times D^3$ using this approach, and generalize the ideas to $\natural_m S^1\times D^3$ in Chapter 5. In Section 4.1, we discuss a generalization of the Pontryagin-Thom construction. Next, we introduce cohorizontal and collinear submanifolds and how those can help with computations in Section 4.2. We then compute $(W_3)_1$ for the barbells $\theta_k(v,w)$ introduced in Example \ref{firstexample} in Section 4.3, and presents infinite families of linearly independent barbell diffeomorphisms. Finally, we give a recipe for all unknotted barbells in $S^1\times D^3$ in Section 4.4.

\subsection{Pontryagin-Thom construction for a wedge of spheres}

In this section, we describe a generalization of the Pontryagin-Thom construction to the setting of oriented, framed cobordism groups and the stable homotopy groups of wedge product of spheres. Although it seems to be widely known, the author is not aware of a well-written source of this construction. This will be needed in later sections. Readers who are familiar with this can skip this section.

We first briefly recall the usual Pontryagin-Thom construction. Let $[M,N]$ denote the homotopy classes of maps between $M$ and $N$ and $[M,N]_*$ denote the homotopy classes of based maps. For $N$ simply-connected, e.g.~$N$ is a sphere, the forgetful map $[M,N]_{*}\to [M,N]$ is a bijection. Let  $\Omega_{m-q:M}^{fr}$ denote the group of framed cobordism classes of framed submanifolds in $M$ of codimension $q$. Given a manifold $M$ of dimension $m$, a map $f\colon M\to S^q$, and two regular values $p_0$ and $p_1$ in $S^q$, the submanifolds $f^{-1}(p_0)$ and $f^{-1}(p_1)$ are framed bordant by transporting the framings along a path. 

The \textit{Pontryagin-Thom collapse} map is defined as follows: given a framed submanifold $Y$ of codimension $q$ in $M$, we take a tubular neighbourhood $U\cong Y\times \mathbb{R}^q$ (which is diffeomorphic to some disk bundle) and construct a map from $M$ to $S^q$ that sends everything outside of $U$ to a basepoint, and inside $U$ we project to $\mathbb{R}^q$ and take the quotient of this disk bundle that sends the boundary sphere boundary to the basepoint too. This gives a map from $M$ to $S^q$ that is well defined up to framed cobordism.

\begin{theorem}[Pontryagin-Thom]
Given an $m$-manifold $M$, there is an isomorphism $[M,S^q]\to \Omega_{m-q:M}^{fr}$ given by taking the preimage of a regular value with the Pontryagin-Thom collapse as an inverse.
\end{theorem}


\begin{figure}
    \centering
    \includegraphics[width=0.4\textwidth]{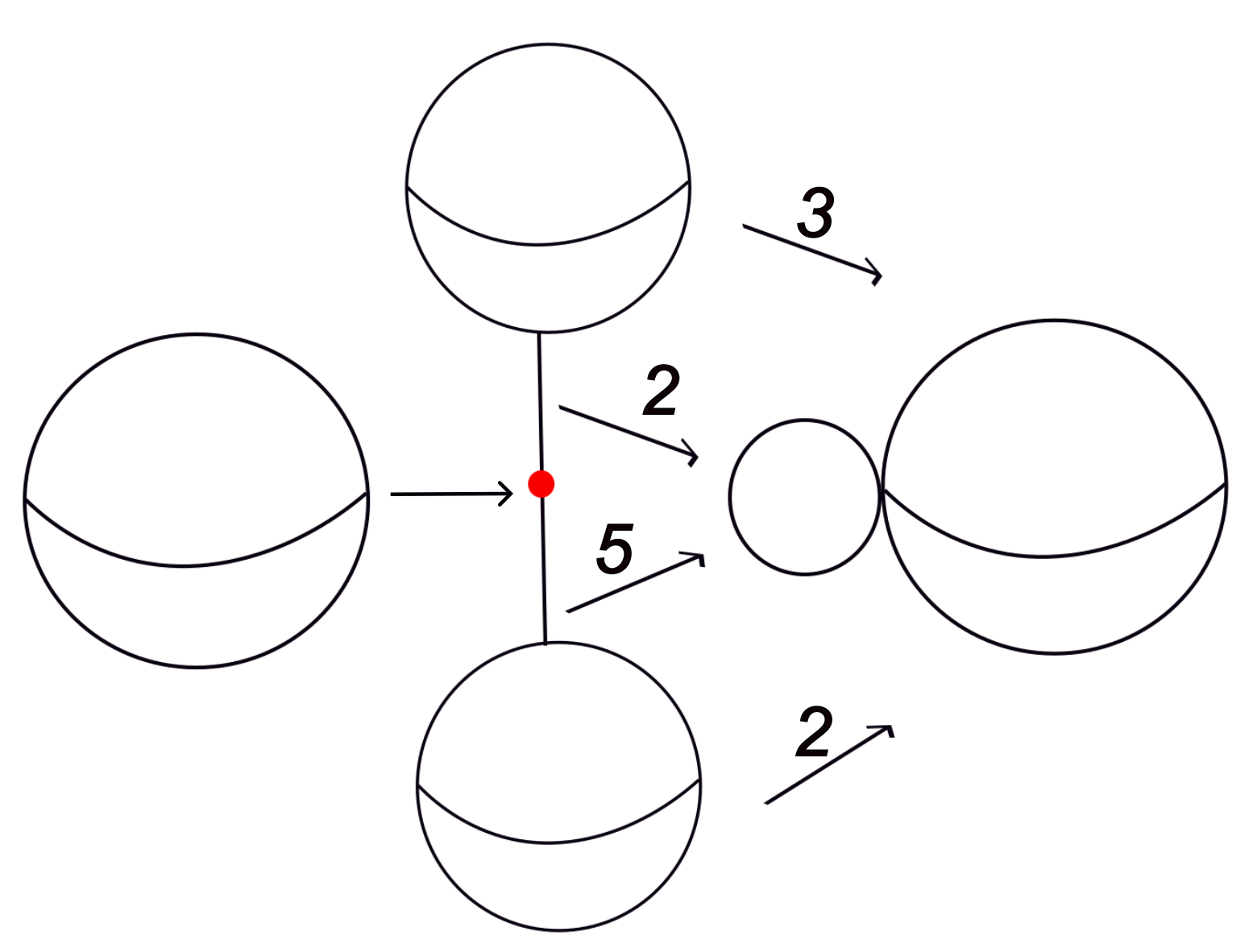}

    \caption{
       A representation of the element $2t^5+3t^2$.}
    \label{pic2}
\end{figure}

\begin{figure}
    \centering
    \includegraphics[width=0.4\textwidth]{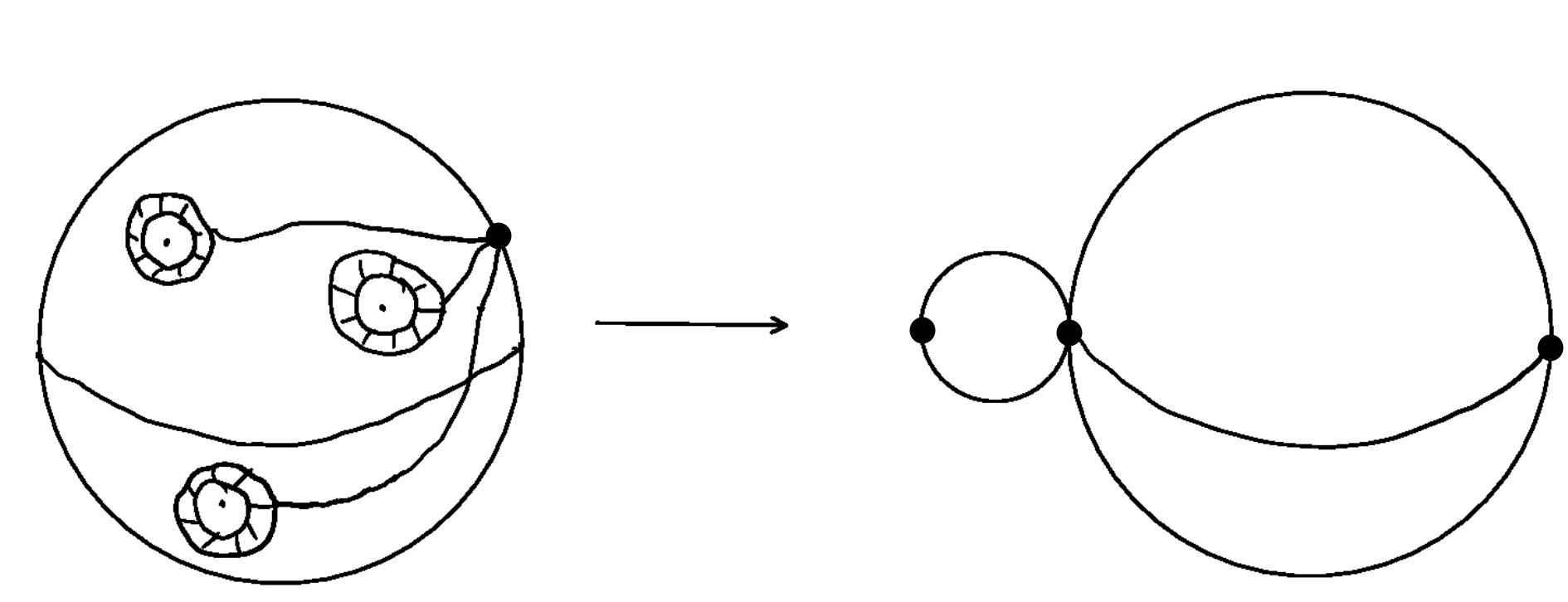}

    \caption{
        The backward map for $M=S^2$.}
    \label{pic}
\end{figure}
We now consider the case when $N=S^1\vee S^n$ and $M=S^m$.
\begin{theorem}
There is an isomorphism $$[S^m,S^1\vee S^n]\cong \Theta^{fr}_N.$$ Here $\Theta^{fr}_N$ is the set of finite collections of cobordism classes of disjoint framed codimension $n$ submanifolds in $S^1\vee S^n$. The map is given by taking the preimage of a regular value in $S^n\subset S^1\vee S^n$ and the preimage of a regular value in $S^1\subset S^1\vee S^n$ where the former gives rise to disjoint codimension $n$ submanifolds in $S^m$ and the latter gives rise to disjoint codimension 1 spheres in $S^m$ separating the codimension-$n$ submanifolds. Furthermore, each codimension $n$ submanifold is assigned an integer $a_k$, called the degree which is well defined up to an overall additive constant, and is given by the signed number of codimension 1 spheres that separate it from the chosen basepoint.
\end{theorem}
\begin{proof}
Let the wedge point of $S^1\vee S^n$ be $y_0$ which we fix as our basepoint. Given a map $f\colon S^m\to S^1\vee S^n$, we take a point $y_1$ in the circle $S^1\subset S^1\vee S^n$ that is not $y_0$. The preimage $f^{-1}(y_1)$ is a union of finitely many codimension-1 (closed) oriented submanifolds in $S^m\setminus y_0$ which we denote by $M_k$, $k=1,\dots,K$ so $f^{-1}(S^1\vee S^n\setminus y_1)$ is a copy of $S^m$ with a finite number of codimension-1 oriented submanifolds removed. Note that since we can locally homotope $f$ such that the restriction to $f^{-1}(S^1\subset S^1\vee S^n)$ becomes smooth, then any point in $S^1$ other than $y_0$ is a regular point by dimension reasons. In addition, take any regular point $y$ in $S^n\subset S^1\vee S^n$. The preimage is a disjoint union of codimension-$n$ submanifolds $N_l$ in $S^m\setminus f^{-1}(y_0)$, $l=1,\dots,L$.

The neighbourhood of each $M_k$ gets mapped in the following way: take the disk bundle of its framing $M_k\times D^1\subset M_k\times \mathbb{R}$ and project to the $\mathbb{R}$ factor, and take the quotient that sends the boundary of $D^1$ to the basepoint to get a map between circles. The manifolds $N_l$ are separated by $M_k$ in the following way: collapsing the neighbourhoods $M_k\times D^1$ into $\{\mathrm{pt}\} \times D^1$ gives a quotient of $S^m$ as a wedge of $L$ long strings at $y_0$ with the other end of each long string connected with an $m$-dimensional sphere, and each long string consists a sequence of at least 1 arcs such that the total number of arcs $K$ is the total number of $M_k$, and $L$ is the total number of $N_l$. Furthermore, each $N_l$ is contained in an $S^m$ in this space. We record the number of arcs each long string contains, leading to degrees corresponding to $N_l$'s. See Figure \ref{pic2} for an example when $m=2$ (the red dot in denotes a choice of basepoint). Note that all of the circles and spheres are oriented, and the orientations are recorded by the numbers in the figure. This gives the forward direction. Up to homotopies of $f$, different choices of regular values lead to framed cobordisms between the corresponding preimages by choosing paths between these choices, and homotopies between two maps gives framed cobordisms.

Conversely, given a collection of disjoint representatives $M_k$ of codimension $n$ framed submanifolds of $S^m$, each of which is assigned a degree $a_k$, we can choose a (closed) tubular neighbourhood $U_k\cong M_k\times D^n$ for each of them and choose a bigger tubular $U_k'$ neighbourhood for each such that they are still disjoint. Then define a map $f'\colon S^m\to S^1\vee S^n$ as follows. In the tubular neighbourhoods $M_k\times D^n$, $f'$ first projects onto the $D^n$ part and then sends $0\in D^n$ to the antipodal point $y_0'$ of $y_0$ in $S^n\subset S^1\vee S^n$. The remaining part of $\mathrm{Int}D^n$ (which is an annulus) is then mapped to $S^n\setminus\{y_0,y_0'\}$ by the identity map between the annuli. It follows that $\partial D^n$ is mapped to $y_0$. Furthermore, the remaining part of the larger tubular neighbourhoods $U_k'\setminus U_k\cong M_k \times S^{n-1}\times I $ is projected to $I$ and then the end points of $I$ are assigned to $y_0$ with $\mathrm{Int} I$ wrapping around the circle part $a_k$ times. Finally, the rest of $S^n$ is sent to $y_0$. This gives the backward direction. One checks that they are inverses of each other. 
\end{proof}
\begin{remark}
    Alternatively, one can lift a map $f\colon S^m\to S^1\vee S^n$ to the universal cover of $S^1\vee S^n$ which is $\mathbb{R}\vee_\infty S^n$ and perform the above computations by taking preimages of points in the various lifted $S^n$.
\end{remark}

In this case, we have $\pi_n(S^1\vee S^n)\cong \mathbb{Z}[t,t^{-1}]$. This can be seen by passing to the universal cover. One way of representing the generators is as follows: given a polynomial $h=a_0+a_1t^{k_1}+\dots+ a_m t^{k_m}$, squash $S^n$ into a sequence of $m+1$ copies of $S^n$ (labeled by $S^n_{i}$) connected by $m$ copies of arcs (labeled by $I_{k_i}$), and map $S^n_{i}$$\subset S^1\vee S^n$ to $S^n$ by a degree $a_i$ map. Further, the arcs $I_{k_i}$ are mapped to $S^1\subset S^1\vee S^n$ by a degree $k_i$ map (with endpoints being sent to $y_0$).

Figure \ref{pic2} and \ref{pic} give an example when $m=n=2$. In this case, everything can be easily visualised. For $n=2$, for such a map $f$, the preimage $f^{-1}(S^1\vee S^n\setminus y_1)$ is a copy of $S^2$ with $n$ punctures. The preimage of any regular value in $S^2$ gives $m$ points in the domain, each of which is separated from the rest by a circle. Conversely, given a collection of points in $S^2$ and corresponding indices, we take a disk around each of the points and an annulus outside the disk. The resulting map is now easy to describe: the boundary of disks are quotiented to the wedge point of $S^1\vee S^2$ and the annuli part (the dashed part in Figure \ref{pic}) is collapsed into an interval with endpoints being sent to $y_0$ and the interior naturally circles around $S^1$. Note that the paths in the domain $S^2$ from the basepoint to each tubular neighbourhood are not necessary but just auxiliary: everything outside of the tubular neighbourhoods get mapped to the wedge point $y_0$.

The above can be generalized to the case of $(\vee_k S^1)\vee (\vee_l S^n)$ by picking up a regular point in each of the spheres (and circles) away from the wedge point and taking a preimage to obtain collections of disjoint submanifolds up to the sloping point. In the following sections, we will make use of the case when $n=3$, $k=1,2,3$ and $l=1,2$.

\subsection{Collinear submanifolds and intersection theory}

In this section, we discuss ideas for calculating the map
$$\overline{ev}_3\colon \pi_2 \mathrm{Emb}(I,S^1\times D^3;I_0)\to \Lambda.$$
For an element $[f]\in \pi_2 \mathrm{Emb}(I,S^1\times D^3)$, we start by analyzing (the homotopy class represented by) the map  $$\overline{ev}_2(f)\colon D^1\times D^1 \times C_2[I]\to C_2[S^1\times D^3]\simeq S^1\times (S^1\vee S^3).$$
Budney--Gabai proposed a method of constructing cobordism classes of certain families of disjoint, framed submanifolds that detect such maps as we shall disclose now.

\begin{figure}
  \centering
  \begin{minipage}[b]{0.4\textwidth}
    \includegraphics[scale=0.4]{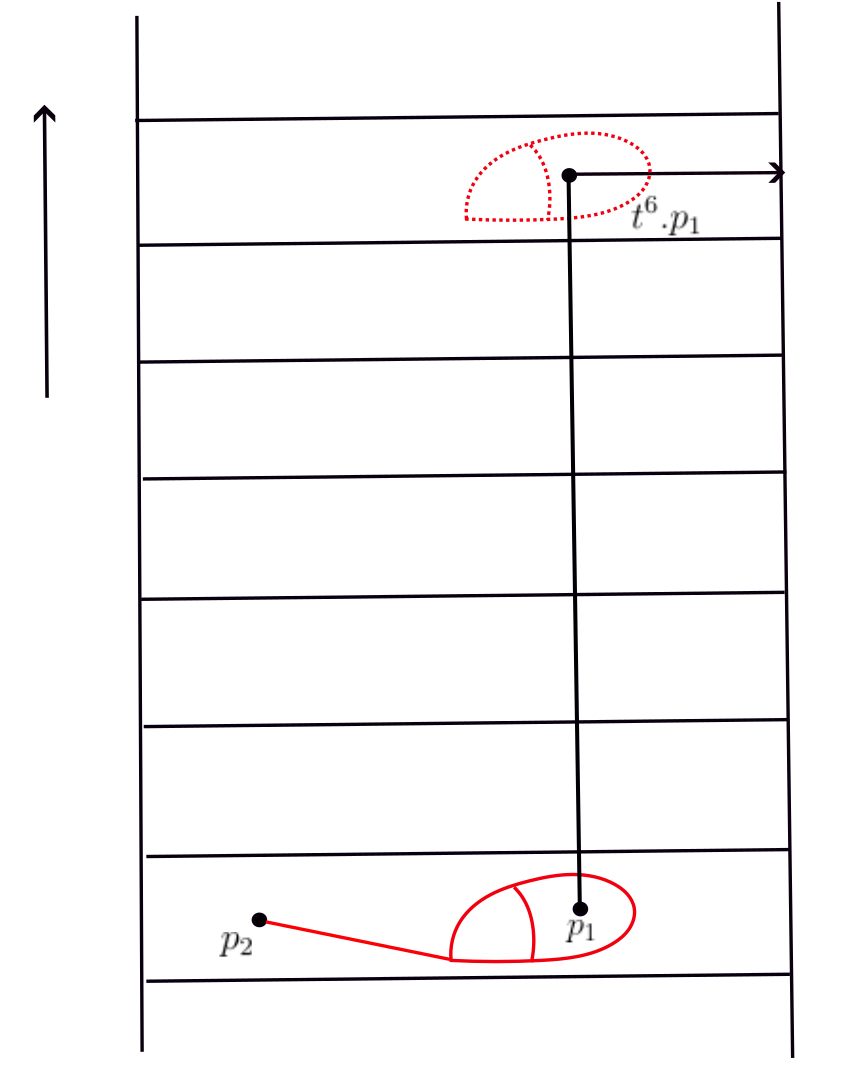}
    \caption{The cohorizontal submanifold $t^6\mathrm{Co}_1^2$ in $\mathbb{R}\times D^3$ detects $t^6.\omega_{12}$.}
    \label{cohoriz}
  \end{minipage}
  \hfill
  \begin{minipage}[b]{0.4\textwidth}
    \includegraphics[scale=0.35]{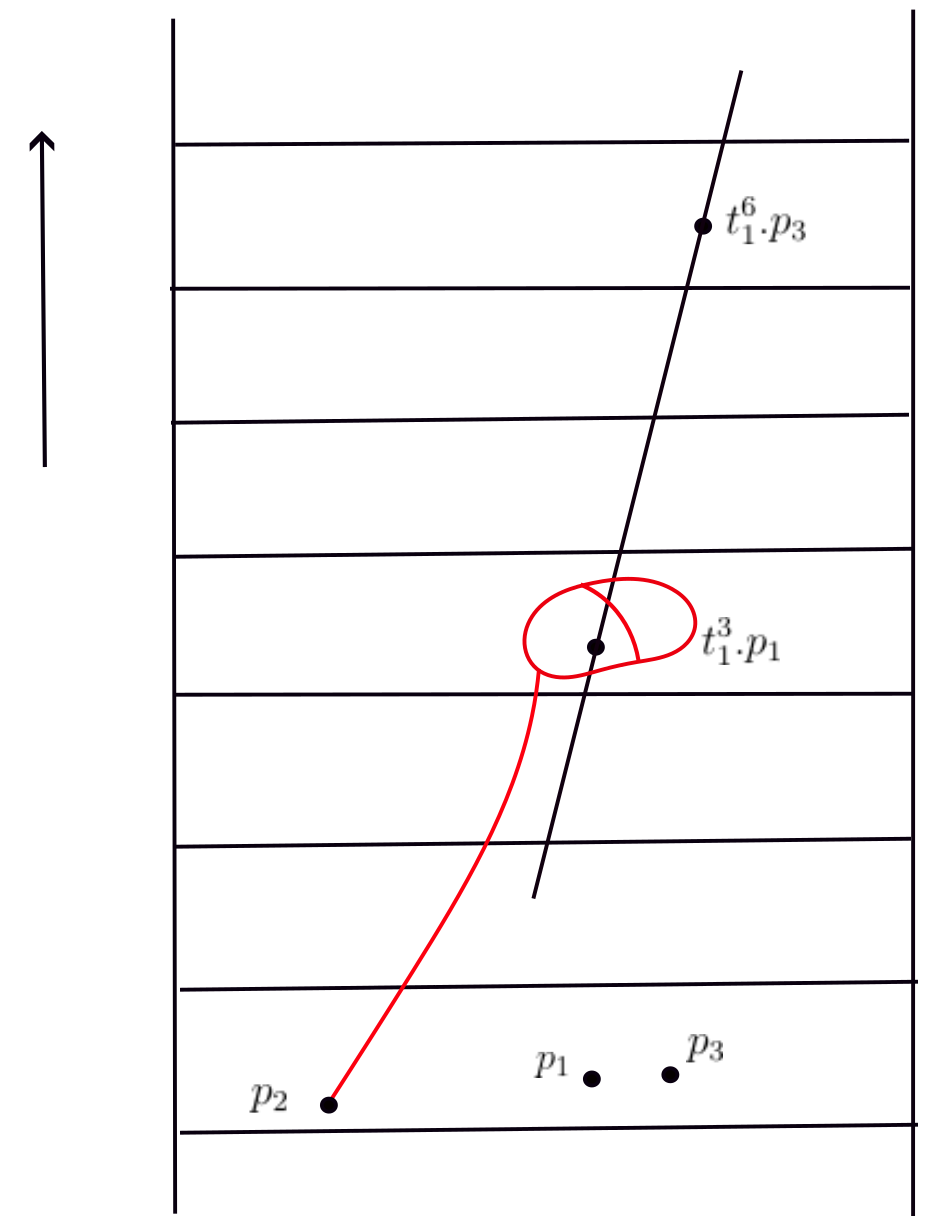}
    \caption{The collinear submanifold $\mathrm{Col}_{3,6}^1$ in $\mathbb{R}\times D^3$ detects $t^3.\omega_{12}$. }
    \label{collinearexample}
  \end{minipage}
\end{figure}

\begin{definition}
    
Let $\delta\in \partial D^3$ be a fixed unit direction. For $i\in \mathbb{Z}$, we define the \textit{cohorizontal submanifold} 
$$t^i\mathrm{Co}_1^2=\{(p_1,p_2)\in C_2[\mathbb{R}\times D^3]\colon t^i.p_2-p_1=\lambda \delta \text{ for some } \lambda>0\}.$$

\end{definition}
This submanifold intersects the image of  $t_1^i.\omega_{12}\colon S^3\to C_2[S^1\times D^3]$ transversely at exactly one point; namely, suppose that we pick a representative of the lift $\overline{t_1^i.\omega_{12}}$ of $t_1^i.\omega_{12}$ to $C_2[\mathbb{R}\times D^3]$ by shooting the first point from the basepoint say $(0,0)$ to $ (i,0)\in \mathbb{R}\times D^3$, and letting the second point loop around $(0,0)$, then the unique intersection point between this image and $t^i\mathrm{Co}_1^2$ is obtained by taking the halfline to the $\delta$-direction that starts from $(i,0)$. See Figure \ref{cohoriz} for a cartoon picture.  The horizontal rectangles in Figure \ref{cohoriz} represent fundamental domains, and the solid red sphere (drawn as a 2-sphere but aiming to represent a 3-sphere) loops around the point $p_1$ and creates a cohorizontal pair. The dotted red sphere represents the image of the second point under the action of $t^6$. The horizontal arrow represents the cohorizontal direction $\delta$. Therefore, we have the following.

\begin{lemma}
\label{cohorizontaldetect}
    The preimage of $t^i\mathrm{Co}_1^2$ under the map $\overline{t_1^i.\omega_{12}}\colon S^3\to C_2[\mathbb{R}\times D^3]$ is a codimension-3 oriented submanifold, which in this case is a single point. For $j\neq i$, the preimage of $t^i\mathrm{Co}_1^2$ under $\overline{t_1^j.\omega_{12}}$ is empty. Therefore, we say that $t^i\mathrm{Co}_1^2$ detects $t^i.\omega_{12}$.
\end{lemma}
Thus, we can use cohorizontal submanifolds to detect the map $f$. In particular, we lift it to $C_2[\mathbb{R}\times D^3]$ and take the preimage of the cohorizontal submanifolds $t^i\mathrm{Co}_1^2$ under this lift. The preimage is given by a disjoint union of (signed) points for all $i$. Counting the signed number of points gives the coefficient of $t^i\mathrm{Co}_1^2$, and taking the sum gives $[f]$.

Next, we consider the third stage map $\overline{ev}_3(f)\in \Lambda$. Such a class is a linear combination of Whitehead brackets $[t_2^{\alpha}.\omega_{12},t_2^{\beta}.\omega_{23}]=t_1^{\alpha}t_3^{\beta} [\omega_{12},\omega_{23}]$ with coefficients detected by the following two collections of \textit{collinear submanifolds}:
$$\mathrm{Col}_{\alpha,\beta}^1=\{(p_1,p_2,p_3)\in \overline{C_3[S^1\times D^3]}\colon (p_2,t^{\alpha}.p_1,t^{\beta}.p_3) \text{ lie on a straight line in $\mathbb{R}\times D^3$ in this order}\}$$
$$\mathrm{Col}_{\alpha,\beta}^3=\{(p_1,p_2,p_3)\in \overline{C_3[S^1\times D^3]}\colon (t^{\alpha}.p_1,t^{\beta}.p_3,p_2) \text{ lie on a straight line in $\mathbb{R}\times D^3$ in this order}\}.$$
See Figure \ref{collinearexample} for a picture of $\mathrm{Col}_{\alpha,\beta}^1$ with $\alpha=3$ and $\beta=6$.

The universal cover $\overline{C_3[S^1\times D^3]}$ can be interpreted as a subspace of $C_3[\mathbb{R}\times D^3]$ in the sense that it contains points of the latter space with disjoint $\mathbb{Z}$-orbits (also see Figure 9 of \cite{Budney-gabai} for details). The next lemma summarises the the situation.

\begin{lemma}[\cite{Budney-gabai2}]
\label{collineardetect}
The two collections of collinear submanifolds $\mathrm{Col}_{\alpha,\beta}^1$ and $\mathrm{Col}_{\alpha,\beta}^3$ are disjoint. If we consider $t_2^{\alpha}.\omega_{12}$ and $t_2^{\beta}.\omega_{23}$ as elements of $\pi_3C_3[S^1\times D^3]$ by taking their images under the inclusion maps $C_2[S^1\times D^3]\to C_3[S^1\times D^3]$ given by
$$\mathrm{Inc}_{12}\colon (p_1,p_2)\to (p_1,p_2,p_0)$$
$$\mathrm{Inc}_{23}\colon (p_2,p_3)\to (p_0,p_2,p_3)$$
for some $p_0\in S^1\times D^3$ which we can choose to be distinct to the other coordinates, then the former collinear submanifold $\mathrm{Col}_{\alpha,\beta}^1$ detects $t_2^{\alpha}.\omega_{12}$ and the latter collinear submanifold $\mathrm{Col}_{\alpha,\beta}^3$ detects $t_2^{\beta}.\omega_{23}$ in the sense that the (transverse) intersections between $\mathrm{Col}_{\alpha,\beta}^1$ and $t_2^{\alpha}.\omega_{12}$ (and $\mathrm{Col}_{\alpha,\beta}^3$ and $t_2^{\beta}.\omega_{23}$) are exactly one point, and for $\alpha'\neq \alpha$ and $\beta'\neq \beta$, the (transverse) intersections between $\mathrm{Col}_{\alpha,\beta}^1$ and $t_2^{\alpha'}.\omega_{12}$ (and $\mathrm{Col}_{\alpha,\beta}^3$ and $t_2^{\beta'}.\omega_{23}$) are zero.

In addition, the preimage of the pair $(\mathrm{Col}_{\alpha,\beta}^1,\mathrm{Col}_{\alpha,\beta}^3)$ under (the lift of) $[t_2^{\alpha}.\omega_{12},t_2^{\beta}.\omega_{23}]$ is a 2-dimensional Hopf link with linking number 1 in $S^5$. Therefore, we say that the pair $(\mathrm{Col}_{\alpha,\beta}^1,\mathrm{Col}_{\alpha,\beta}^3)$ detects the Whitehead product $[t_2^{\alpha}.\omega_{12},t_2^{\beta}.\omega_{23}]$.
\end{lemma}

\begin{proof}
   The property of being disjoint follows directly from the orderings of the triples. We show that the former collinear submanifold $\mathrm{Col}_{\alpha,\beta}^1$ detects $t_2^{\alpha}.\omega_{12}$. The latter one can be argued in a similar way. 

By choosing an appropriate representative of $t_2^{\alpha}.\omega_{12}$, we can arrange its image to be $p_1\times (\gamma\vee S)\times p_0$ where $p_1$ is a fixed point in $0\times D^3\subset \mathbb{R}\times D^3$, and $\gamma$ is a null homotopic path from a fixed point $p_2\in 0\times D^3 \subset \mathbb{R}\times  D^3$ to a 3-sphere $S$ that loops around $t_1^{\alpha}.p_1\in \alpha\times D^3 \subset \mathbb{R}\times D^3$. Again, see Figure \ref{collinearexample} for a picture of $\mathrm{Col}_{\alpha,\beta}^1$ with $\alpha=3$ and $\beta=6$. Now, we can deduce that the intersection between this image and the submanifold $\mathrm{Col}_{\alpha,\beta}^1$ consists of a unique point. Namely, $t_3^\beta. p_0$ and $t_1^\alpha.p_1$ determines a straight line which intersects $S$ at a unique point. 

To see the last statement, we recall that $C_3[S^1\times D^3]$ is homotopy equivalent to $S^1\times (S^1\vee S^3_{12})\times (S^1\vee S_{12}^3 \vee S_{23}^3)$. The elements $\omega_{12}$ and $\omega_{23}$ represent the homotopy classes of the inclusion maps of the 3-spheres $S_{12}^3$ and $S_{23}^3$ respectively. The actions of $t_2^{\alpha}$ and $t_2^{\beta}$ reflect the actions of the second and third circles in the wedge. Fixing a basepoint appropriately, we can take the Whitehead product $[t_2^{\alpha}.\omega_{12},t_2^{\beta}.\omega_{23}]$ by composing the wedge map with the attaching map of the top cell of $S^3\times S^3$:
$$S^5\to S^3\vee S^3 \xrightarrow{t_2^{\alpha}.\omega_{12}\vee t_2^{\beta}.\omega_{23}} C_3[S^1\times D^3].$$
The preimage of the pair $(\mathrm{Col}_{\alpha,\beta}^1,\mathrm{Col}_{\alpha,\beta}^3)$ under this composition is the same as the preimage of two points: say $x_1$ in $S_{12}^3$ and $x_2$ in $S_{23}^3$ that are distinct from the wedge points. Since the attaching map 
$$S^5\cong \partial D^6\cong \partial (D^3\times D^3)\cong \partial D^3 \times D^3 \cup D^3 \times \partial D^3\to S^3\vee S^3$$is given by $(\psi_1 \times \overline{\psi_2}\cup \overline{\psi_2}\times \psi_1)$, where $\psi_i$ is the attaching map of the 3-cell the $i$-th 3-sphere and $\overline{\psi_i}$ is the corresponding characteristic map, one concludes that the preimage of the two points $x_1$ and $x_2$ is a disjoint union of two 2-spheres with linking number 1.
\end{proof}

The preimage of the pair $(\mathrm{Col}_{\alpha,\beta}^1,\mathrm{Col}_{\alpha,\beta}^3)$ under (the lift of) $\overline{ev}_3(f)$ is a disjoint union of pairs of codimension 3 (thus dimension 2) oriented submanifolds in $S^5$, and the linking numbers of these pairs determines the coefficient of $[t_2^{\alpha}.\omega_{12},t_2^{\beta}.\omega_{23}]$ in the linear combination representing $\overline{ev}_3(f)$.

In practice, however, sorting out collinear submanifolds (and their preimages) in a 5-dimensional space like $I\times I\times C_3[I]$ is not a straightforward task. The following lemma helps us simplify calculations involving collinear submanifolds by reducing the above linking number computation strategy to linking number computations that only involve cohorizontal submanifolds which are much easier to manage and visualize.
\begin{lemma}[\cite{Budney-gabai}]
\label{cohorizontal calculation}
    Given a smooth map $f\colon S^5\to C_3[S^1\times D^3]$, we can assume generically that it has no cohorizontal triples. As before, we fix $\delta\in \partial D^3$. Define the cohorizontal manifold for $k,j\in \{1,2,3\}$
    $$t^i \mathrm{Co}_j^k=\{(p_1,p_2,p_3)\in C_3[S^1\times D^3]\colon t^i.p_j -p_k=\lambda \delta \text{ for some } \lambda>0\}.$$
In other words, we allow one of the three coordinates to be free with the other two being cohorizontal. Then the linking numbers of the preimage of the pair $(\mathrm{Col}_{\alpha,\beta}^1,\mathrm{Col}_{\alpha,\beta}^3)$ under the lift of $f$ to the universal cover $\overline{f}\colon S^5\to \overline{C_3[S^1\times D^3]}$ are the same as the linking numbers of the preimage of the pair $(t^{\alpha}\text{Co}^1_2-t^{\alpha-\beta}\text{Co}^1_3, t^{\beta-\alpha}\text{Co}^3_1-t^\beta\text{Co}^3_2)$ under $f$ for all $\alpha,\beta\in \mathbb{Z}$.
\end{lemma}

This lemma allows one to perform the calculations at the third stage entirely by using the idea of cohorizontal submanifolds rather than running into the difficulties of visualizing collinear triples. We will present some ideas behind the proof of this lemma in Section 5.2 in a more general setting.

\subsection{Computing $W_3$ for $\theta_k(e_i,e_j)$}

\begin{figure}
    \centering
    \includegraphics[width=0.55\textwidth]{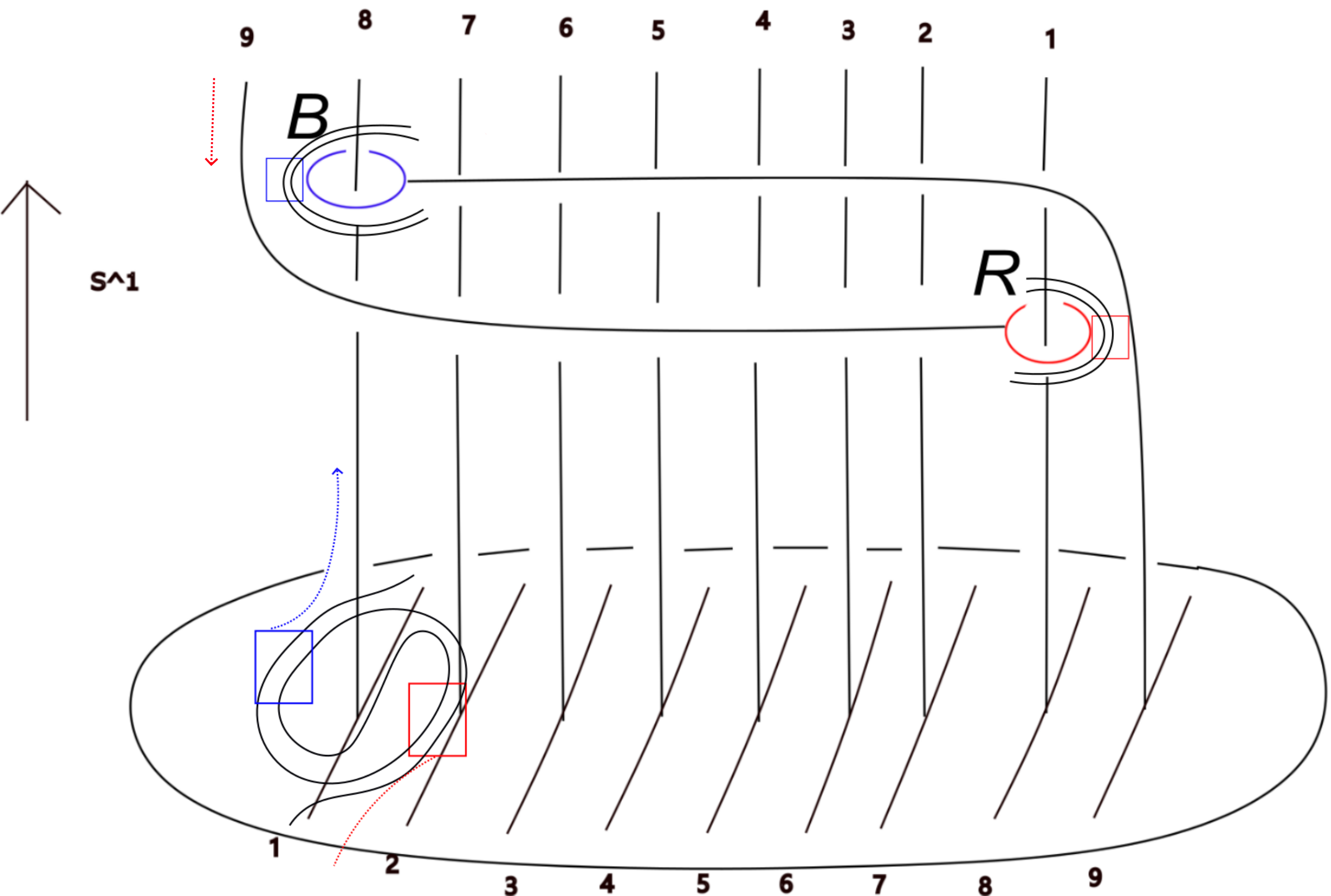}
    \caption{The embedded barbell $\delta'_{10}=\theta_{10}((0,\dots,0,1,0),(0,\dots,0,1,0))$.}
    \label{barbelltheta92}
\end{figure}

Recall from Chapter 2 (cf. discussion around Example \ref{firstexample}) that unknotted barbells in $S^1\times D^3$ are determined by the double coset $$\langle \nu_B \rangle\backslash\langle \nu_B, \nu_R,t\rangle/\langle \nu_R \rangle.$$ Consider the special classes of barbells studied by Budney--Gabai denoted $\theta_k(v,w)$, with $k\in \mathbb{Z}^{+}$, $(v)_i,(w)_j\in \mathbb{Z}^{k-1}$ from Example \ref{firstexample} (refer to Figures \ref{barbelltheta9} and \ref{barbelltheta92} for examples). Let $\delta'_k$ denote the barbell in $S^1\times D^3$ whose defining vectors are $v= w=(0,0,\dots,0,1,0)$, i.e.~
$$\delta'_k\coloneqq \theta_k((0,0,\dots,0,1,0),(0,0,\dots,0,1,0))=\theta_k(e_{k-2},e_{k-2}).$$
See Figure \ref{barbelltheta92} (which is a repetition of Figure \ref{barbelltheta9} for the case $k=10$, with more annotations which we will explain below). We explore the computation of $W_3$ for this barbell as an important example before diving into more general statements.

\begin{figure}
\centering
       \includegraphics[width=0.4\textwidth]{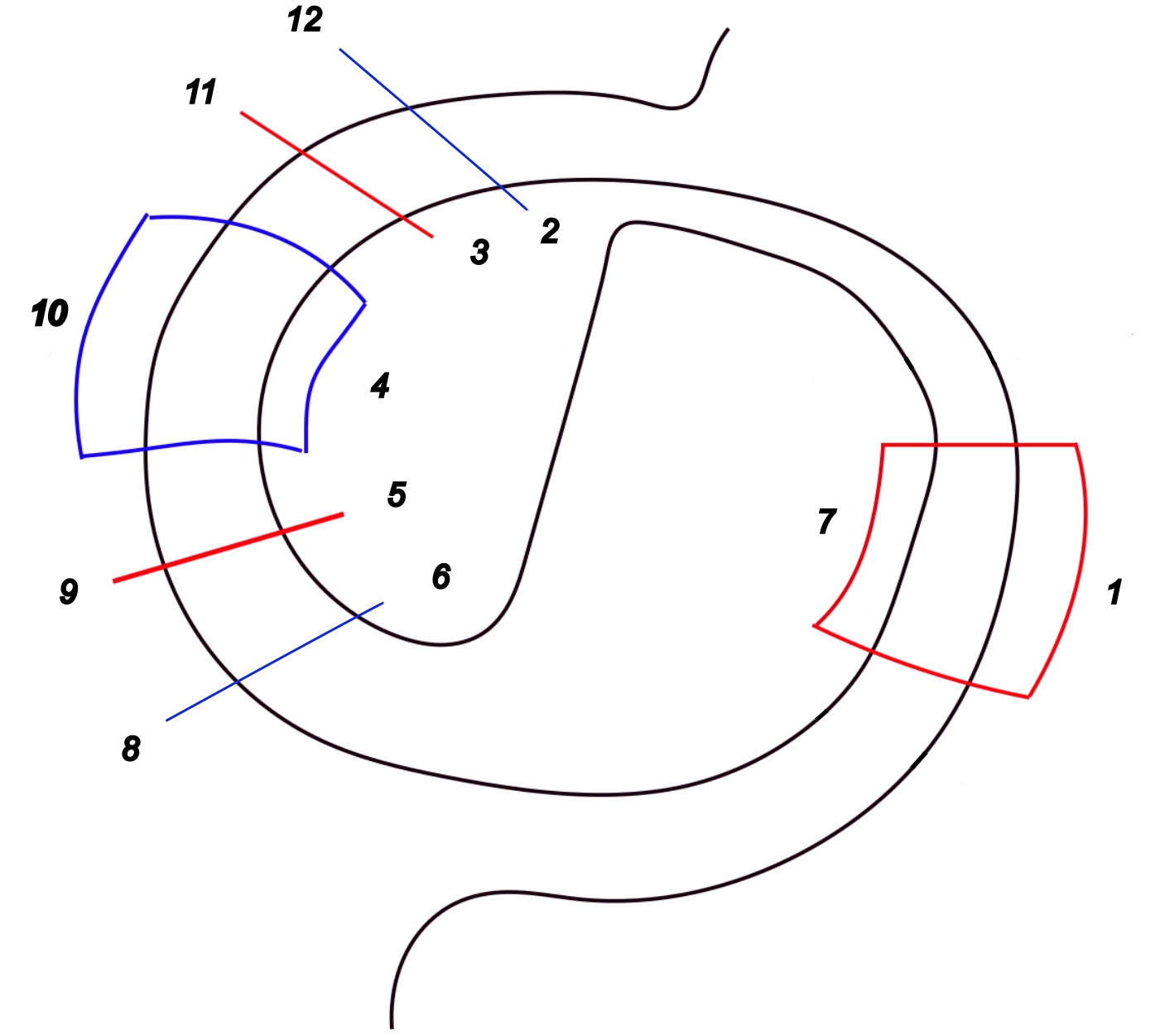}
        \caption{The Type 1 intersection for point 1.}
         \label{twisted arc}
\end{figure}

Recall that to define the scanning map $s_{s_0}$, we choose a fixed point $s_0\in S^1$ and take the disk $s_0\times D^3\subset S^1\times D^3$. For a class $[\Phi]\in \pi_0\mathrm{Diff}(S^1\times D^3,\partial)$, the image $s_{s_0}([\Phi])$ can be analyzed by a case-by-case method by choosing a suitable representative $\Phi$ and looking at $\Phi(s_0\times D^3)$ in detail. The first step is analyzing the image $s_{s_0}([\Phi_{\delta_k'}])$ under the scanning map $s_{s_0}$, where $\Phi_{\delta_k'}$ is the induced barbell diffeomorphism from $\delta_k'$. By assumption, $\delta_k'$ intersects the 3-ball $ \{s_0\} \times D^3$ transversely at the bar at $k-1$ isolated points but not at the cuff spheres. At each of the $k-1$ intersection points on the bar, following the approach in \cite{Budney-gabai,Budney-gabai2} (in particular, see Proposition 6.3 of \cite{Budney-gabai} and Proposition 2.2 of \cite{Budney-gabai2}), we twist the corresponding embedded $k-1$ arcs (drawn as short arcs in the horizontal plane in Figure \ref{barbelltheta92}) controlled by the two parameter families (blue and red) artificially into the shape shown in Figure \ref{twisted arc} (ignoring the colored arcs on sides of the colored boxes for now). Virtually, there are five strands in Figure \ref{twisted arc}, two in the blue box, two in the red box, and one in the middle. The following lemma states an alternative way of interpreting the image for scanning barbell diffeomorphisms.

\begin{lemma}[\cite{Budney-gabai,Budney-gabai2}]
\label{alternativebarbell}
    The barbell diffeomorphism $\Phi_{\delta_k'}$ is isotopic to a diffeomorphism whose restriction to the 3-ball $ \{s_0\} \times D^3$ has the following effects on each twisted arc (as in Figure \ref{twisted arc}) from each of the $k-1$ intersection points: it grabs the two strands in the blue box on the left side, following the bar (as the dotted blue arrow indicated in Figure \ref{barbelltheta92} for the first interval) until getting to the blue cuff, loops around the blue cuff sphere and comes back, and at the same time, grabs the two strands in the red box on the right side, following the bar  (as the dotted red arrow indicated in Figure \ref{barbelltheta92} for the first interval) until getting to the red cuff, loops around the red cuff sphere and comes back. Therefore, the image $s_{s_0}([\Phi_{\delta_k'}])$ is given by a sum denoted by 
    $$[\delta_k']_1+[\delta_k']_2+\dots +[\delta_k']_{k-1}$$
    of $k-1$ 2-parameter families of embedded intervals, each of which is given by twisting one of the $k-1$ arcs to the shape of Figure \ref{twisted arc}, following the above description, and finally untwisting.
\end{lemma}

By Lemma \ref{alternativebarbell}, the second stage map can be written as a sum
$$\overline{ev}_2(s_{s_0}([\Phi_{\delta_k'}]))=\overline{ev}_2(s_{s_0}([\Phi_{\delta_k'}]_1)+\dots+\overline{ev}_2(s_{s_0}([\Phi_{\delta_k'}]_{k-1})$$
with each term $\overline{ev}_2(s_{s_0}([\Phi_{\delta_k'}]_i)$, $i=1,2,\dots,k-1$ representing a map
$$D^1\times D^1 \times C_2[I]\to C_2[S^1\times D^3].$$
which we will apply Lemma \ref{cohorizontaldetect} to separately. Suppose that we fix the direction pointing horizontally to the left side being our choice of the cohorizontal direction $\delta$, then the points contained in the cohorizontal submanifolds $t^i\mathrm{Co}_1^2$ happen near the the cuffs when the strands of the bar inside the blue and red boxes (see Figure \ref{twisted arc}) get dragged to the farthest point respectively, i.e.~looping around the blue and red cuffs in Figure \ref{barbelltheta92}. At these moments, the vertical strands inside the red and blue boxes create cohorizontal pairs.

We now analyze the first intersection point (counting from the left to the right as shown in Figure \ref{barbelltheta92}) of $\delta_k'$ which we call a \textit{Type 1 point} whose name will become clear in due time. The two blue and two red arcs in Figure \ref{twisted arc} correspond to the moment when the farthest point is reached for both red and blue cuff spheres. We parameterize the interval by $[0,13]$. The two sets of strands contained in the red and blue boxes are looping around the red and blue circles, and each creates 8 cohorizontal points (so 16 in total). In this Type 1 case, there are no arcs on the red box side because the strands do not cross the spanning disks of the cuff spheres along the way of reaching the farthest point at the red cuff sphere. The next lemma describes cohorizontal points arising from this Type 1 intersection point.

\begin{lemma}
\label{secondstagepoints}
     The preimage of the cohorizontal points arising from the first intersection point of $\delta_k'$ are as shown in Figure \ref{twisted arc2}. The square on the left hand side is parameterized by $[0,1]\times [0,1]$, and the axes of the triangle picture are parameterized by $[0,13]$. The cohorizontal points appearing at the centre of the square are listed below.
 \begin{center}
\begin{tabular}{ c c c c}
 \textcolor{red}{(1,3)} &  \textcolor{red}{(3,7)}&  \textcolor{blue}{(2,4)}  &\textcolor{blue}{(2,10)} \\ 
  \textcolor{red}{(1,5)} &  \textcolor{red}{(5,7)} & \textcolor{blue}{(4,6)}& \textcolor{blue}{(6,10)} \\  
 \textcolor{red}{ (1,9)} &  \textcolor{red}{(7,9)}& \textcolor{blue}{(4,8)} & \textcolor{blue}{(8,10)} \\
  \textcolor{red}{(1,11)} &  \textcolor{red}{(7,11)}&\textcolor{blue}{(4,12)}& \textcolor{blue}{  (10,12)}.
\end{tabular}
\end{center}
 \end{lemma}
 
\begin{proof}
Figure \ref{twisted arc2} is a Budney--Gabai style picture as Figures 12 and 13 of \cite{Budney-gabai2} but we interpret it in a slightly different way. In the square, there is a cross of two intervals, one in blue and one in red. Cohorizontal points only appear in these two intervals. The two axes of the square are annotated as red and blue, meaning the former one controls the two strands in the red box, and the latter one controls the two strands in the blue box, in Figure \ref{twisted arc}.

\begin{figure}
    \centering
    \includegraphics[width=0.7\textwidth]{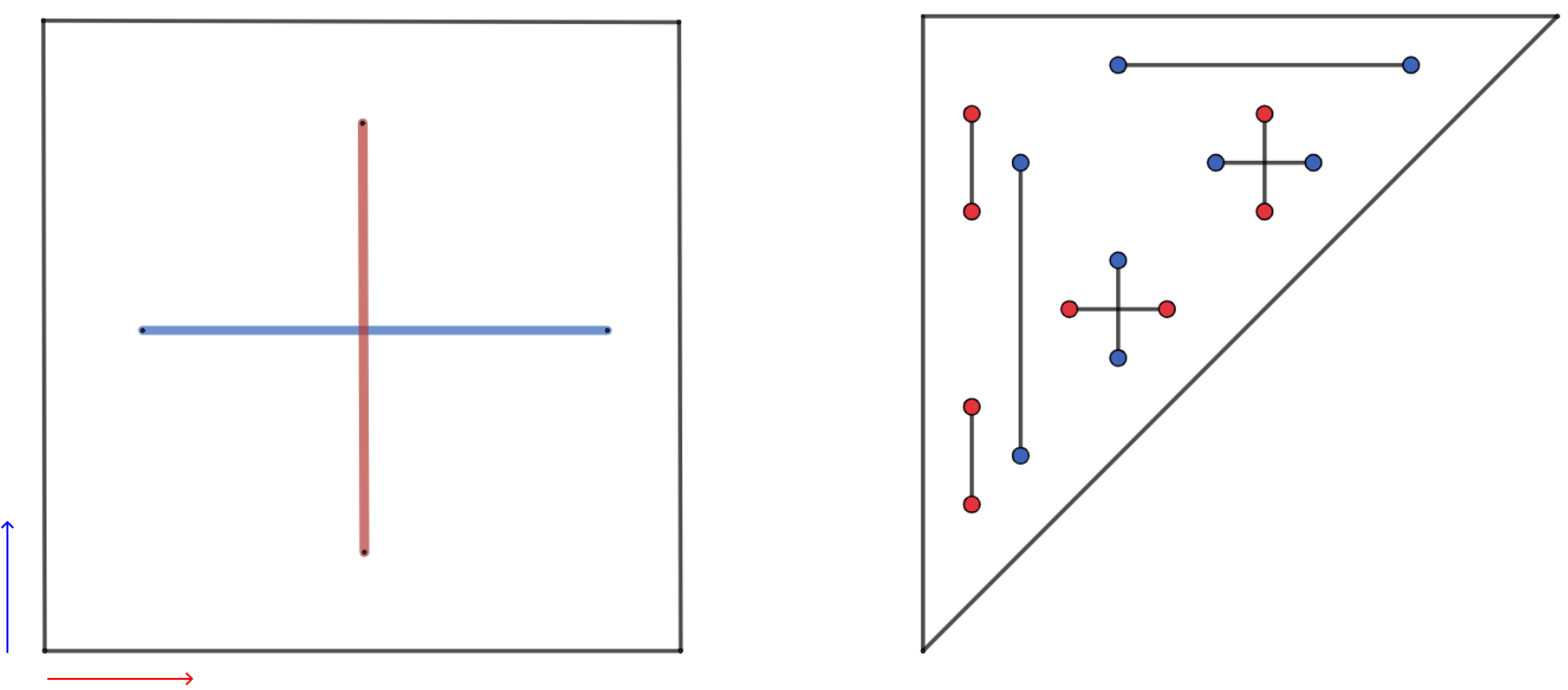}
    \caption{The preimage of the cohorizontal manifolds for the second stage map $D^1\times D^1\times C_2[I]\to C_2[S^1\times D^3]$ at a Type 1 point.}
    \label{twisted arc2}
\end{figure}

At the center of the square they intersect. This point corresponds to the moment when both of the two sets of strands contained in the red and blue boxes are at the farthest position, looping around the two cuff spheres. At this point, both of the red and blue dots in the triangle picture are present, corresponding to the double points mentioned in the paragraph before this lemma. 

When fixing the blue coordinate at $0.5$ and freely moving the red coordinate, i.e.~on the blue interval in the square, the blue cohorizontal points in the triangle picture stay. The 4 pairs of red cohorizontal points $$\{(1,3), (1,5)\},\{(1,9),(1,11)\},\{(3,7),(5,7)\},\{(7,9),(7,11)\}$$ gradually come together as we approache the two end points of the blue interval, tracing out a circle for each pair. This is shown as intervals connecting each pair in the triangle picture. After passing these two endpoints, no cohorizontal points are present. To see this, imagine that we fix the two strands that come from the blue box, keeping them at the farthest position, and move the two strands in the red box back to the initial position, then when these two strands cross the vertical strand that the red cuff $R$ links, the 4 pairs of red cohorizontal points merge together in pairs, and then disappear.

Similarly, varying the blue coordinate along the red interval traces out 4 circles of blue cohorizontal points corresponding to the pairs
$$\{(2,4),(2,10)\},\{(4,6),(4,8)\},\{(6,10),(8,10)\},\{(4,12),(10,12)\}.$$\end{proof}

Pulling back the orientations of the normal bundles of the cohorizontal submanifolds $t^i\mathrm{Co}_1^2$ (codimension 3) along $\overline{ev}_2(s_{s_0}([\Phi_{\delta_k'}]_1))$ gives orientations to the 8 circles. As every 1-manifold in $D^4$ is unlinked, this disjoint union of unlinked codimension 3 (1-dimensional) oriented submanifolds in $D^4$ is null-cobordant. This implies that the second stage map $\overline{ev}_2(s_{s_0}([\Phi_{\delta_k'}]_1))$ is null homotopic. 

Next, we move forward to detect the third stage map
$$\overline{ev}_3(s_{s_0}([\Phi_{\delta_k'}]_1))\colon D^1\times D^1\times C_3[I]\to C_3'[S^1\times D^3]$$ as an element in $\pi_5 C_3'[S^1\times D^3]$. By Lemma \ref{cohorizontal calculation}, we need to understand the preimage of the following pair: 
\begin{equation*}
    (t^{\alpha}\text{Co}^1_2-t^{\alpha-\beta}\text{Co}^1_3, t^{\beta-\alpha}\text{Co}^3_1-t^\beta\text{Co}^3_2)
\end{equation*}
for each $\alpha,\beta\in \mathbb{Z}$. As defined in Lemma \ref{cohorizontal calculation}, for each of the cohorizontal submanifolds in the above pair, one of the three coordinates is free to move with the remaining two being cohorizontal (remember that we assume no cohorizontal triples as in Lemma \ref{cohorizontal calculation}). This implies that each cohorizontal submanifold will only appear on two of the three facets of $C_3[I]$. In fact, these are represented by a disjoint union of 2-spheres embedded in $D^1\times D^1\times C_3[I]$. We draw $C_3[I]$ as a tetrahedron with coordinates $t_1,t_2$ and $t_3$. We follow the convention in \cite{Budney-gabai2} and draw a tetrahedron as a triangle with the $t_2$ coordinate pointing out of the page. We omit the square $D^1\times D^1$ in our pictures as it will stay the same as the square in Figure \ref{twisted arc2}. As outlined in Example 3.6 of \cite{Budney-gabai}, the third stage pictures are obtained from the second stage pictures by connecting the cohorizontal points in the preimage of the same cohorizontal submanifold at the same position on the two facets that are involved, by an interior arc, and closing the tetrahedron by attaching the null homotopies to the facets. The interior arc corresponds to varying the free coordinate (cf. Lemma \ref{cohorizontal calculation}).

Again, following the convention in \cite{Budney-gabai2}, we draw the preimage in pairs with different colours (that may produce non-trivial linking numbers) as in Figures \ref{calculate1}, \ref{calculate2}, \ref{calculate3}, \ref{calculate4}, \ref{calculate5} and \ref{calculate6}. These pictures show cohorizontal points at the centre of the square, i.e.~when we are at the point $(0,0)\in D^1\times D^1$. Along the boundary facets, we attach null homotopies of the second stage map with the convention that we close the red spheres first and the blue spheres afterwards. In our pictures, this gives rise to a collection of circles, and when we vary in the two directions of $D^1\times D^1$, these circles traces out a disjoint union of 2-spheres. One can refer to the construction of Figures 16 and 17 of \cite{Budney-gabai2} for comparison. Following the approach in \cite{Budney-gabai2}, we pull back the orientation of the normal bundles of the cohorizontal submanifolds (codimension 3) to give orientations to the 2-spheres (codimension 3). We now apply Lemma \ref{cohorizontal calculation}. For (the preimage of) each cohorizontal submanifold in the pair 
$(t^{\alpha}\text{Co}^1_2-t^{\alpha-\beta}\text{Co}^1_3, t^{\beta-\alpha}\text{Co}^3_1-t^\beta\text{Co}^3_2)$, the power of $t$ is determined by the number of times one needs to travel (along the bar) along the circle direction of $S^1\times D^3$ from one point of a cohorizontal pair to the other. We then read off the linking numbers from Figures \ref{calculate1}, \ref{calculate2}, \ref{calculate3}, \ref{calculate4}, \ref{calculate5} and \ref{calculate6} and take the sum to get the following polynomial: 
$$\overline{ev}_3(s_{s_0}([\Phi_{\delta_k'}]_1))=-t_1^{k-2}t_3^{k-2}-t_1^{k-2}t_3^{k-2}+t_1^{2-k}t_3^{0}+t_1^{0}t_3^{2-k}.$$
For example, in Figure \ref{calculate1}, there is only one linking pair (drawn in the middle of the picture). It takes $k-2$ times along the circle factor to go from the blue cuff sphere, along the bar, to where the cohorizontal pairs occur near the blue cuff (for $k=10$, this is near strand 8, counting from right to left on the top in Figure \ref{barbelltheta92}). Thus we have $\textcolor{blue}{\alpha=k-2}$. Similarly, we have $\textcolor{red}{\beta-\alpha=2-k}$. This implies that $\textcolor{red}{\beta=0}$, so this picture contributes $\pm t_1^{2-k}t_3^{0}$ with the sign being determined by the orientations of this linking pair. To determine the orientations, note that the null-homotopy denoted by the blue arc attached to the $t_2=t_3$ facet is determined by the pair of cohorizontal points $(4,6)$ and $(4,8)$ which admit negative and positive normal bundle orientations respectively, with reference to the pull-back orientation of the normal bundle of $\textcolor{blue}{t^\alpha\mathrm{Co}_2^1}$ in $C_3[S^1\times D^3]$ and with the direction
 pointing horizontally to the left side being our choice of the cohorizontal direction $\delta$, and thus is oriented by the direction pointing from $(4,6)$ to $(4,8)$. Similarly, the red arc representing a null-homotopy attached to the same $t_2=t_3$ facet is oriented by pointing from $(3,7)$ to $(5,7)$, leading to a positive linking. Therefore, Figure \ref{calculate1} gives rise to the monomial $t_1^{2-k}t_3^{0}$. The contributions from the remaining figures can be calculated in this manner as well.

    \begin{figure}
  \centering
  \begin{minipage}[b]{0.4\textwidth}
    \includegraphics[width=\textwidth]{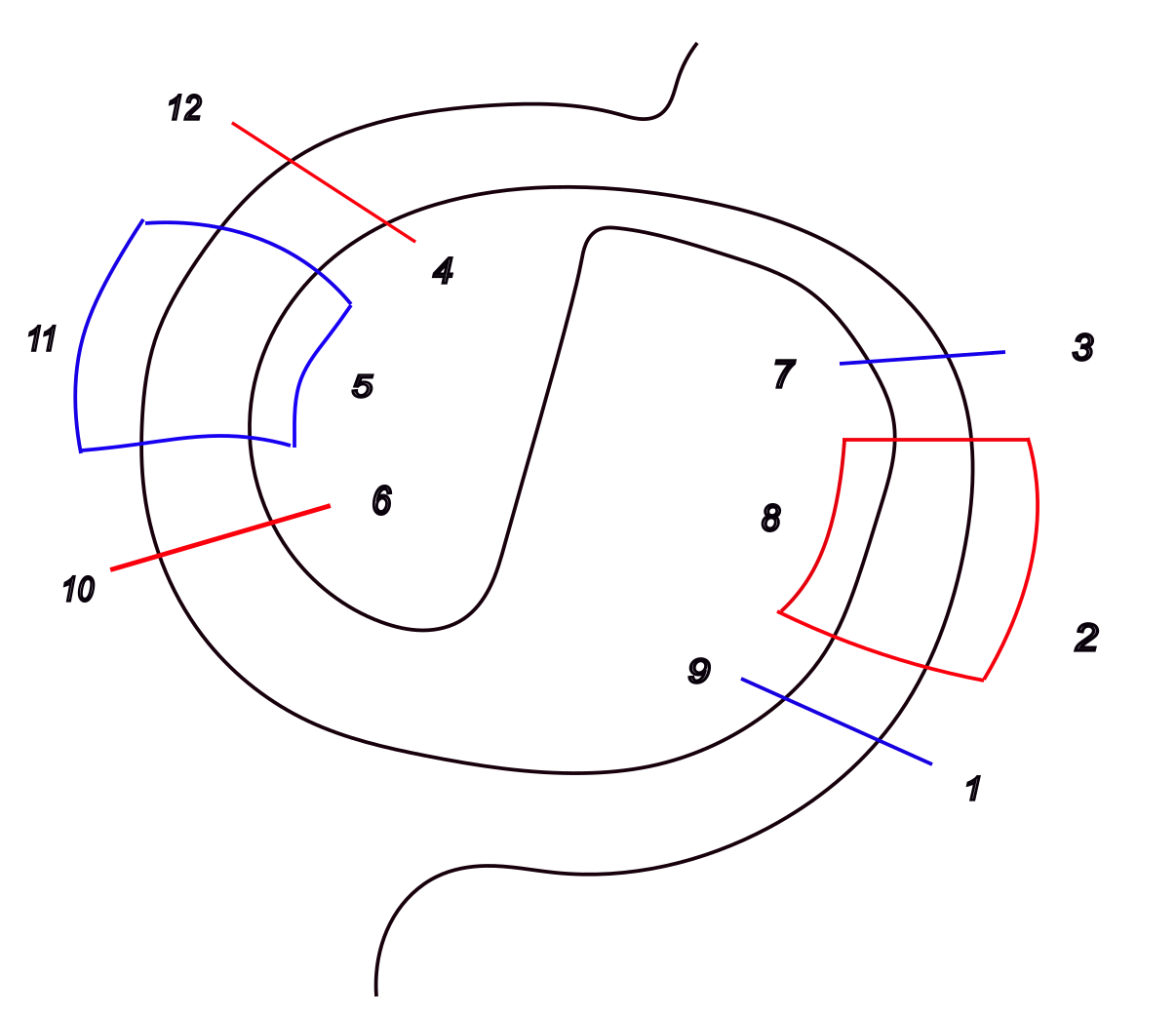}
    \caption{The Type 2 intersection for points 2 to $k-2$.}
    \label{type2qqq}
  \end{minipage}
  \hfill
  \begin{minipage}[b]{0.4\textwidth}
    \includegraphics[width=\textwidth]{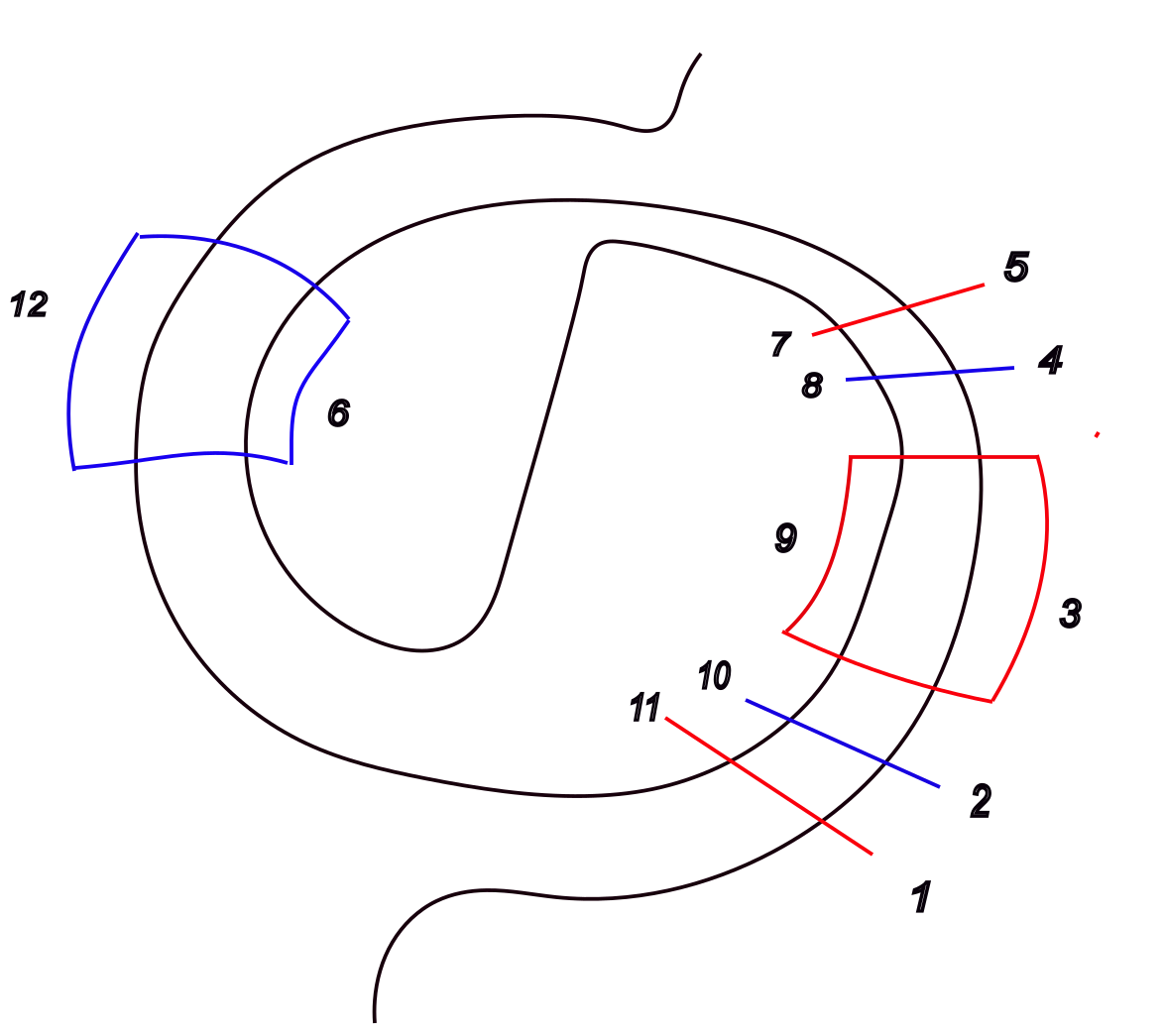}
    \caption{The Type 3 intersection for point $k-1$.}
    \label{type3qqq}
  \end{minipage}
\end{figure}

    \begin{figure}[!tbp]
  \centering
  \begin{minipage}[b]{0.46\textwidth}
    \includegraphics[width=\textwidth]{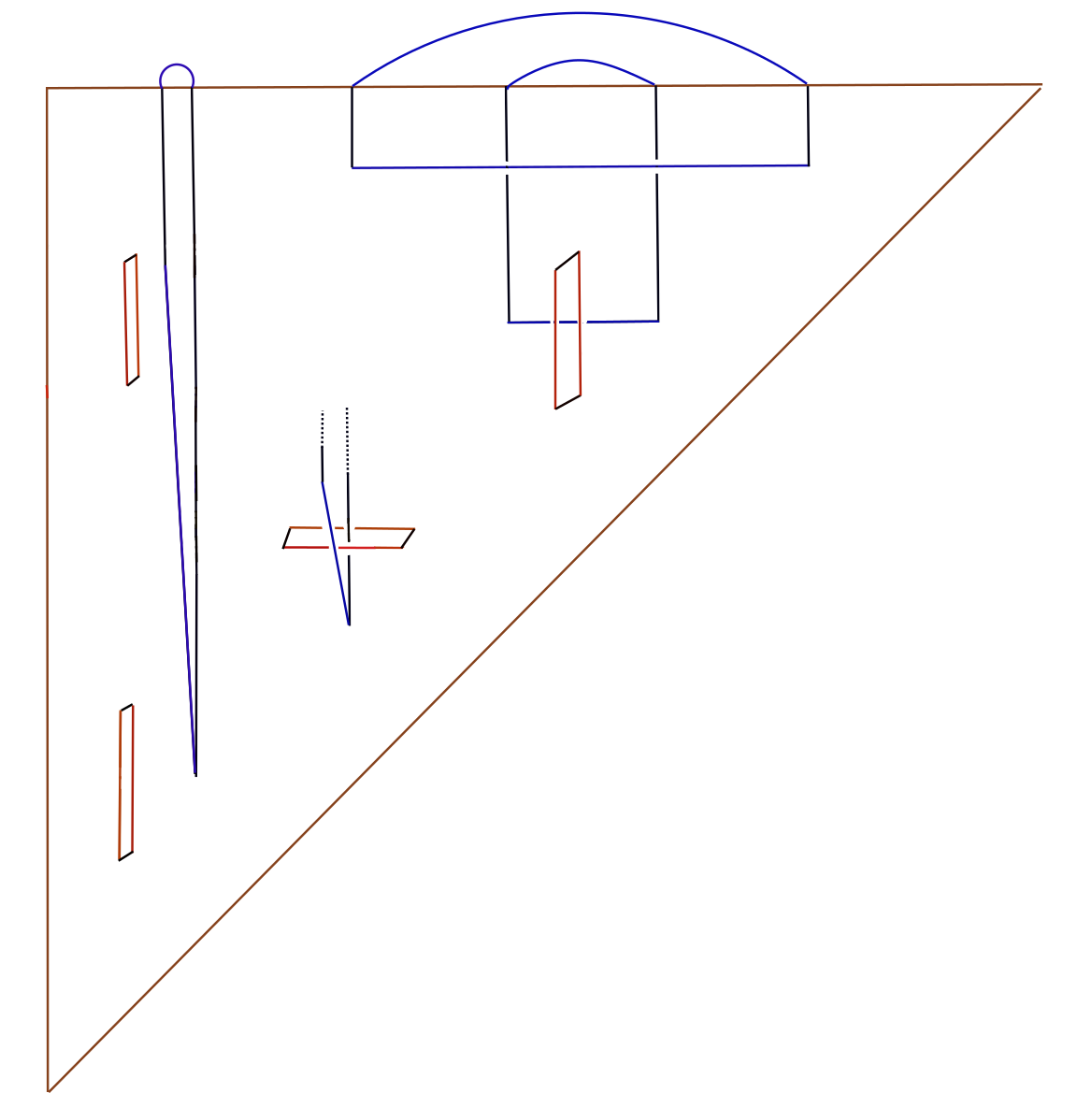}
    \caption{(\textcolor{blue}{$t^\alpha\mathrm{Co}_2^1$},\textcolor{red}{$t^{\beta-\alpha}\mathrm{Co}_1^3$}), contributes $t_1^{2-k}t_3^{0}$.}
      \label{calculate1}
  \end{minipage}
  \hfill
  \begin{minipage}[b]{0.44\textwidth}
    \includegraphics[width=\textwidth]{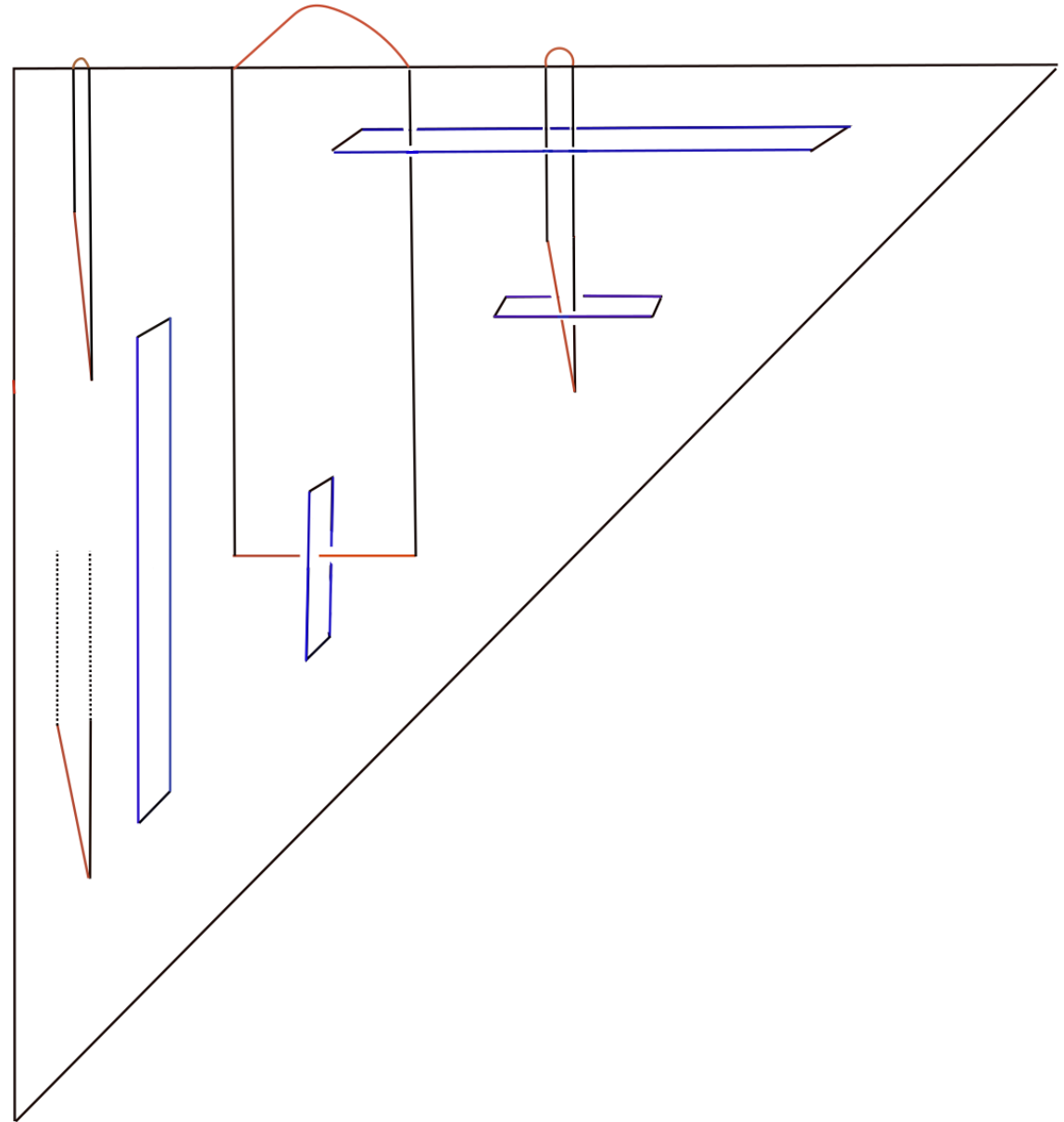}
    \caption{ (\textcolor{red}{$t^\alpha\mathrm{Co}_2^1$},\textcolor{blue}{$t^{\beta-\alpha}\mathrm{Co}_1^3$}), contributes 0.}
      \label{calculate2}
  \end{minipage}

\end{figure}

   \begin{figure}[!tbp]
  \centering
  \begin{minipage}[b]{0.5\textwidth}
    \includegraphics[width=\textwidth]{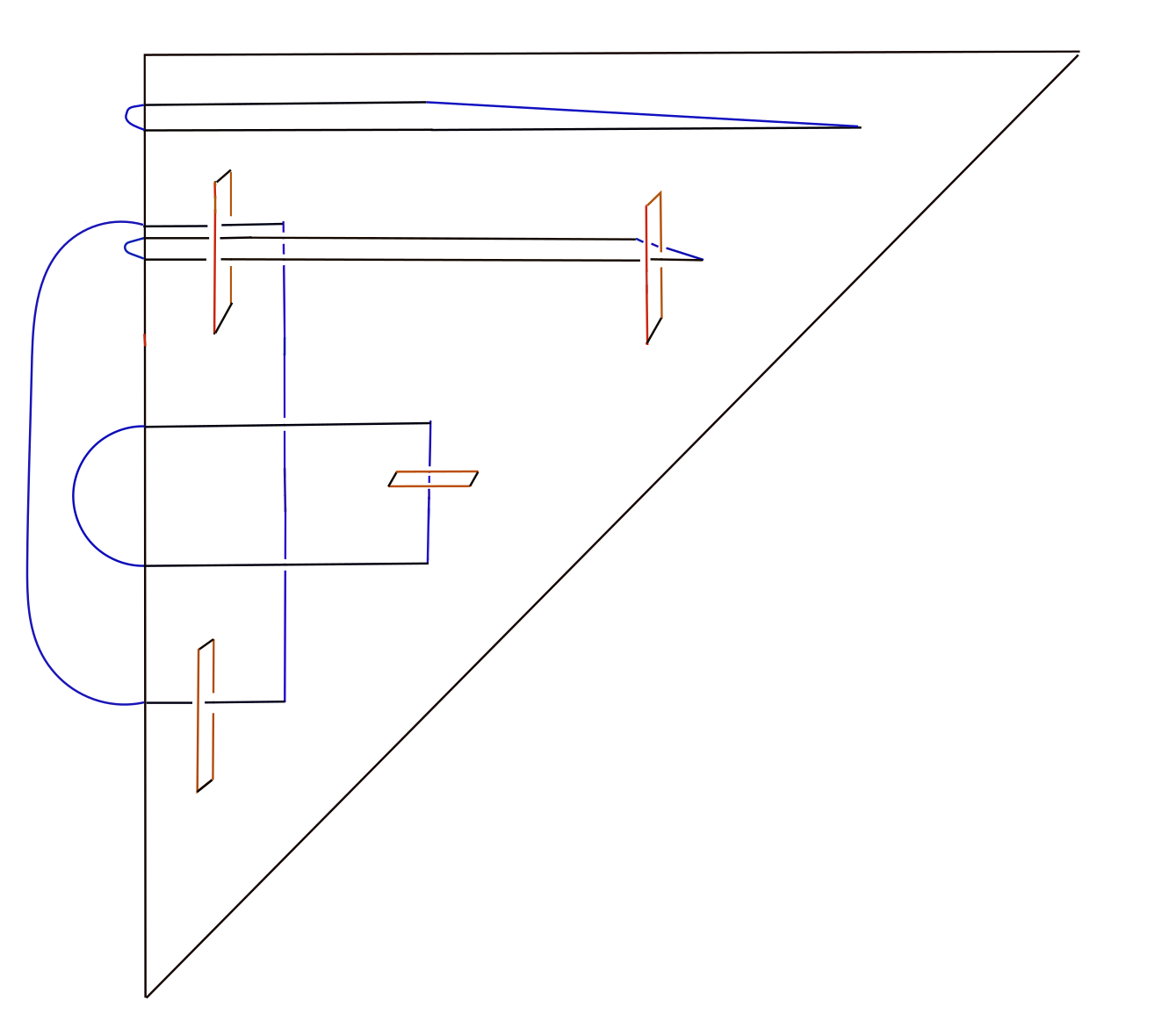}
    \caption{ (\textcolor{red}{$t^{\alpha-\beta}\mathrm{Co}_3^1$},\textcolor{blue}{$t^{\beta}\mathrm{Co}_2^3$}), contributes $t_1^0t_3^{2-k}$.}
      \label{calculate3}
  \end{minipage}
  \hfill
  \begin{minipage}[b]{0.48\textwidth}
    \includegraphics[width=\textwidth]{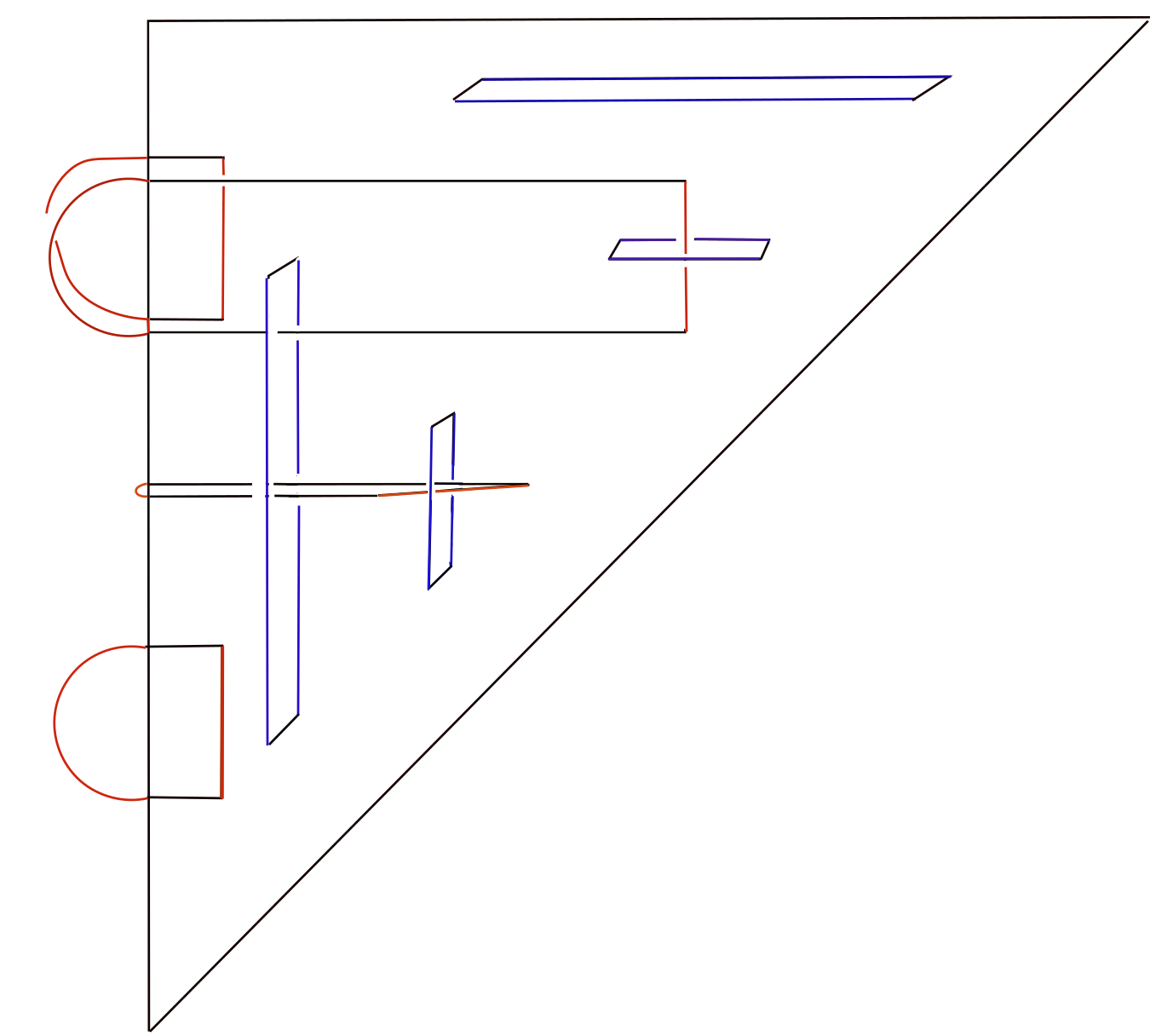}
    \caption{(\textcolor{blue}{$t^{\alpha-\beta}\mathrm{Co}_3^1$},\textcolor{red}{$t^{\beta}\mathrm{Co}_2^3$}), contributes 0.}
      \label{calculate4}
  \end{minipage}

\end{figure}

 \begin{figure}[!tbp]
  \centering
  \begin{minipage}[b]{0.47\textwidth}
    \includegraphics[width=\textwidth]{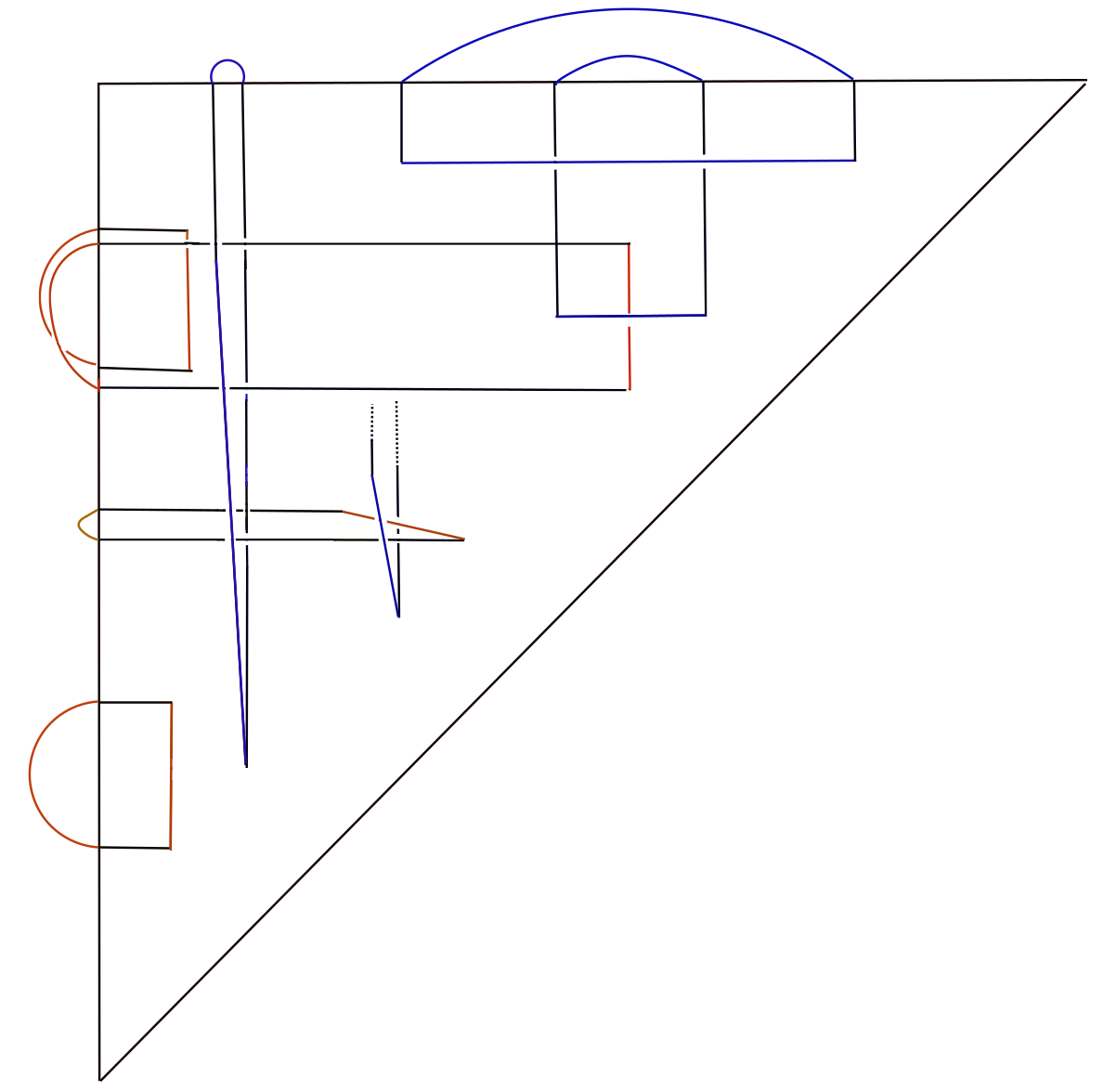}
    \caption{(\textcolor{blue}{$t^\alpha\mathrm{Co}_2^1$},\textcolor{red}{$t^{\beta}\mathrm{Co}_2^3$}), contributes $-t_1^{k-2}t_3^{k-2}$.}
      \label{calculate5}
  \end{minipage}
  \hfill
  \begin{minipage}[b]{0.49\textwidth}
    \includegraphics[width=\textwidth]{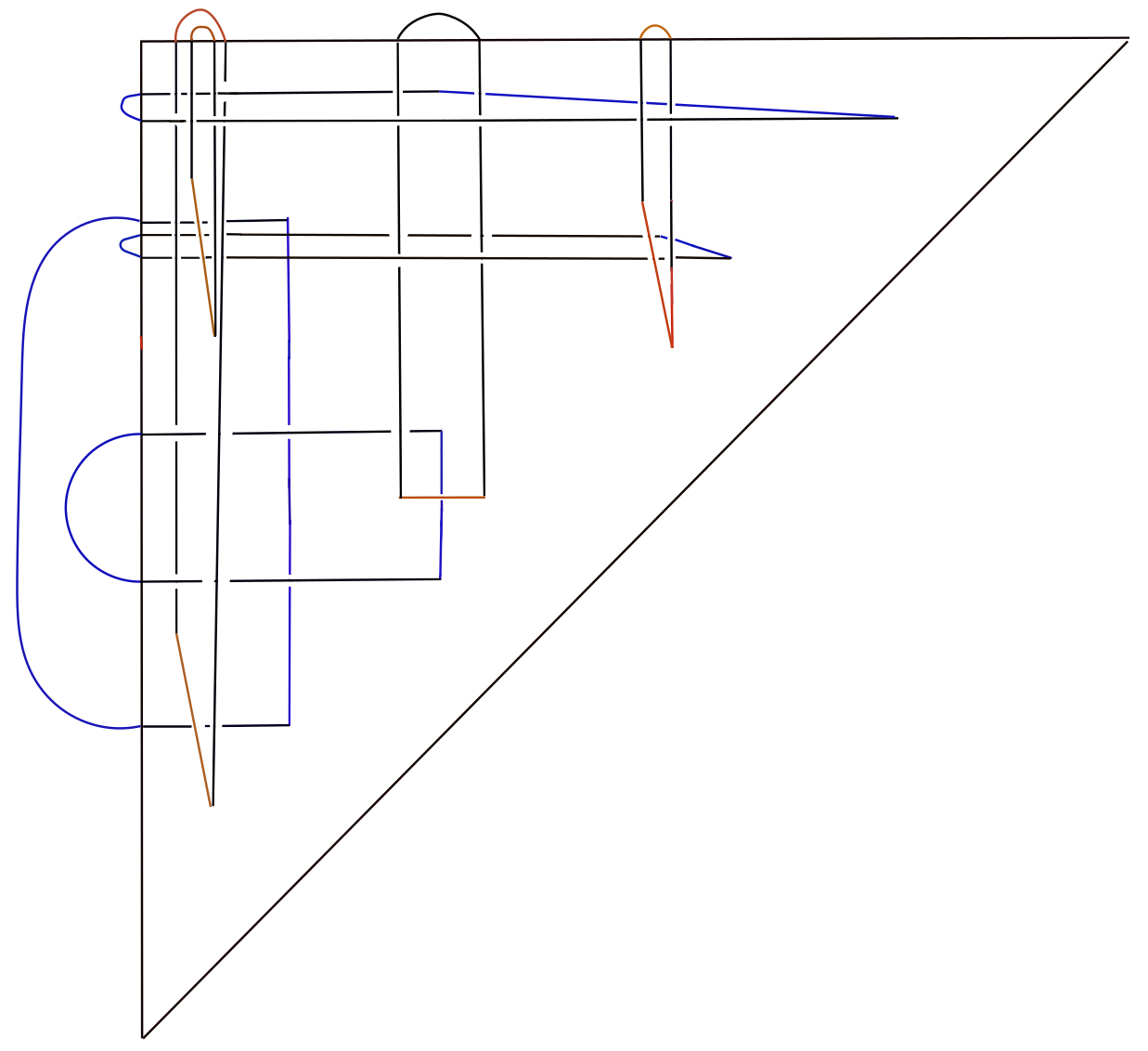}
    \caption{(\textcolor{red}{$t^\alpha\mathrm{Co}_2^1$},\textcolor{blue}{$t^{\beta}\mathrm{Co}_2^3$}), contributes $-t_1^{k-2}t_3^{k-2}$.}
      \label{calculate6}
  \end{minipage}

\end{figure}

To finish the calculation of $W_3([\Phi_{\delta_k'}])$, we need to calculate the contributions from the remaining $k-2$ intersection points. We first make the following important observation on intersection points of the bar with a chosen scanning disk for $\theta_k(e_i,e_j)$ with $i+j\geq k$.

\begin{lemma}
\label{typeintersection}
Assume that a barbell $\theta_k(e_i,e_j)$ in $S^1\times D^3$, $k\geq 1$ with $i+j\geq k$ is isotoped to a position such that the cuff spheres are disjoint from $\{s_0\}\times D^3\subset S^1\times D^3$ for some point $s_0\in S^1$, and the bar intersects $\{s_0\} \times D^3$ at $k-1$ points transversely, then these intersection points are classified into three types, based on the patterns of the cohorizontal pairs shown in Figures \ref{twisted arc}, \ref{type2qqq} and \ref{type3qqq}. The colored arcs on the sides of the blue and red boxes indicate intersections with the spanning disks of the two cuff spheres, i.e.~the natural 3-balls they bound. Furthermore, the barbell $\theta_k(e_i,e_j)$ produces $(k-i-1)$ Type 1 intersection points, $(-k+j+i+1)$ Type 2 intersection points and $(k-j-1)$ Type 3 intersection points.

\end{lemma}

\begin{proof}
For an intersection point (that is contained in an interval in $\{s_0\}\times D^3$) arising from $\theta_k(e_i,e_j)$, the corresponding barbell diffeomorphism creates intersections between this interval and the spanning disks of the two unknotted cuff spheres, i.e.~the natural 3-ball they bound in $S^1\times D^3$. Since $i+j\geq k$, the blue sphere is on the left of the red sphere, meaning that the nearest intersections that appear on the sphere of one colour are always of the other color. In fact, one can check point by point that there are only three possible patterns of such intersections as shown in Figures \ref{twisted arc}, \ref{type2qqq} and \ref{type3qqq}.

Since the nonzero coordinate of the first vector $e_i$ is at the $i$-th position, it follows that the blue cuff links the $i$-th vertical strand, counting from right to left. This indicates that counting from the left to right, the first $(k-i-1)$ intersection points are of Type 1. Similarly, the last $(k-j-1)$ intersection points are of Type 3 since $e_j$ indicates that the red cuff links the $j$-th vertical strand, counting from right to left. The remaining $(k-1)-(k-i-1)-(k-j-1)=(-k+j+i+1)$ intersection points in the middle are of Type 2.
\end{proof}

The following result by \cite{Budney-gabai2}) concerns the barbell $\delta_k=\theta_k(e_{k-1},e_{k-2})$ (see p.26 of \cite{Budney-gabai2}).

\begin{lemma}[\cite{Budney-gabai2}]
\label{deltabudney}
    The barbell $\delta_k=\theta_k(e_{k-1},e_{k-2})$ with $k\geq 4$ has one Type 3 intersection point that contributes
    $$-t_1^{2-k}t_3^{1-k}-t_1^{2-k}t_3^{1-k}-t_1^{2-k}t_3^{1}+t_1^{1-k}t_3^{-1}+t_1^{k-1}t_3^{1}+t_1^{1}t_3^{k-1}+t_1^{-1}t_3^{1-k}-t_1^1t_3^{2-k}$$
    and $k-2$ Type 2 intersection points, each of which contributes 
    $$t_1^{-1}t_3^{1-k}+t_1^{1-k}t_3^{-1}-t_1^{2-k}t_3^1-t_1^1t_3^{2-k}.$$
    
\end{lemma}

Lemma \ref{deltabudney} concerns the barbell $\delta_k$ which has no Type 1 but only Type 2 and Type 3 intersection points. Although $\delta_k$ and $\delta_k'$ are different barbells, by Lemma \ref{typeintersection}, the Type 2 intersection points of both barbells can be calculated following the same process, but with modifications to the exponents of $t_1$ and $t_3$ terms based on the number of times one needs to travel (along the bar) along the circle direction of $S^1\times D^3$ from one point of a cohorizontal pair to the other, as mentioned before. For example, in the top left picture of Figure 16 of \cite{Budney-gabai2}, the linking number of the preimage of the pair (\textcolor{red}{$t^\alpha\mathrm{Co}_2^1$},\textcolor{blue}{$t^{\beta-\alpha}\mathrm{Co}_1^3$}) was calculated. We can use the same picture to calculate the same linking number for $\delta_k'$ instead of $\delta_k$ by changing the exponents of $t$ by changing $1-k$ to $2-k$, and $k-1$ to $k-2$, leading to $-t_1^{2-k}t_3^0$ rather than $-t_1^{2-k}t_3^1$ in \cite{Budney-gabai2}. Therefore, using Figures 16 and 17 of \cite{Budney-gabai2}, we can calculate $W_3$ of $\delta_k'$ in full. In particular, one can verify that the $k-3$ Type 2 intersection points contribute 0 (more to discuss below in Lemma \ref{mainformula}) and the last intersection point contributes:
$$t_1^0t_3^{k-2}+t_1^{k-2}t_3^0-t_1^{2-k}t_3^{2-k}-t_1^{2-k}t_3^{2-k}.$$
Again, changes to the $t$ powers need to be made to the formulas in \cite{Budney-gabai2} based on the number of times a cohorizontal point needs to travel (along the bar) along the circle direction to get to the cuff sphere with the same colour (with signs accounted by the pulled back orientations as in \cite{Budney-gabai2}). Taking the sum of the above two polynomials and applying the Hexagon relation twice, we deduce that $W_3(\delta_k')=0$ for $k\geq 3$.

\begin{proposition}
    For $k\geq 3$, we have 
    \begin{align*}
        & W_3(\delta_k')=
        -t_1^{k-2}t_3^{k-2}-t_1^{k-2}t_3^{k-2}+t_1^{2-k}t_3^{0}+t_1^{0}t_3^{2-k}+
        t_1^0t_3^{k-2}+t_1^{k-2}t_3^0-t_1^{2-k}t_3^{2-k}-t_1^{2-k}t_3^{2-k}\\
        &=-(t_1^{k-2}t_3^{k-2}+t_1^{k-2}t_3^{k-2})-(t_1^{2-k}t_3^{2-k}+t_1^{2-k}t_3^{2-k})+(t_1^{2-k}t_3^{0}+t_1^0t_3^{k-2})+(t_1^{0}t_3^{2-k}+
        t_1^{k-2}t_3^0)=0.
    \end{align*}
\end{proposition}

We now extend the calculations in to $\theta_k(e_i,e_j)$ with $i,j,k\in\mathbb{Z}^{+}$. By Lemma \ref{typeintersection}, we observe that for $\theta_k(e_i,e_j)$ with $1\leq i\leq j\leq k-1$ and $k\geq 3$, the shapes and patterns of the triangle and tetrahedron pictures arising from the intersection between the bar and the scanning disk coincide with Figures \ref{calculate1}, \ref{calculate2}, \ref{calculate3}, \ref{calculate4}, \ref{calculate5} and \ref{calculate6} we presented in the last section (together with Figures 16 and 17 presented in Section 3 of \cite{Budney-gabai2}). If we use the convention that $e_i$ controls the position of the red cuff and $e_j$ controls the position of the blue cuff, then the exponents $\alpha,\beta,\alpha-\beta,\beta-\alpha$ of $t$ in each of the cohorizontal submanifolds in the pair 
$$(t^{\alpha}\text{Co}^1_2-t^{\alpha-\beta}\text{Co}^1_3, t^{\beta-\alpha}\text{Co}^3_1-t^\beta\text{Co}^3_2),$$
can be described in terms of $i$ and $j$ (via colour constituents as proposed on p.26 of \cite{Budney-gabai2}): $\pm i$ for the red constituent part and $\pm j$ for the blue constituent part. Therefore, we have the following theorem.

\begin{figure}
  \centering
  \begin{minipage}{0.45\textwidth}
   \centering
    \includegraphics[width=\linewidth]{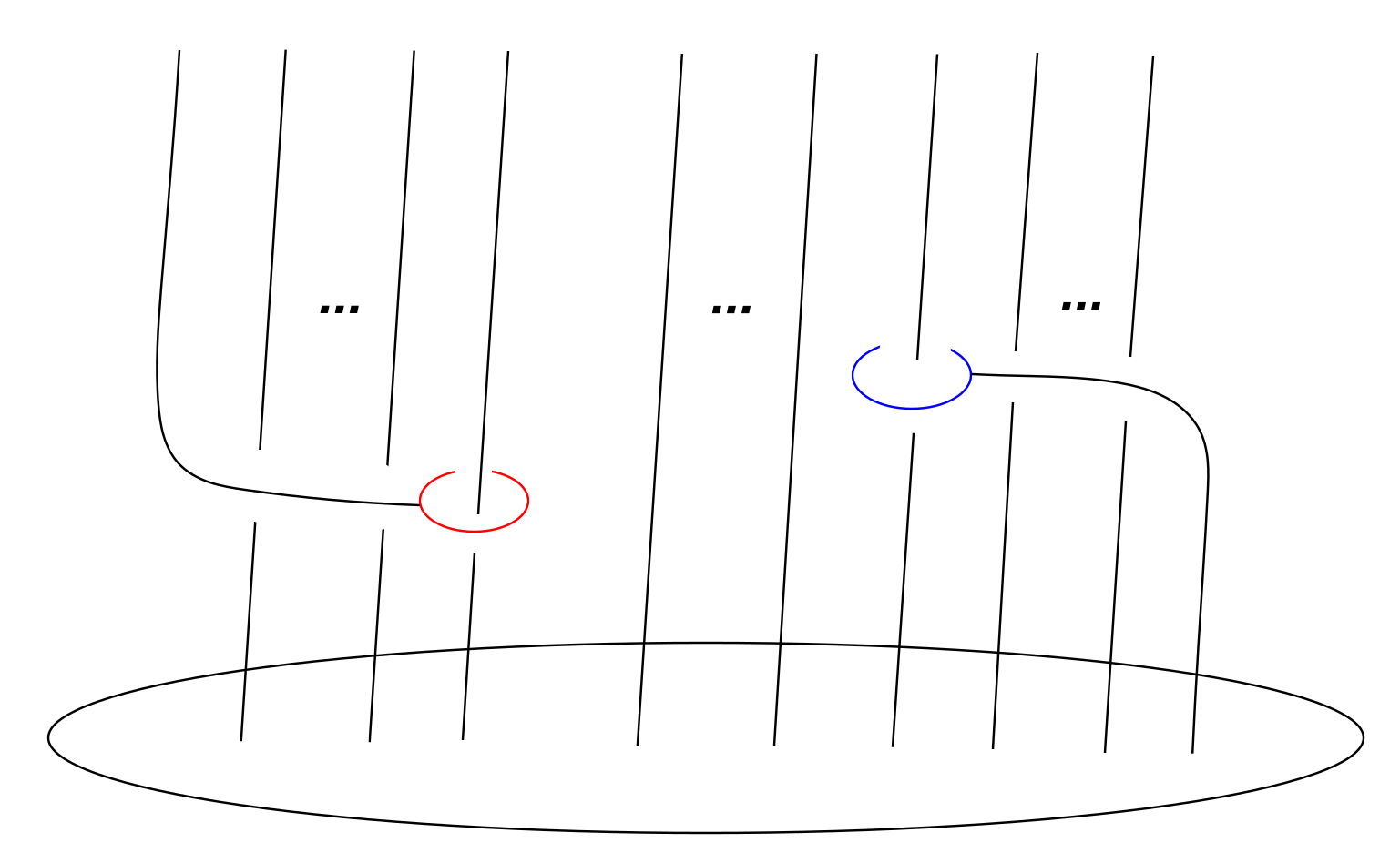}
        \caption{$\theta_k(e_i,e_j)$ with $i+j<k$.}
        \label{newargument2}
  \end{minipage}
  \hfill
  \begin{minipage}{0.45\textwidth}
   \centering
        \includegraphics[width=\linewidth]{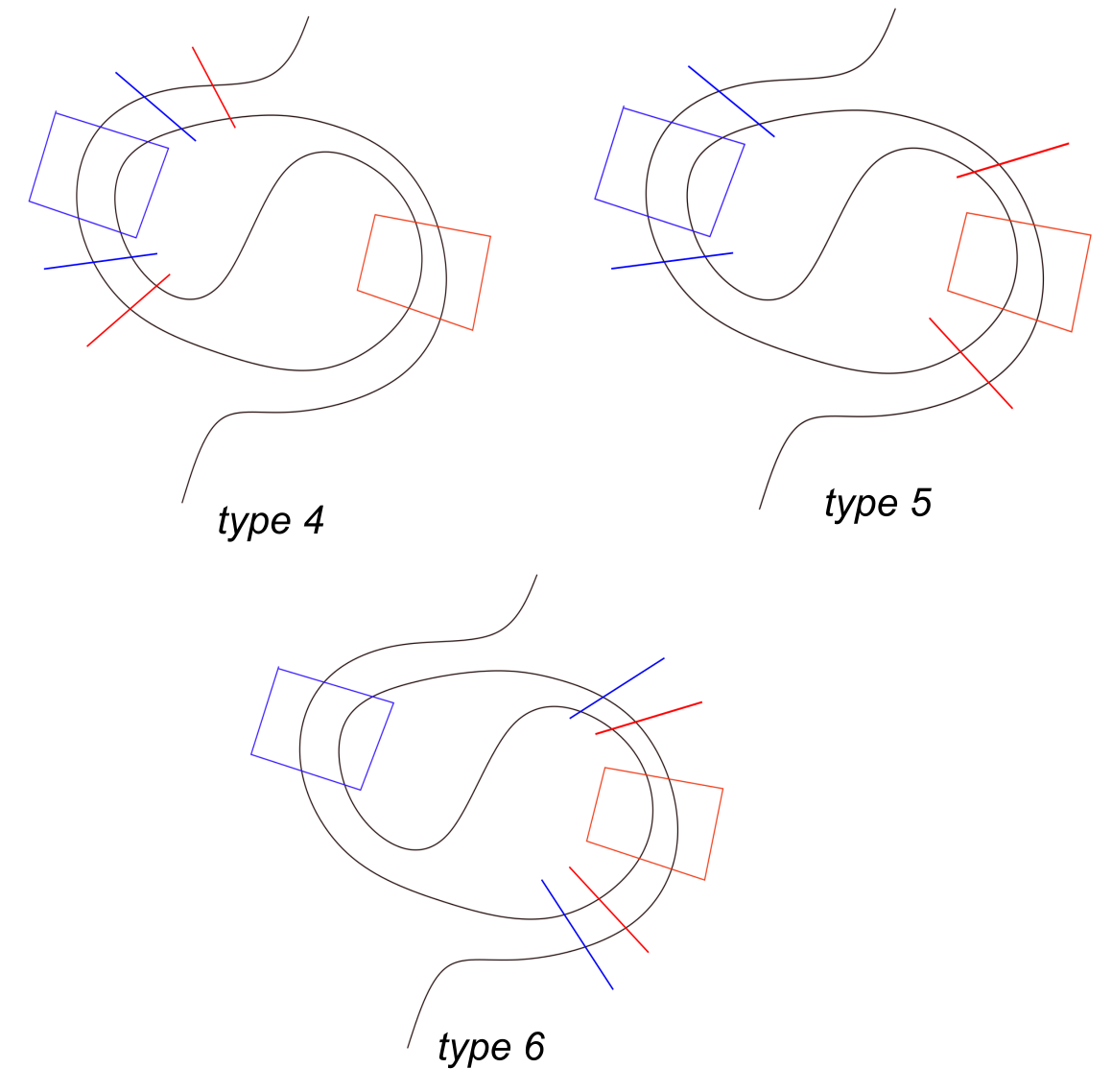}
        \caption{The extra three points of intersection points for $\theta_k(e_i,e_j)$ with $i+j<k$.}
        \label{newtypes}
  \end{minipage}

\end{figure}

 \begin{theorem}
 \label{mainformula}

     For $k\geq 3$, the polynomial we obtain from a Type 1 intersection from $\theta_k(e_i,e_j)$ with $i+j\geq k$ is
$$\mathrm{T}_1(i,j)=-t_1^{j}t_3^{i}-t_1^{i}t_3^{j}+t_1^{-i}t_3^{j-i}+t_1^{j-i}t_3^{-i}.$$
Based on Budney--Gabai's calculations, we also have the formula for a Type 2 intersection point:
$$\mathrm{T}_2(i,j)=-t_1^{-j}t_3^{i-j}+t_1^{-i}t_3^{j-i}+t_1^{j-i}t_3^{-i}-t_1^{i-j}t_3^{-j}$$
and for a Type 3 intersection point:
$$\mathrm{T}_3(i,j)=-t_1^{-j}t_3^{-i}-t_1^{-i}t_3^{-j}-t_1^{-j}t_3^{i-j}+t_1^{-i}t_3^{j-i}+t_1^{i}t_3^{i-j}+t_1^{i-j}t_3^{i}+t_1^{j-i}t_3^{-i}-t_1^{i-j}t_3^{-j}.$$
 \end{theorem}
 Combining Lemma \ref{typeintersection} and Theorem \ref{mainformula} allows us to write down the $W_3$ invariant for any $\theta_k(e_i,e_j)$ with $k\geq 3$ and $i+j\geq k$. For example, one can verify that when $i=j=k-1$, 
 $$W_3(\Phi_{\delta_k'})=W_3(\Phi_{\theta_k(e_{k-1},e_{k-1})})=0.$$ 
So far we have calculated $W_3$ of $\theta_k(e_i,e_j)$ with $i+j\geq k$. We now consider the case where $\theta_k(e_i,e_j)$ with $j\leq k-1-i$ (see Figure \ref{newargument2}) following a similar strategy. These are barbells such that there is no ``overlap'' between the two cuff spheres, and there are three more types of intersection points that need to be considered.  In particular, in Figure \ref{newargument2}, if we use the same scanning disk as before, then points on the left side of the red cuff (including the one just below it) belong to Type 4, points on the right side of the blue cuff (including the one just below it) belong to Type 6, and points in between the two cuff spheres belong to Type 5. The scanning pictures of these three extra types are shown in Figure \ref{newtypes}. We further observe the following lemma which simplifies our calculations by eliminating Type 5 intersection points.

    \begin{figure}
  \centering
  \begin{minipage}{0.38\textwidth}
   \centering
    \includegraphics[width=\textwidth]{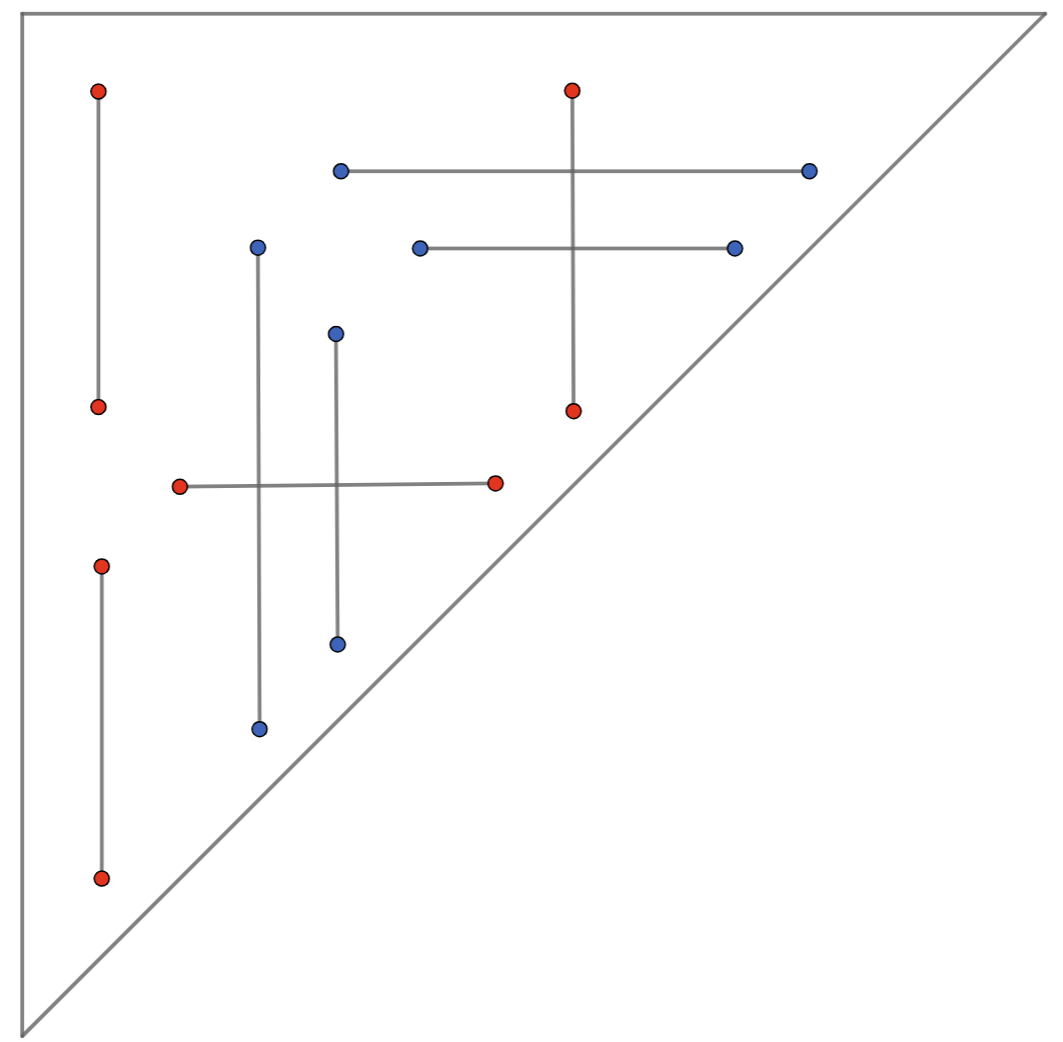}
    \caption{Cohorizontal manifolds in the second stage map for a Type 4 intersection point.}
    \label{type4}
  \end{minipage}
  \hfill
  \begin{minipage}{0.4\textwidth}
   \centering
    \includegraphics[width=\textwidth]{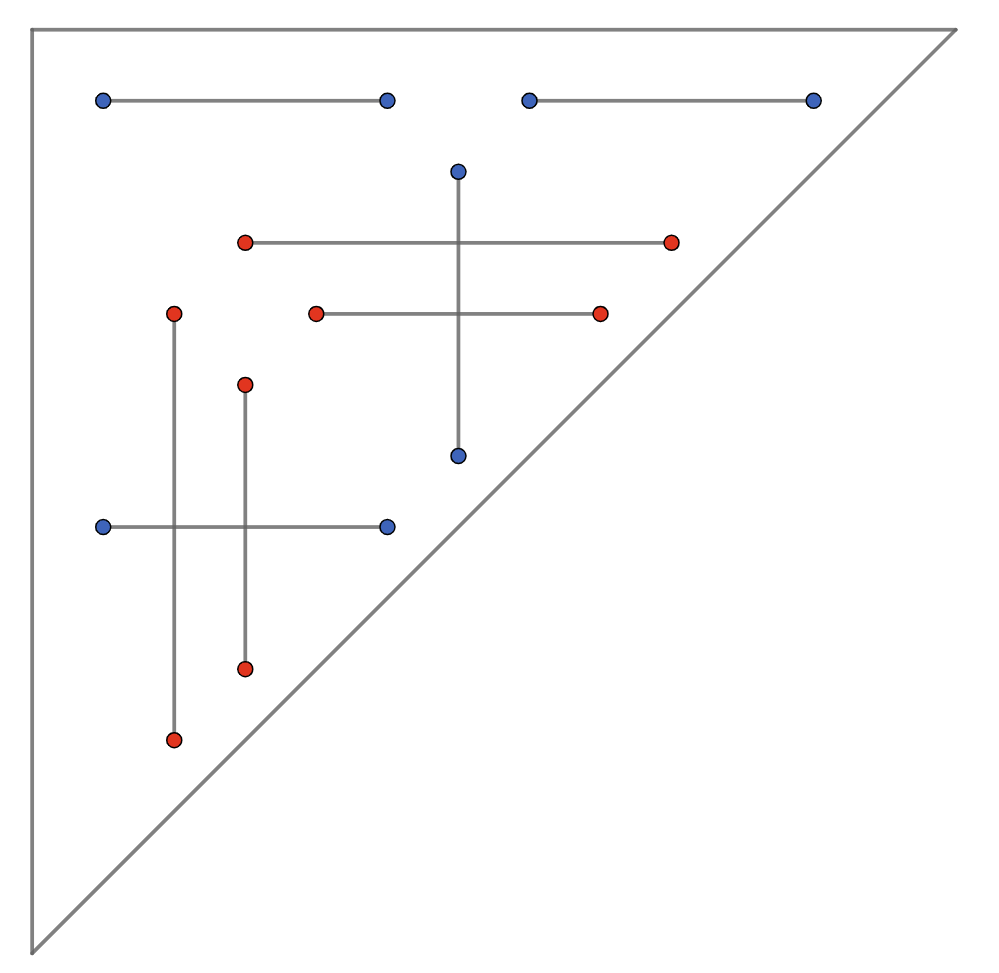}
    \caption{Cohorizontal manifolds in the second stage map for a Type 6 intersection point.}
    \label{type6}
  \end{minipage}

\end{figure}

  \begin{lemma}
  \label{type4and6barbell}
      For $j\leq k-1-i$ and $i,j,k\in \mathbb{Z}^{+}$, $\theta_k(e_i,e_j)=\theta_{i+j+1}(e_i,e_j)$. In addition, $\theta_{i+j+1}(e_i,e_j)$ has $j$ Type 4 intersection points and $i$ Type 6 intersection points.
  \end{lemma}
  \begin{proof}
      Since $k-i-j\geq 1$, the blue cuff (controlled by $e_j$ by convention) can be isotoped along the bar through the negative circle direction $k-i-j-1$ times until the blue cuff and the red cuff link the same vertical strand (with the blue on top). The lemma then follows from renaming the indices.
  \end{proof}

  In particular, we can avoid Type 5 points in the calculations and only consider Type 4 and Type 6 intersection points. We now describe the cohorizontal submanifolds for Type 4 and 6 intersection points in the second stage. We omit the square $D^1\times D^1$ part in our pictures as it will stay the same as the square in Figure \ref{twisted arc2}. Again, we parametrize the edges of the triangles by $[0,13]\times [0,13]$. The next lemma can be obtained in the same manner as for Lemma \ref{secondstagepoints}.

  \begin{lemma}
      
  In the second stage, the cohorizontal points for a Type 4 intersection point in the domain of the map
    $$D^1\times D^1 \times C_2[I]\to C_2[S^1\times D^3]$$
    are given by \begin{center}
\begin{tabular}{ c c c c}
 \textcolor{red}{(1,2)} &  \textcolor{red}{(1,6)}&  \textcolor{blue}{(3,4)}  &\textcolor{blue}{(3,10)} \\ 
  \textcolor{red}{(1,8)} &  \textcolor{red}{(1,12)} & \textcolor{blue}{(4,5)}& \textcolor{blue}{(4,9)} \\  
 \textcolor{red}{ (2,7)} &  \textcolor{red}{(6,7)}& \textcolor{blue}{(5,10)} & \textcolor{blue}{(9,10)} \\
  \textcolor{red}{(7,8)} &  \textcolor{red}{(7,12)}&\textcolor{blue}{(4,11)}& \textcolor{blue}{  (10,11)}.
\end{tabular}
\end{center}
and is illustrated by Figure \ref{type4}. Similarly, the cohorizontal points for a Type 6 intersection point 
are given by \begin{center}
\begin{tabular}{ c c c c}
 \textcolor{red}{(2,3)} &  \textcolor{red}{(2,9)}&  \textcolor{blue}{(1,12)}  &\textcolor{blue}{(5,12)} \\ 
  \textcolor{red}{(3,4)} &  \textcolor{red}{(3,8)} & \textcolor{blue}{(7,12)}& \textcolor{blue}{(11,12)} \\  
 \textcolor{red}{ (3,10)} &  \textcolor{red}{(9,10)}& \textcolor{blue}{(1,6)} & \textcolor{blue}{(5,6)} \\
  \textcolor{red}{(4,9)} &  \textcolor{red}{(8,9)}&\textcolor{blue}{(6,7)}& \textcolor{blue}{  (6,11)}.
\end{tabular}
\end{center}
and is illustrated in Figure \ref{type6}. Fixing one of the colored cuffs and moving the other leads to the black arcs connecting pairs of cohorizontal points just as in cases for types 1, 2 and 3, resulting in in 8 oriented circles as cohorizontal submanifolds.
\end{lemma}

\begin{figure}[!tbp]
  \centering
  \begin{minipage}[b]{0.43\textwidth}
    \includegraphics[width=\textwidth]{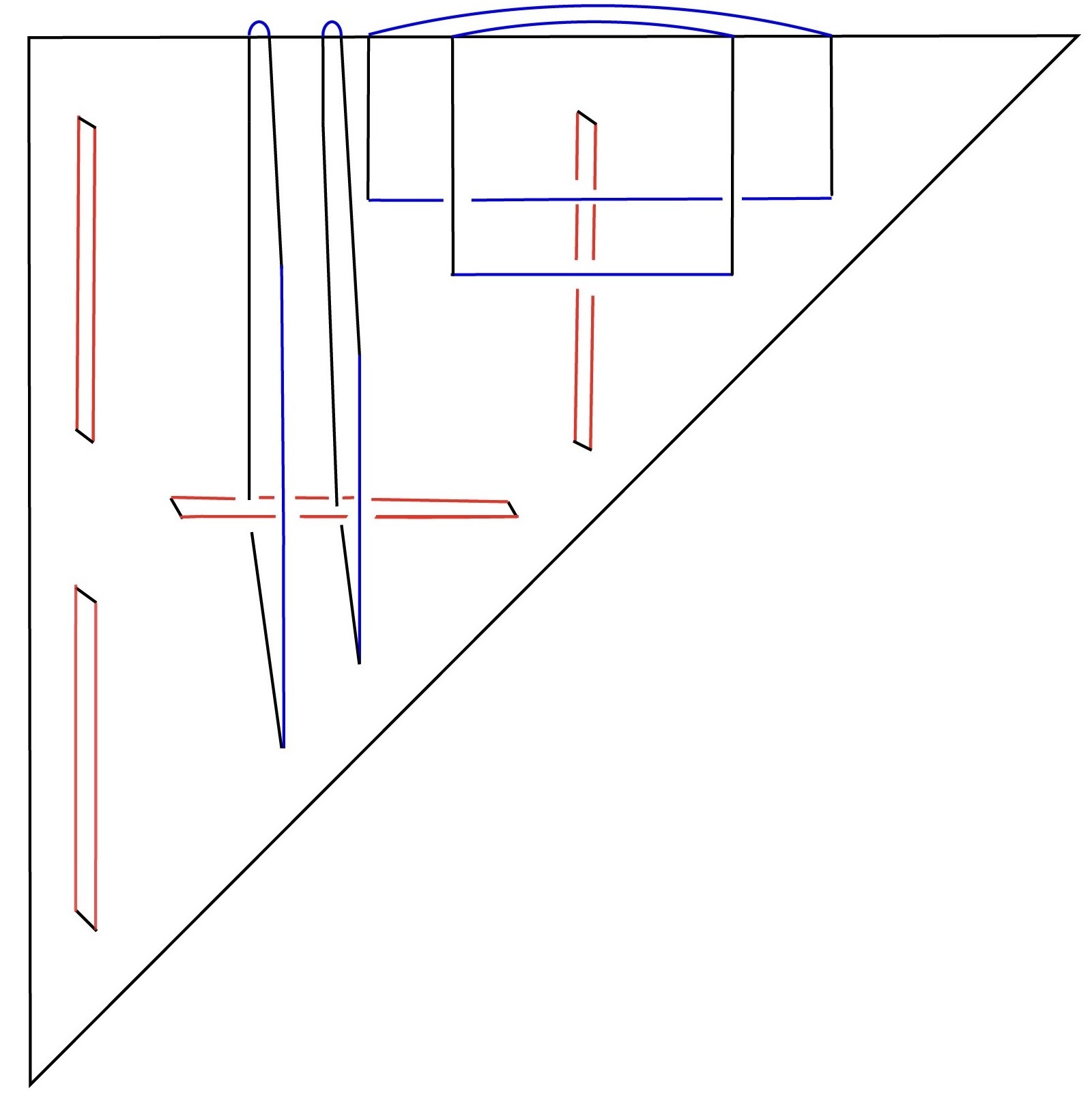}
    \caption{Type 4 (\textcolor{blue}{$t^\alpha\mathrm{Co}_2^1$},\textcolor{red}{$t^{\beta-\alpha}\mathrm{Co}_1^3$}), contributes $t_1^it_3^{j+i}+t_1^{-i}t_3^{j-i}$.}
      \label{calculate1type4}
  \end{minipage}
  \hfill
  \begin{minipage}[b]{0.44\textwidth}
    \includegraphics[width=\textwidth]{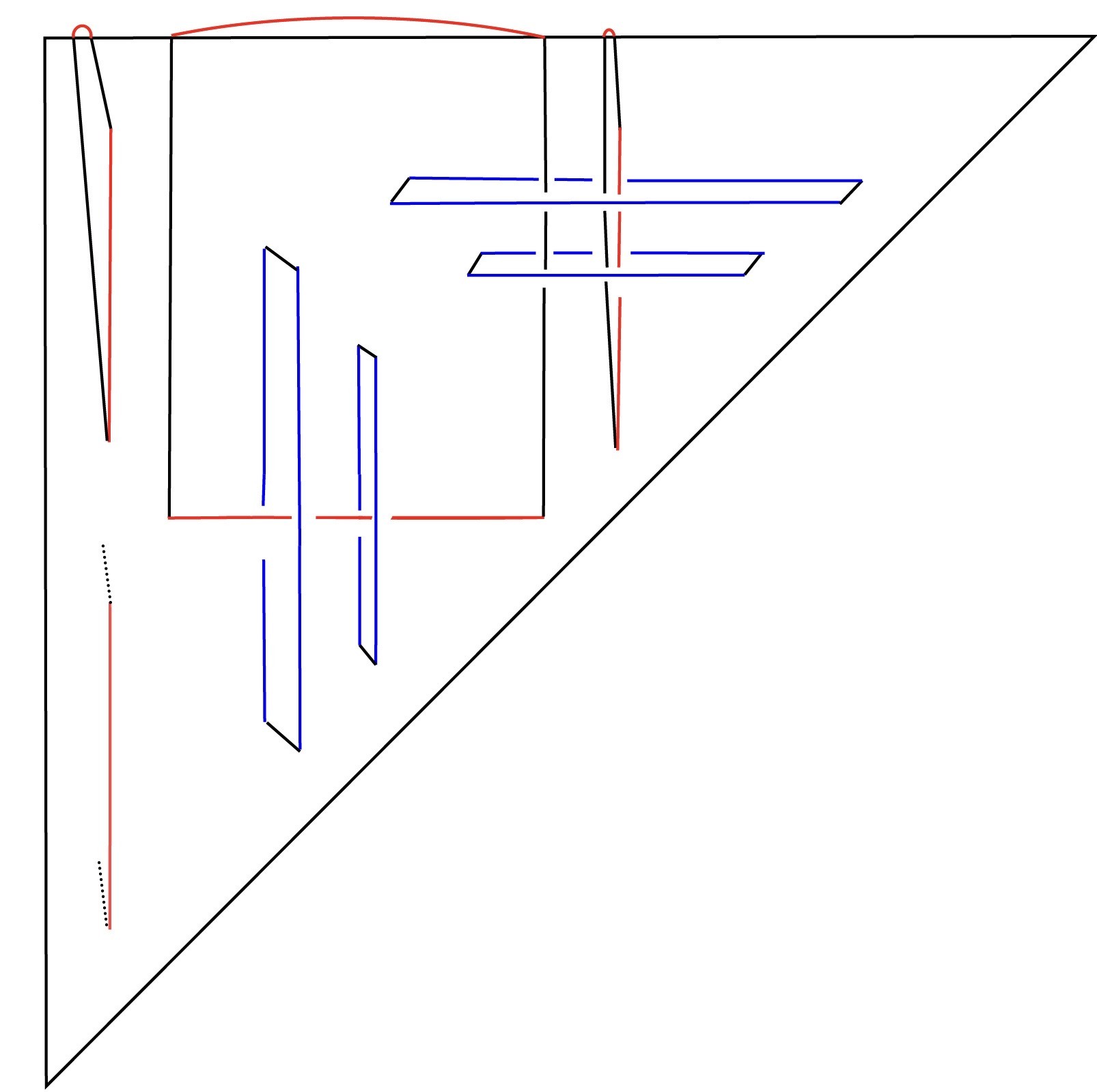}
    \caption{Type 4 (\textcolor{red}{$t^\alpha\mathrm{Co}_2^1$},\textcolor{blue}{$t^{\beta-\alpha}\mathrm{Co}_1^3$}), contributes 0.}
      \label{calculate2type4}
  \end{minipage}

\end{figure}

\begin{figure}[!tbp]
  \centering
  \begin{minipage}[b]{0.44\textwidth}
    \includegraphics[width=\textwidth]{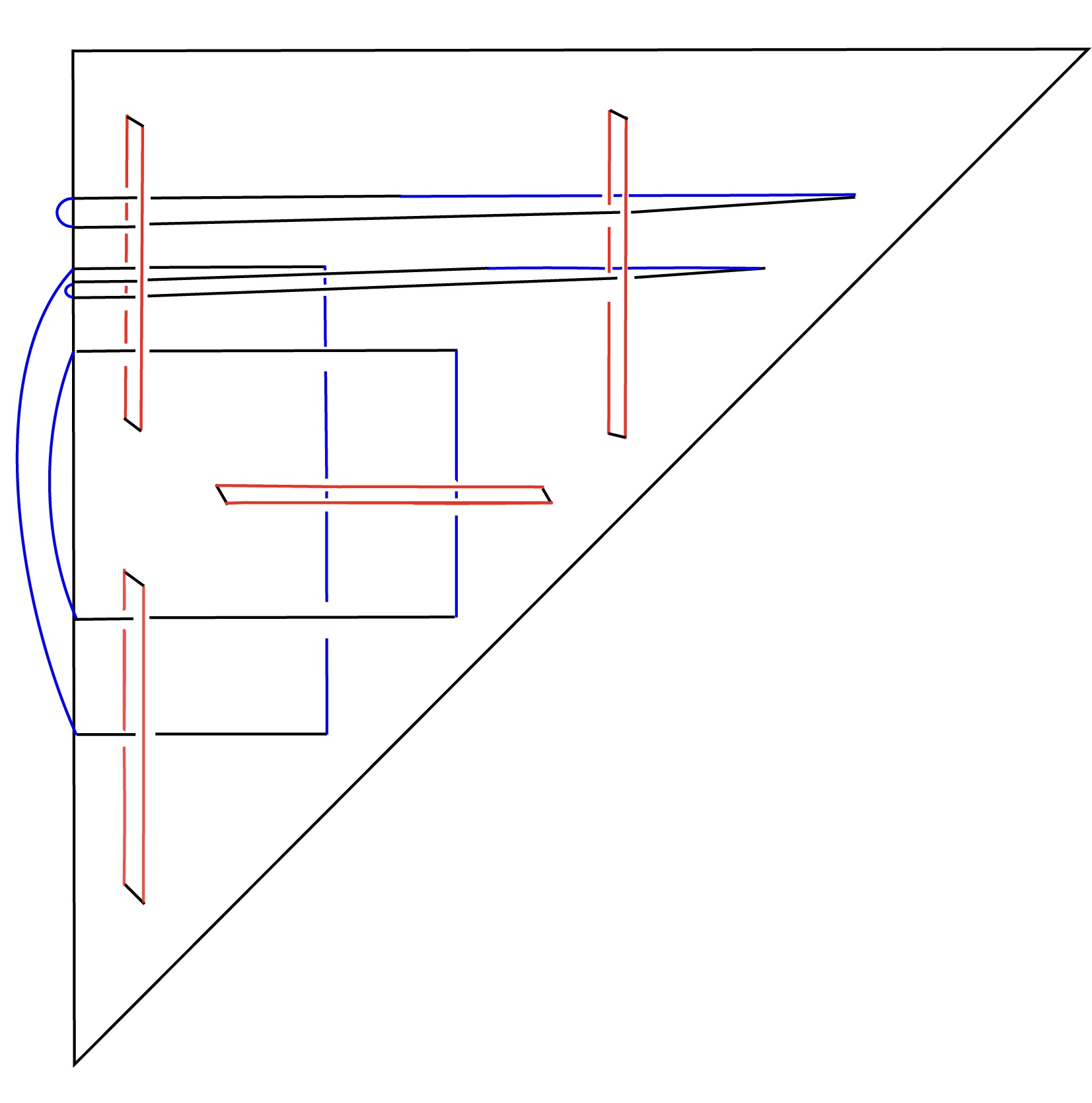}
    \caption{Type 4 (\textcolor{red}{$t^{\alpha-\beta}\mathrm{Co}_3^1$},\textcolor{blue}{$t^{\beta}\mathrm{Co}_2^3$}), contributes 
    $-t_1^{j+i}t_3^i-t_1^{j-i}t_3^{-i}$.}
      \label{calculate3type4}
  \end{minipage}
  \hfill
  \begin{minipage}[b]{0.44\textwidth}
    \includegraphics[width=\textwidth]{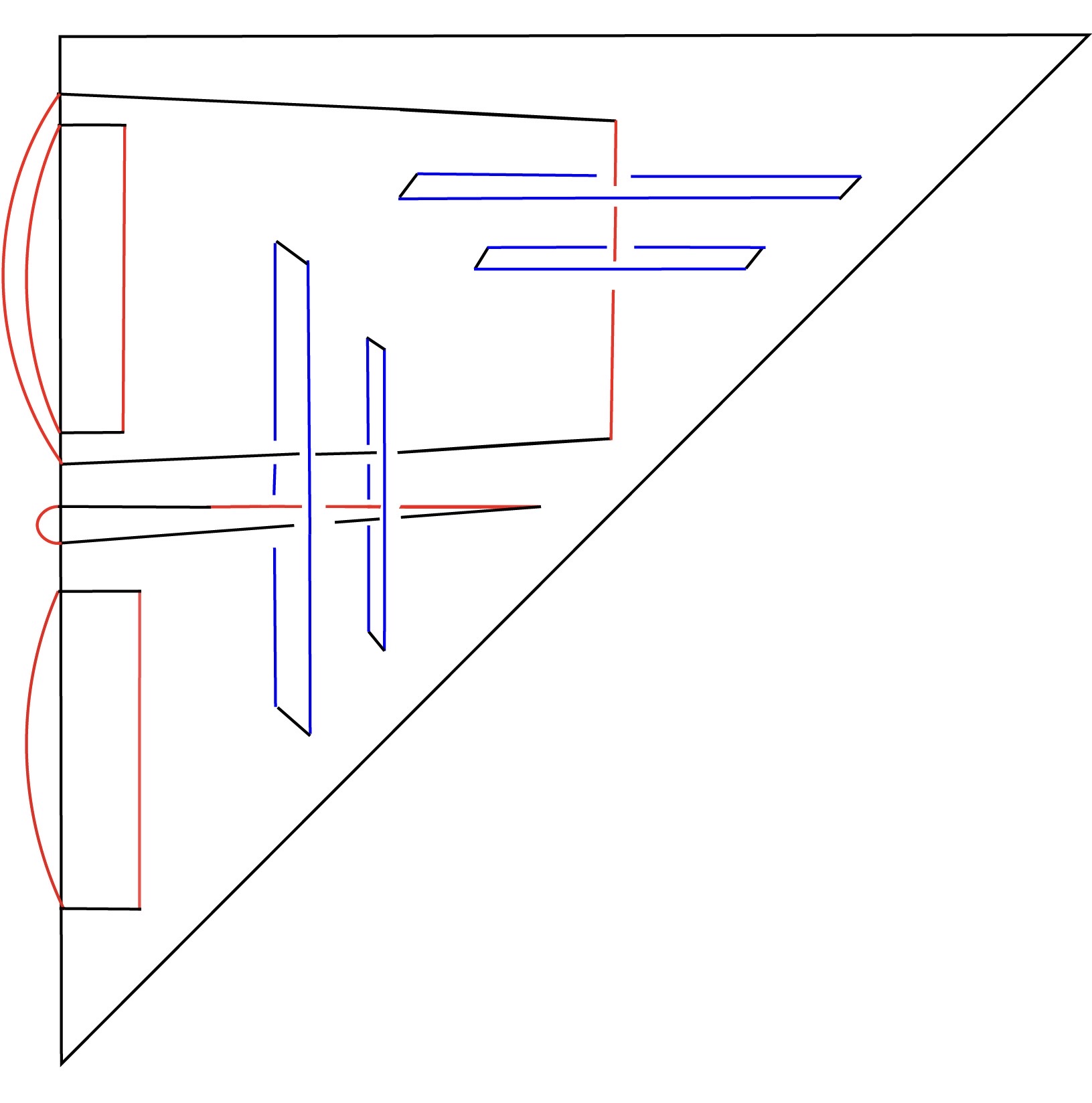}
    \caption{Type 4 (\textcolor{blue}{$t^{\alpha-\beta}\mathrm{Co}_3^1$},\textcolor{red}{$t^{\beta}\mathrm{Co}_2^3$}) contributes 0.}
      \label{calculate4type4}
  \end{minipage}
\end{figure}

\begin{figure}[!tbp]
  \centering
  \begin{minipage}[b]{0.46\textwidth}
    \includegraphics[width=\textwidth]{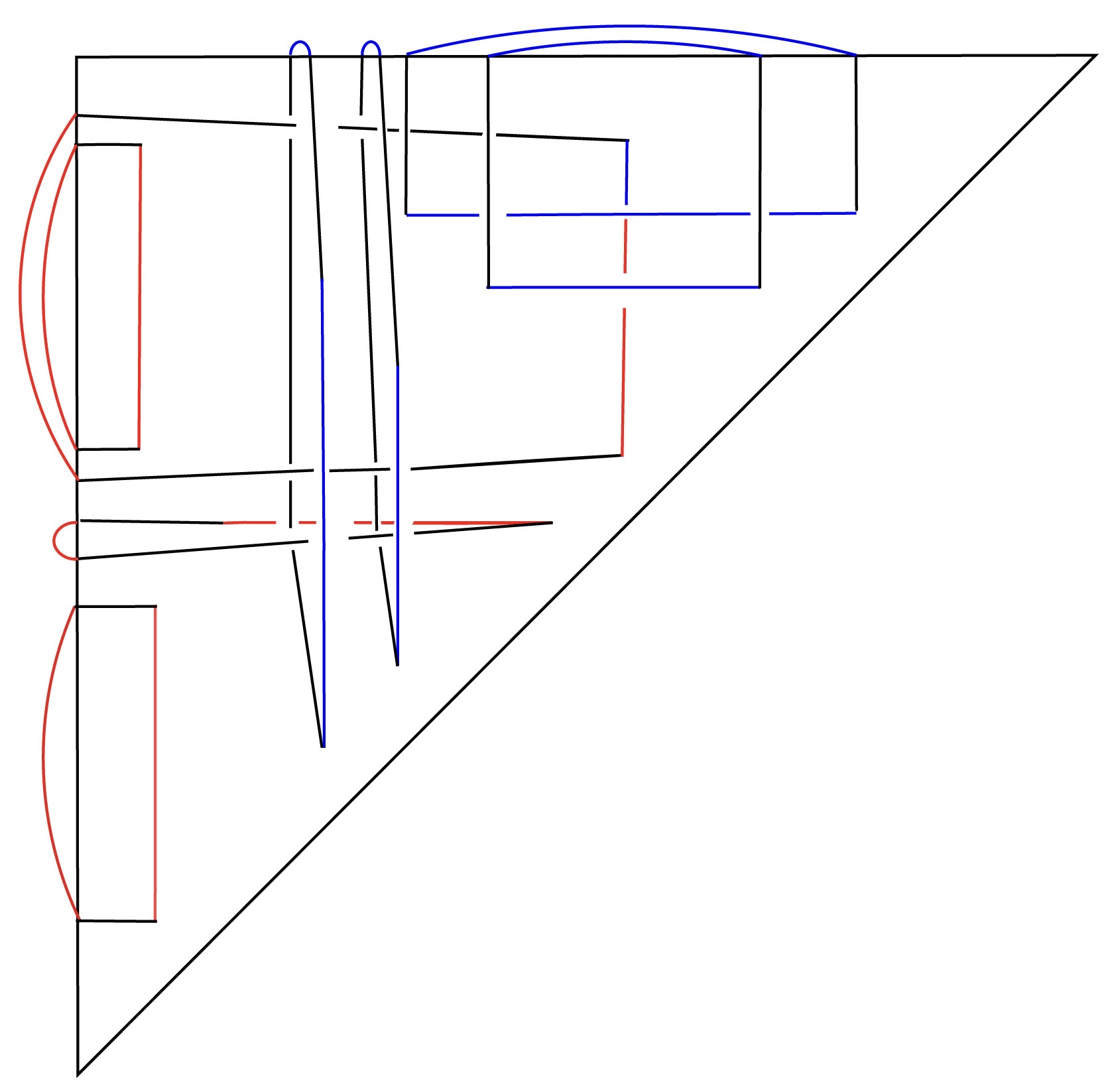}
    \caption{Type 4 (\textcolor{blue}{$t^\alpha\mathrm{Co}_2^1$},\textcolor{red}{$t^{\beta}\mathrm{Co}_2^3$}), contributes $t_1^it_3^j+t_1^{-i}t_3^j$.}
      \label{calculate5type4}
  \end{minipage}
  \hfill
  \begin{minipage}[b]{0.45\textwidth}
    \includegraphics[width=\textwidth]{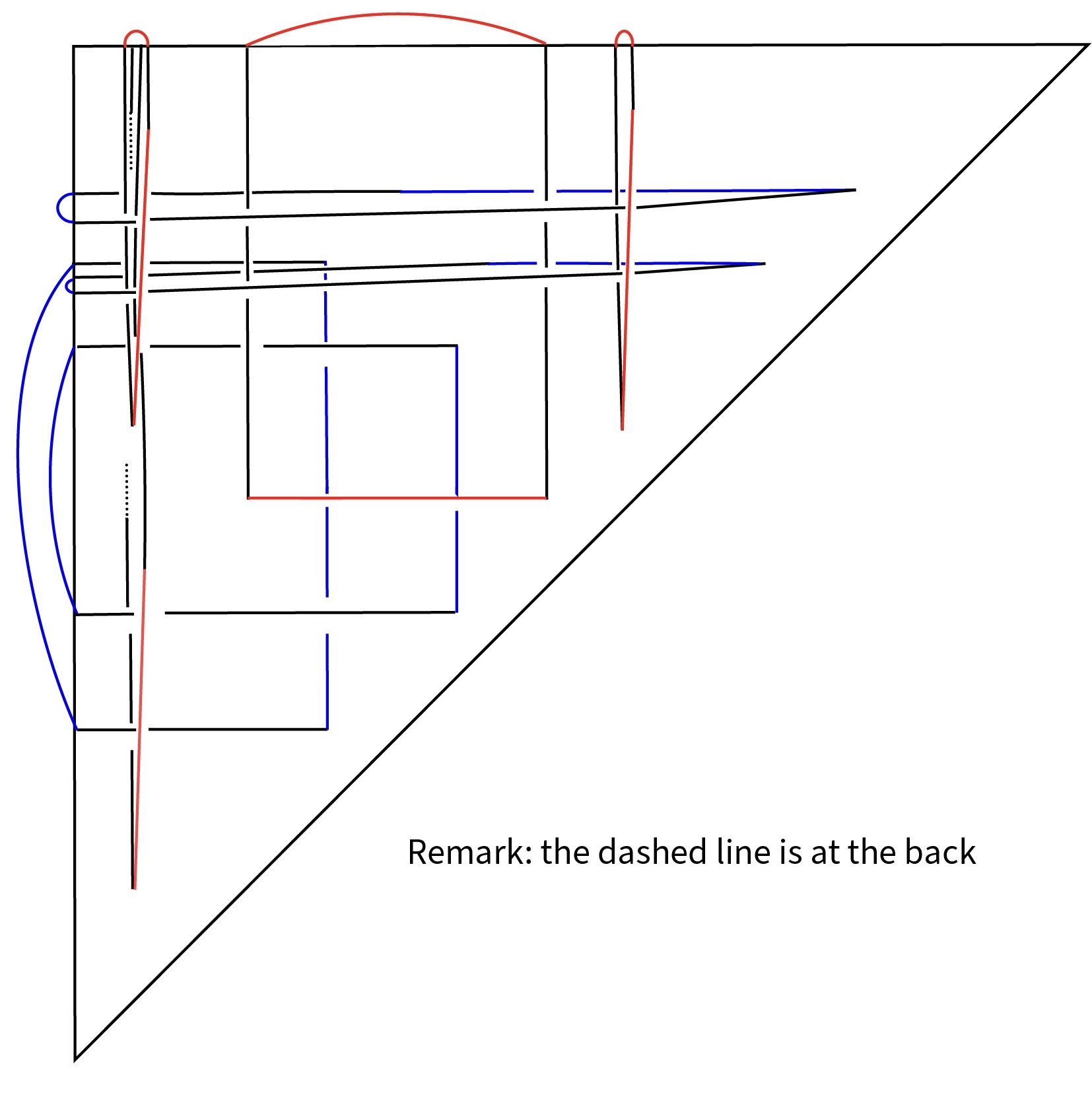}
    \caption{Type 4 (\textcolor{red}{$t^\alpha\mathrm{Co}_2^1$},\textcolor{blue}{$t^{\beta}\mathrm{Co}_2^3$}), contributes $-t_1^{j} t_3^{i}-t_1^{j}t_3^{-i}$.}
      \label{calculate6type4}
  \end{minipage}
\end{figure}

\begin{figure}[!tbp]
  \centering
  \begin{minipage}[b]{0.43\textwidth}
    \includegraphics[width=\textwidth]{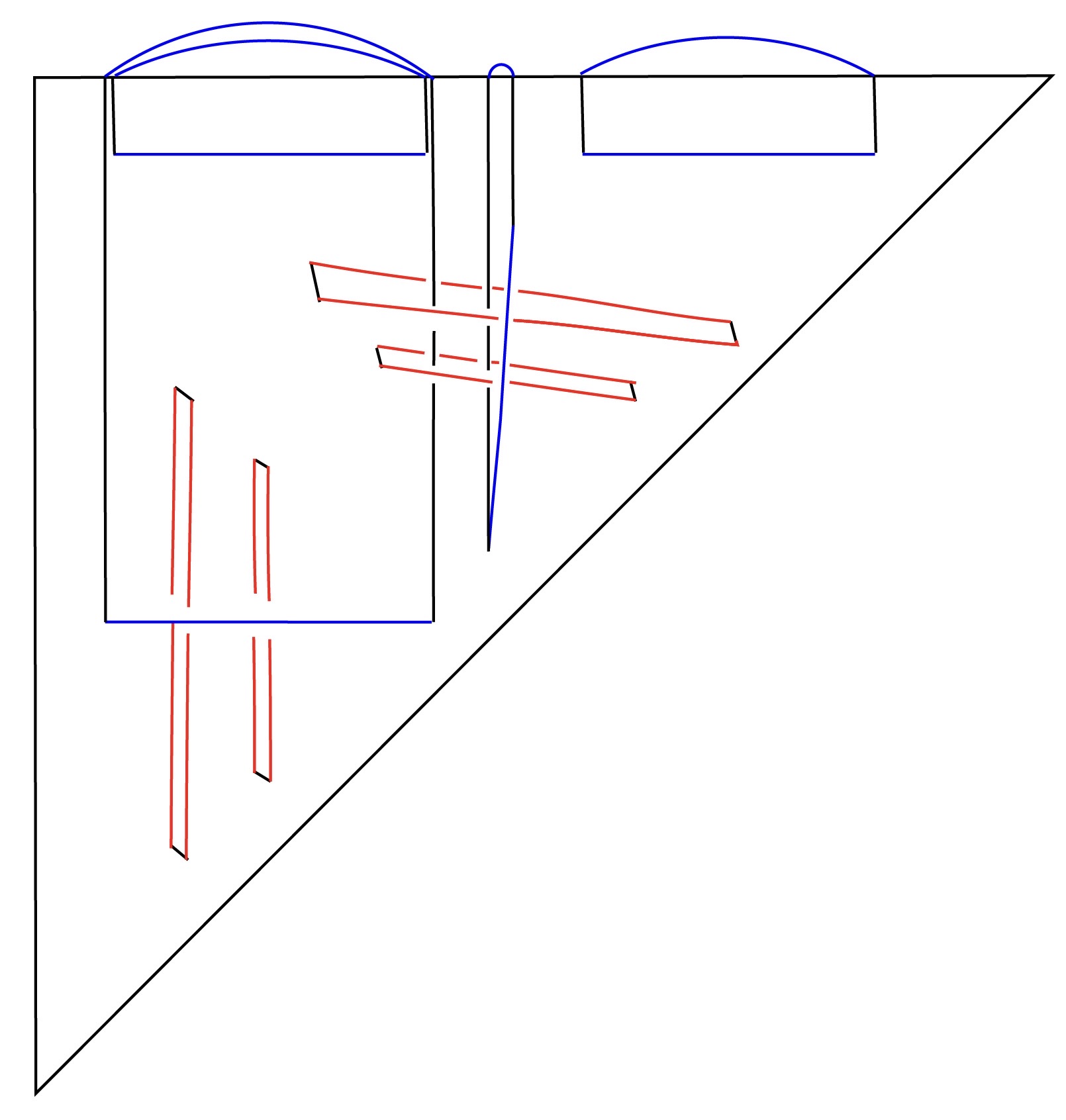}
    \caption{Type 6 (\textcolor{blue}{$t^\alpha\mathrm{Co}_2^1$},\textcolor{red}{$t^{\beta-\alpha}\mathrm{Co}_1^3$}), contributes $-t_1^{-i}t_3^{-j-i}-t_1^{-i}t_3^{j-i}+t_1^it_3^{i-j}+t_1^it_3^{j+i}$.}
      \label{calculate1type6}
  \end{minipage}
  \hfill
  \begin{minipage}[b]{0.45\textwidth}
    \includegraphics[width=\textwidth]{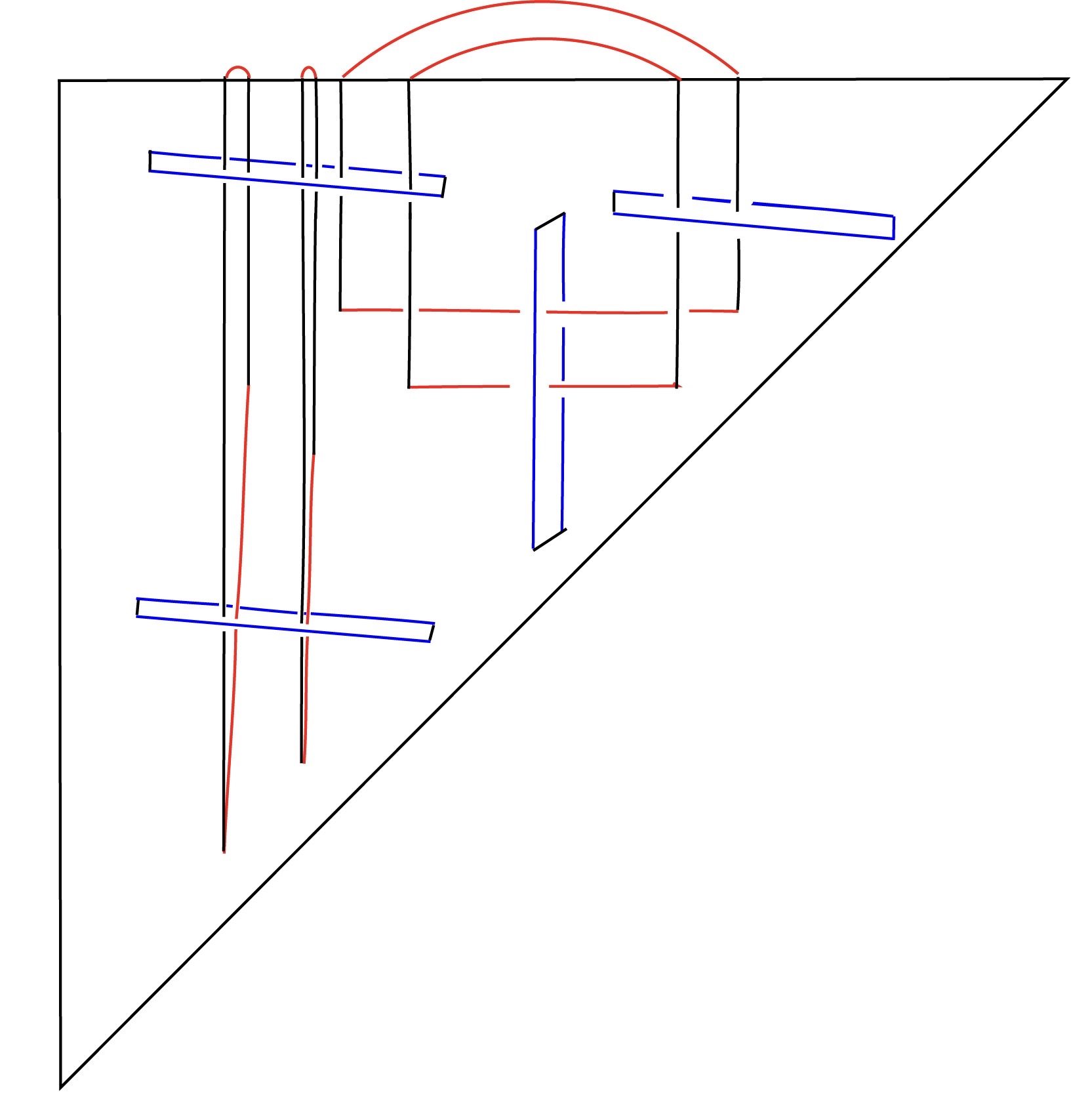}
    \caption{Type 6 (\textcolor{red}{$t^\alpha\mathrm{Co}_2^1$},\textcolor{blue}{$t^{\beta-\alpha}\mathrm{Co}_1^3$}), contributes $t_1^jt_3^{j+i}+t_1^{-j}t_3^{i-j}$.}
      \label{calculate2type6}
  \end{minipage}

\end{figure}

\begin{figure}[!tbp]
  \centering
  \begin{minipage}[b]{0.45\textwidth}
    \includegraphics[width=\textwidth]{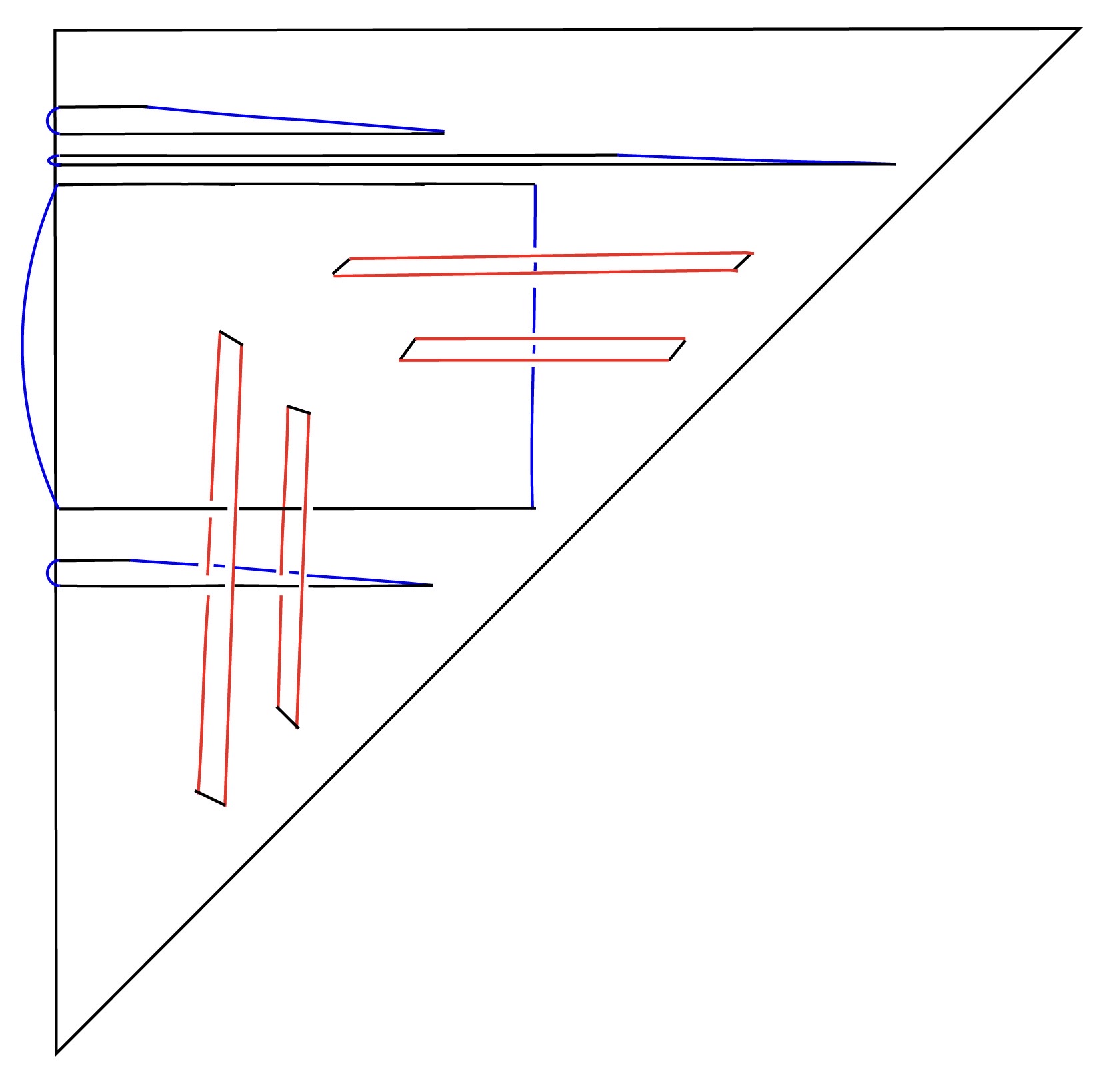}
    \caption{Type 6 (\textcolor{red}{$t^{\alpha-\beta}\mathrm{Co}_3^1$},\textcolor{blue}{$t^{\beta}\mathrm{Co}_2^3$}), contributes 
    $-t_1^{j+i}t_3^{i}-t_1^{i-j}t_3^i+t_1^{j-i}t_3^{-i}+t_1^{-j-i}t_3^{-i}$.}
      \label{calculate3type6}
  \end{minipage}
  \hfill
  \begin{minipage}[b]{0.46\textwidth}
    \includegraphics[width=\textwidth]{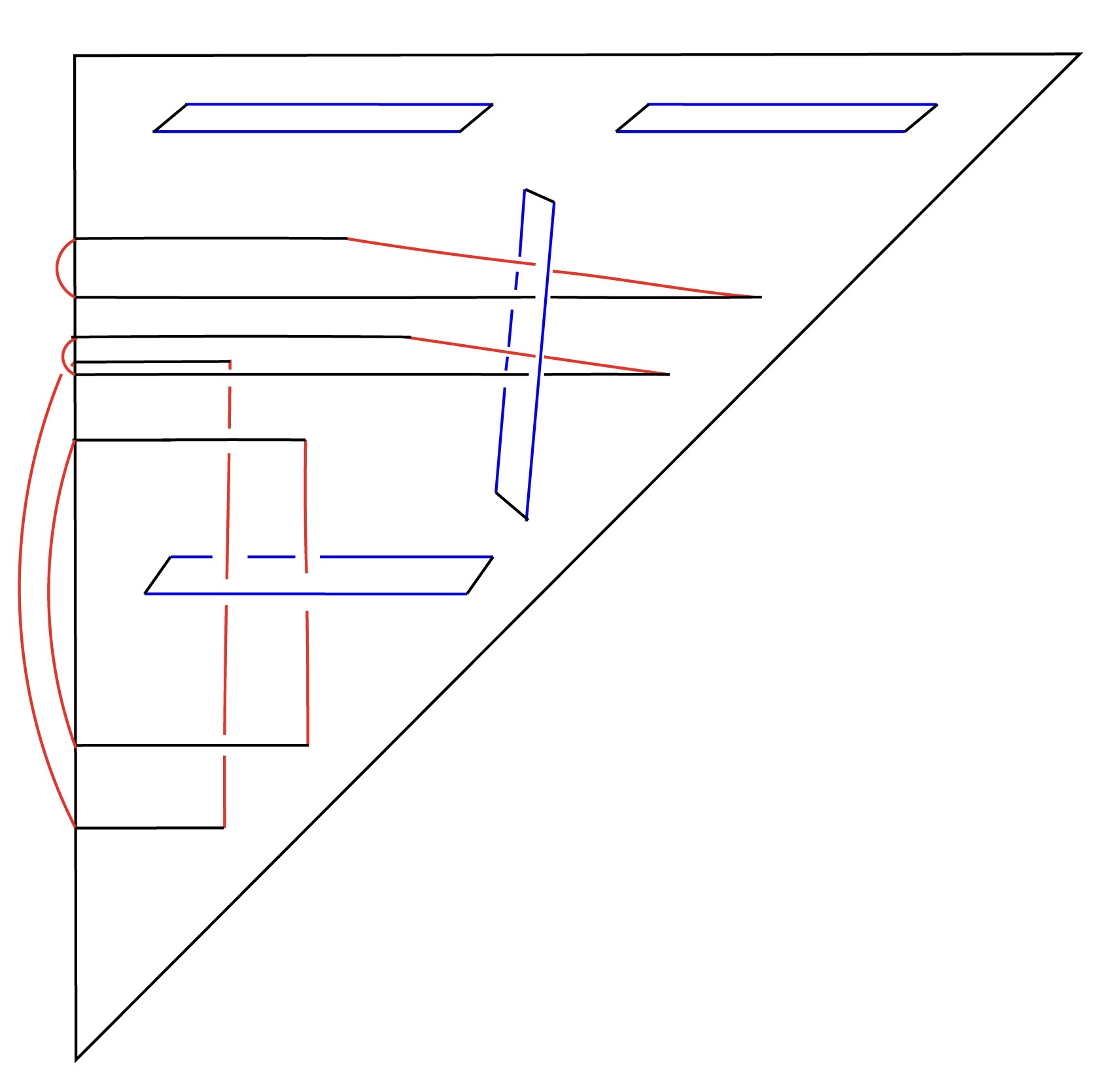}
    \caption{Type 6 (\textcolor{blue}{$t^{\alpha-\beta}\mathrm{Co}_3^1$},\textcolor{red}{$t^{\beta}\mathrm{Co}_2^3$}) contributes $-t_1^{j+i}t_3^j-t_1^{i-j}t_3^{-j}$.}
      \label{calculate4type6}
  \end{minipage}
\end{figure}

\begin{figure}[!tbp]
  \centering
  \begin{minipage}[b]{0.46\textwidth}
    \includegraphics[width=\textwidth]{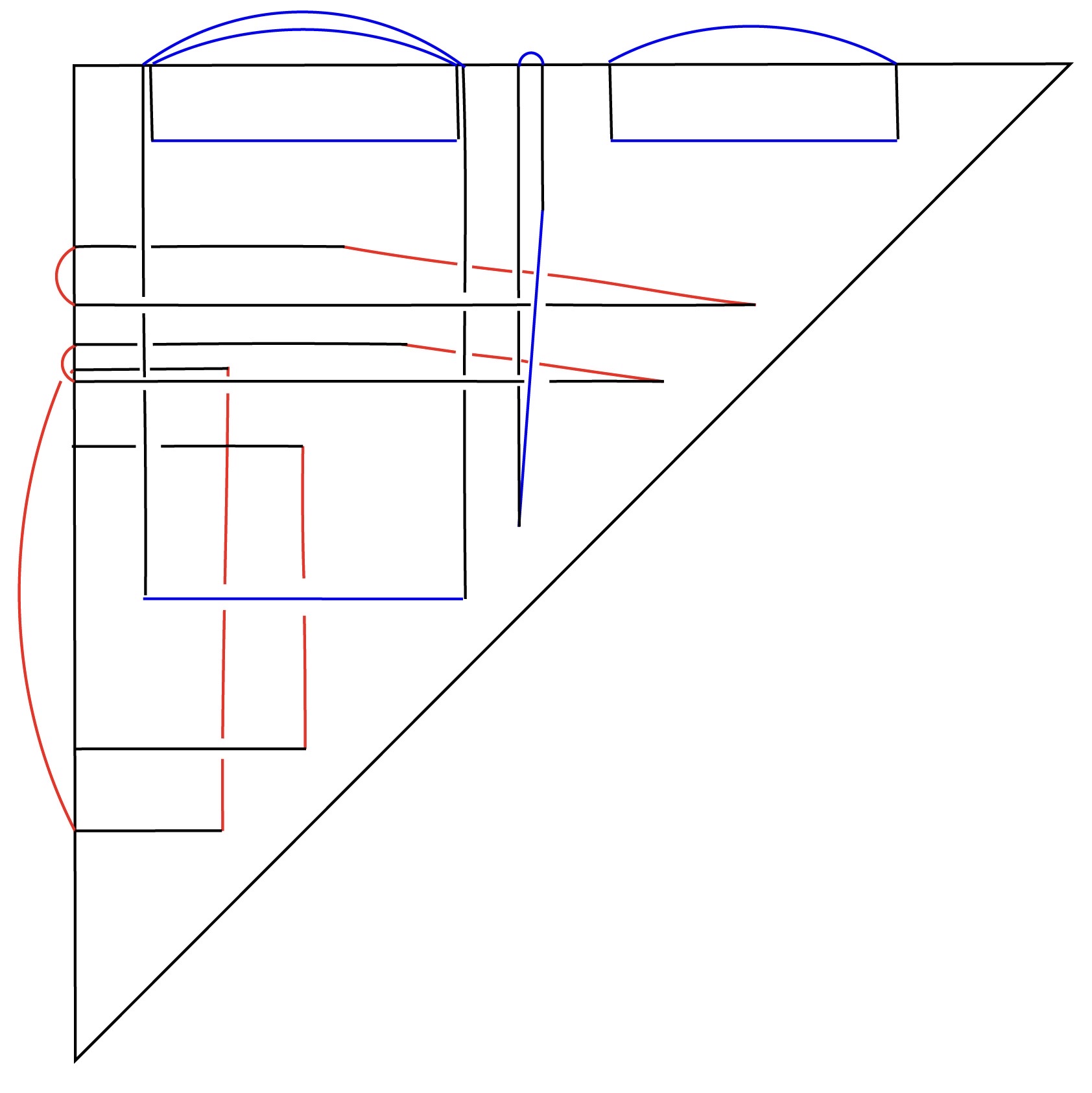}
    \caption{Type 6 (\textcolor{blue}{$t^\alpha\mathrm{Co}_2^1$},\textcolor{red}{$t^{\beta}\mathrm{Co}_2^3$}), contributes $-t_1^j t_3^i-t_1^{-j}t_3^i$.}
      \label{calculate5type6}
  \end{minipage}
  \hfill
  \begin{minipage}[b]{0.45\textwidth}
    \includegraphics[width=\textwidth]{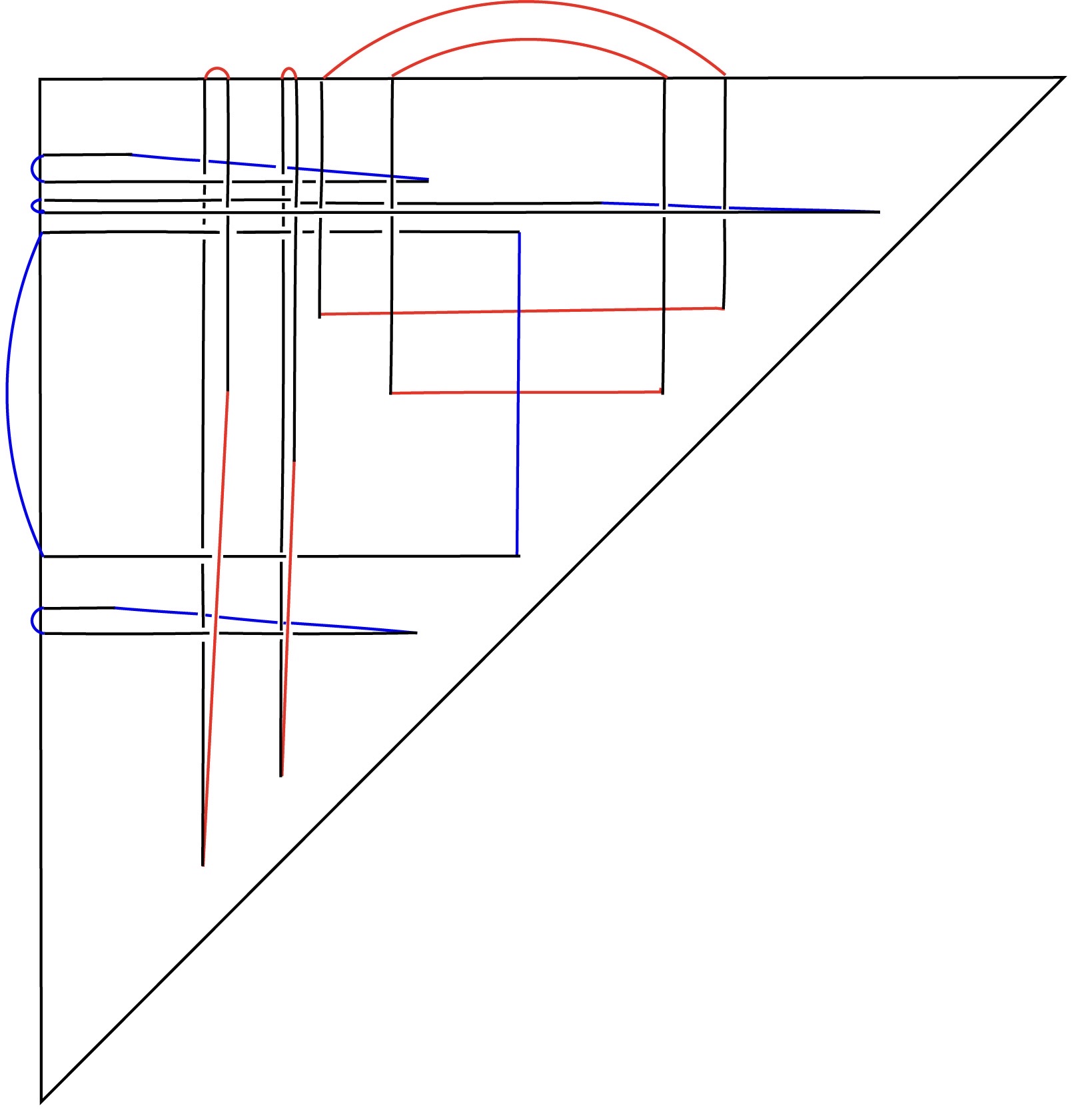}
    \caption{Type 6 (\textcolor{red}{$t^\alpha\mathrm{Co}_2^1$},\textcolor{blue}{$t^{\beta}\mathrm{Co}_2^3$}), contributes $t_1^{i} t_3^{j}+t_1^{i}t_3^{-j}$.}
      \label{calculate6type6}
  \end{minipage}
\end{figure}

As before, we present the relevant cohorizontal submanifolds based on their constituent parts, depending on which cuff the points are being mapped to (via colours) as drawn in Figures \ref{calculate1type4}, \ref{calculate2type4}, \ref{calculate3type4}, \ref{calculate4type4}, \ref{calculate5type4} and \ref{calculate6type4} for Type 4 intersection points and \ref{calculate1type6}, \ref{calculate2type6}, \ref{calculate3type6}, \ref{calculate4type6}, \ref{calculate5type6} and \ref{calculate6type6} for Type 6 intersection points, with the contributions by each constituent pair noteded in the figures as well. Taking the sum and applying Lemma \ref{type4and6barbell} leads to the following result.

\begin{theorem}
\label{newcasecomplete}
     For $i,j,k\in \mathbb{Z}^{+}$ with $i+j\leq k-1$, the polynomial obtained from a Type 4 intersection point of $\theta_{i+j+1}(e_i,e_j)$ is 
    $$\mathrm{T}_4(i,j)= t_1^it_3^{j+i}+t_1^{-i}t_3^{j-i}-t_1^{j+i}t_3^i-t_1^{j-i}t_3^{-i}+t_1^it_3^j+t_1^{-i}t_3^j-t_1^{j}t_3^{i}-t_1^{j}t_3^{-i}$$
    and the polynomial obtained from a Type 6 intersection point is 
   $$\mathrm{T}_6(i,j)= -t_1^{-i}t_3^{-j-i}-t_1^{-i}t_3^{j-i}+t_1^it_3^{i-j}+t_1^it_3^{j+i}+t_1^jt_3^{j+i}+t_i^{-j}t_3^{i-j}-t_1^{j+i}t_3^{i}-t_1^{i-j}t_3^i+t_1^{j-i}t_3^{-i}$$
$$+t_1^{-j-i}t_3^{-i}-t_1^{j+i}t_3^j-t_1^{i-j}t_3^{-j}-t_1^j t_3^i-t_1^{-j}t_3^i+t_1^{i} t_3^{j}+t_1^{i}t_3^{-j}.$$
In particular, we have \begin{align*}
        & W_3(\theta_{i+j+1}(e_i,e_j))=
i(t_1^it_3^{j+i}+t_1^{-i}t_3^{j-i}-t_1^{j+i}t_3^i-t_1^{j-i}t_3^{-i}+t_1^it_3^j+t_1^{-i}t_3^j-t_1^{j}t_3^{i}-t_1^{j}t_3^{-i})\\
&+j(-t_1^{-i}t_3^{-j-i}-t_1^{-i}t_3^{j-i}+t_1^it_3^{i-j}+t_1^it_3^{j+i}+t_1^jt_3^{j+i}+t_i^{-j}t_3^{i-j}-t_1^{j+i}t_3^{i}-t_1^{i-j}t_3^i+t_1^{j-i}t_3^{-i}
+\\
&t_1^{-j-i}t_3^{-i}-t_1^{j+i}t_3^j-
t_1^{i-j}t_3^{-j}-t_1^j t_3^i-t_1^{-j}t_3^i+t_1^{i} t_3^{j}+t_1^{i}t_3^{-j}).
    \end{align*} 
\end{theorem}

We mentioned in the beginning of Section 4.2 that Budney--Gabai discovered an alternative method of calculating the $W_3$ invariant for $\theta_k(e_i,e_j)$ for all possible $i$ and $j$ in \cite{Budney-gabai}. In particular, they produced a formula for $\theta_k(e_i,e_j)$ in terms of their fundamental classes $G(p,q)\in \pi_2\mathrm{Emb}(I,S^1\times D^3)$. For $i+j\geq k$, 
\begin{multline*}
    [\theta_k(e_i,e_j)]=(k-i-1)(D(-j,i)-D(i-j,-i))
    +(k-j-1)(D(i,-j)-D(i-j,j))
+(i+j+1-k)D(i,-j)
\end{multline*}
and for $i+j\leq k-1$,
\begin{multline*}
     [\theta_k(e_i,e_j)]=i(D(i,-j)-D(i+j,-j))-D(i-j,j)+D(i,j)
    +j(D(-j,i)-D(-i-j,i)-D(i-j,-i)+D(i,j))
\end{multline*}
where $D(i,j)=-G(j,-i)+G(-j,i)-G(i,-j)+G(-i,j)$
and they calculated that $W_3(G(i,j))=t_1^{i-j}t_3^{-j}$. See Example 3.6 of \cite{Budney-gabai}. One can verify that our formulae in Theorems \ref{mainformula} and \ref{newcasecomplete} coincide with their formulae. For example, by the Hexagon relation, \begin{align*}
  &  W_3(D(-j,i)-D(i-j,-i))  = (-t_1^{i-j}t_3^{-j}+t_1^{j-i}t_3^{j}-t_1^{i-j}t_3^i+t_1^{j-i}t_3^{-i})-(-t_1^{-j}t_3^{i-j}+t_1^jt_3^{j-i}-t_1^{-j}t_3^{-i}+t_1^{j}t_3^i)\\
&   = t_1^i t_3^{j}-t_1^{-j} t_3^{-i}+t_1^{-i}t_3^{j-i}+t_1^{i-j}t_3^{j}=-t_1^{j}t_3^{i}-t_1^{i}t_3^{j}+t_1^{-i}t_3^{j-i}+t_1^{j-i}t_3^{-i}
\end{align*} 
which indicates that the two formulae for a Type 1 intersection point coincide. The remaining cases can be checked in a similar way. We will see in Section 4.4 that our formula has the advantage that it provides direct insights of calculating $W_3$ for general unknotted barbells.

We now discuss some properties of $W_3(\theta_k(v,w))$. The following lemma is a restatement of Corollary 7.18 in \cite{Budney-gabai} (which is phrased in terms of certain fundamental classes $G(p,q)$ mentioned above) but in a different context and can be checked using Theorems \ref{mainformula} and \ref{newcasecomplete} together with the Hexagon relation.

\begin{lemma}
\label{antisym}
    For $k> 0$, $$W_3(\theta_k(e_i,e_j))=-W_3(\theta_k(e_j,e_i)).$$
\end{lemma}
\begin{proof}
This is purely by plugging in the formulae. We verify the case $i+j\geq k$ and the remaining case is similar. By Theorem \ref{mainformula}, we have
    \begin{align*}
       & W_3(\theta_k(e_i,e_j))=(k-i-1)(-t_1^{j}t_3^{i}-t_1^{i}t_3^{j}+t_1^{-i}t_3^{j-i}+t_1^{j-i}t_3^{-i})+\\
       &(-k+j+i-1)(-t_1^{-j}t_3^{i-j}+
       t_1^{-i}t_3^{j-i}+t_1^{j-i}t_3^{-i}-t_1^{i-j}t_3^{-j})+\\
       &(k-j-1)(-t_1^{-j}t_3^{-i}-t_1^{-i}t_3^{-j}-t_1^{-j}t_3^{i-j}+
       t_1^{-i}t_3^{j-i}+
       t_1^{i}t_3^{i-j}+t_1^{i-j}t_3^{i}+t_1^{j-i}t_3^{-i}-t_1^{i-j}t_3^{-j})
    \end{align*} coming from the three types of intersection points. Note that the contribution of a Type 2 intersection point is itself anti-symmetric. We show that the contributions of Type 1 intersection points of $W_3(\theta_k(e_i,e_j))$ cancel with the contributions of Type 3 intersection points of $W_3(\theta_k(e_j,e_i))$:
    \begin{align*}
        &(k-i-1)(-t_1^{j}t_3^{i}-t_1^{i}t_3^{j}+t_1^{-i}t_3^{j-i}+t_1^{j-i}t_3^{-i})+\\
      &  (k-i-1)(-t_1^{-i}t_3^{-j}-t_1^{-j}t_3^{-i}-t_1^{-i}t_3^{j-i}+t_1^{-j}t_3^{i-j}+t_1^{j}t_3^{j-i}+t_1^{j-i}t_3^{j}+t_1^{i-j}t_3^{-j}-t_1^{j-i}t_3^{-i})\\
  &      =(k-i-1)(-(t_1^{-i}t_3^{-j}+t_1^jt_3^i)-(t_1^{-j}t_3^{-i}+t_1^{i}t_3^j)+(t_1^{-j}t_3^{i-j}+t_1^{j-i}t_3^{j})+(t_1^j t_3^{j-i}+t_1^{i-j}t_3^{-j}))\\
   &     =(k-i-1)(-(t_1^{-i}t_3^{-j}+t_1^jt_3^i)-(r_1^{i-j}t_3^{-j}+t_1^jt_3^{j-i})+(t_1^{-j}t_3^{i-j}+t_1^{j-i}t_3^{j})+(t_1^j t_3^{j-i}+t_1^{i-j}t_3^{-j}))=0
\end{align*}
Similarly, the contributions of Type 3 intersection points of $W_3(\theta_k(e_i,e_j))$ cancels with the contributions of Type 1 intersection points of $W_3(\theta_k(e_j,e_i))$.
\end{proof}

\begin{remark}
    In fact, there is a direct way of proving Lemma \ref{antisym} for all $1\leq i,j\leq k-1$. We shrink $\theta_k(e_i,e_j)$ such that it is embedded in a smaller (4-dimensional) solid torus in $S^1\times D^3$. Let $\psi$ denote the rotation of this smaller solid torus in $S^1\times D^3$ by $\pi$. See Figure \ref{newargument3} for a 3-dimensional picture. Now if we conjugate $\Phi_{\theta_k(e_i,e_j)}$ by $\psi$, the claim is we get $\Phi_{\theta_k(e_j,e_i)}^{-1}$. In other words, we have $$\psi \Phi_{\theta_k(e_j,e_i)}^{-1} \psi^{-1} = \Phi_{\theta_k(e_i,e_j)}.$$
To see this, note that the image $\psi(\theta_k(e_i,e_j))$ is equal to $\theta_k(e_j,e_i)$ with the role of the two cuff spheres switched, and this switch corresponds to taking the inverse of the induced barbell diffeomorphism. But now since $\pi_0\mathrm{Diff}(S^1\times D^3,\partial)$ is abelian, we get exactly the antisymmetric relation we want.
    \begin{figure}
        \centering
        \includegraphics[width=0.4\linewidth]{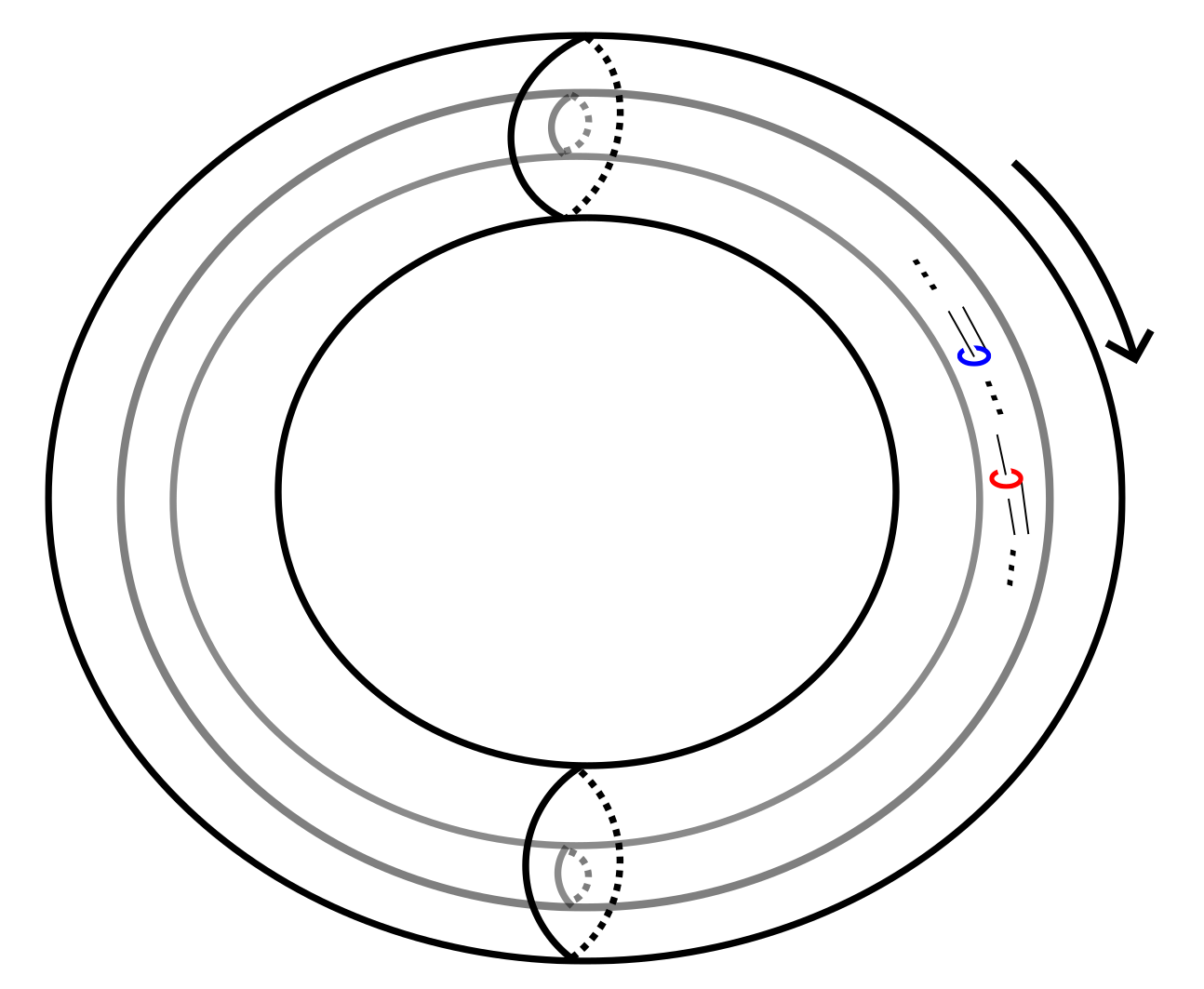}
        \caption{Anti-symmetric relation of $\theta_k(e_i,e_j)$.}
        \label{newargument3}
    \end{figure}
\end{remark}

The following lemma is a key observation for calculating the $W_3$ invariant of $\theta_k(v,w)$.

\begin{lemma}
\label{linear}
    The $W_3$ invariant of $\theta_k(v,w)$ satisfies $$W_3(\theta_k(v,w))=\Sigma v_i w_j W_3(\theta_k(e_i,e_j))$$
\end{lemma}
\begin{proof}

The idea is that we can analyze each set of cohorizontal points that arises from one of the possible pairs of vertical strands (one going through the red cuff and one going through the blue cuff) separately and take the sum. We now elaborate on this. By isotoping $\theta_k(v,w)$ if necessary, we saw that the $W_3$ invariant is determined by the $k-1$ intersections of the bar with the 3-ball $s_0\times D^3$. First, note that 
$$W_3(\theta_k(ae_i,e_j))=aW_3(\theta_k(e_i,e_j))$$
for $a\in \mathbb{Z}$ where $a$ indicates that the blue cuff $B$ wraps around the $i$-th strand $a$ times, counting from right to left. This is because $\theta_k(ae_i,e_j)$ produces $a$ parallel copies of the $8$ cohorizontal spheres, counting with a sign rather than only $8$ cohorizontal spheres, near the blue cuff sphere for $\theta_k(e_i,e_j)$. Similarly, we have
$$W_3(\theta_k(e_i,be_j))=bW_3(\theta_k(e_i,e_j)).$$

Next, if for $1\leq i,j\leq k$, we let $e_{ij}$ be the vector with all entries 0 except the $i$-th and $j$-th place being 1, then we have 
$$W_3(\theta_k(e_{ij},e_l))=W_3(\theta_k(e_i,e_l))+W_3(\theta_k(e_j,e_l)).$$
To see this, we observe that for each of the $k-1$ intersection points, $\theta_k(e_{ij},e_l)$ has blue cohorizontal points near both the $i$-th strand and the $j$-th strand, counting from right to left. Each of these two sets of cohorizontal points pairs with the red cohorizontal points (near the $l$-th strand, counting from left to right) separately, precisely leading to the cohorizontal points of $\theta_k(e_i,e_l)$ and $\theta_k(e_j,e_l)$. Similarly, we have 
$$W_3(\theta_k(e_{l},e_{ij}))=W_3(\theta_k(e_l,e_i))+W_3(\theta_k(e_l,e_j)).$$
More generally, we observe that 
$$W_3(\theta_k(ae_i+be_j,ce_k+de_l)=acW_3(e_i,e_k)+adW_3(e_i,e_l)+bcW_3(e_j,e_k)+bdW_3(e_j,e_l).$$
for $a,b,c,d\in \mathbb{Z}$ since the cohorizontal points that arise from this barbell can be analysed separately in pairs in a similar way, leading to above combination. Now, the theorem follows from the combination of the above.\end{proof}

In the proof of Lemma \ref{linear}, we break up the $W_3$ invariant of $\theta_k(v,w)$ into a sum of the $W_3$ of a set of ``sub-barbells'' $\theta_k(e_i,e_j)$ (where $i=j$ is allowed), each of which gives a parallel copy of a set of tetrahedron pictures. Essentially, Budney--Gabai proved the same result in a slightly different context (by directly homotoping $s_{s_0}(\theta_k(v,w))$ into a sum of their fundamental classes). See Theorem 8.1 of \cite{Budney-gabai}. 

We end this section (and chapter) with a discussion of the question of which barbells give rise to linearly independent barbell diffeomorphisms. First, note that $\pi_0\mathrm{Diff}(S^1\times D^3,\partial)$ is abelian. This can be seen by noting that we can shrink the supports of any two diffeomorphisms to two smaller 4-dimensional solid tori to make them commute. Therefore, it makes sense to talk about linear independence of diffeomorphisms (up to isotopy) in this group. Budney--Gabai (using Lemma \ref{deltabudney}) proved that $W_3(\delta_k)=W_3(\theta_k(e_{k-1},e_{k-2}))$ for $k\geq 4$ are linearly independent. Thus $\delta_k=\theta_k(e_{k-1},e_{k-2})$ for $k\geq 4$ are linearly independent diffeomorphisms of $S^1\times D^3$.

For $k\geq 6$, we can apply Theorem \ref{mainformula} to calculate $W_3(\theta_k(e_{k-1},e_{k-3}))$. By setting all terms $t_1^pt_3^q=0$ except when $p-q=k-1$, we get $$-(k-1)t_1^2t_3^{3-k}$$ which is nonzero. Furthermore, if a finite combination of $\theta_n(e_{n-1},e_{n-3})$, for $1\leq n\leq k$, (i.e.~$k$ is the largest $n$)
$$a_6\theta_6(e_{5},e_{3})+a_7 \theta_7(e_{6},e_{4})+\dots+ a_k \theta_k(e_{k-1},e_{k-3})=0$$
equals 0, then setting all terms $t_1^pt_3^q$ except when $p-q=k-1$ kills $\theta_n(e_{n-1},e_{n-3})$ with $n<k$. This implies that the coefficient $a_k$ before $\theta_k(e_{k-1},e_{k-3})$ is 0. Therefore, we deduce that the set $\{\theta_k(e_{k-1},e_{k-3})\}_{k\geq 6}$ is linearly independent. Furthermore, by Corollary 8.2 of \cite{Budney-gabai}, the same argument applied to $\delta_k=\theta_k(e_{k-1},e_{k-2})$ gives rise to $$-(k-1)t_1^1t_3^{2-k}\neq 0$$
which implies that the set $\{\delta_k,\theta_k(e_{k-1},e_{k-3})\}_{k\geq 6}\}$ is linearly independent.

By applying Theorem \ref{mainformula} again to $W_3(\theta_k(e_{k-1},e_{k-m}))$ for $k>m+1$, and following the same argument above, namely by setting all terms $t_1^pt_3^q$ in $W_3(\theta_k(e_{k-1},e_{k-m}))$ to zero except when $p-q=k-1$, one deduces that the set of barbells $$\{\theta_k(e_{k-1},e_{k-m})\}_{{m\in \{3,\dots,
 [(k-1)/2]-1\}},k\geq 2m-1}$$ is linearly independent. Here $[(k-1)/2]$ denotes the integer part of $(k-1)/2$. Therefore, we have the following theorem.
 \begin{theorem}
\label{chapter4main}
        The elements $\theta_k(e_{k-1},e_{k-3})$ for $k\geq 6$ of $\pi_0\mathrm{Diff}(S^1\times D^3,\partial)$ are linearly independent. Further, these elements are linearly independent to 
        $$\delta_k=\theta_k(e_{k-1},e_{k-2})=\theta_k((0,\dots,0,1),(0,\dots,0,1,0))$$
        for $k\geq 4$.
        
        More generally, there exist linearly independent elements $\theta_k(e_{k-1},e_{k-m})$ of $\pi_0\mathrm{Diff}(S^1\times D^3,\partial)$ for $m\in \{3,4,\dots,[(k-1)/2]-1\}$ with $k\geq 2m-1$. Here $[(k-1)/2]$ is the integer part of $(k-1)/2$.
    \end{theorem}

We can also consider $\theta _k (e_{k-2},e_{k-m})$ with  $k\geq 4$, $m\in \{3,4,\dots,k-2,k-1\}$ and $k-m>2$, and apply Theorem \ref{mainformula}. Setting all terms of $W_3(\theta _k (e_{k-2},e_{k-m}))$ in the form of $t_1^pt_3^q$ to zero except when $p-q=k-2$ gives rise to $$-(k-2)t_1^{m-2}t_3^{m-k}$$ which is non-zero if $k-m\neq m-2$, i.e.~$2m\neq k+2$. By the Hexagon relation together with Theorem \ref{mainformula},
 $$W_3(\theta_k(e_{k-2},e_{k-m})+\theta_k(e_{k-2
},e_{m-2}))=0$$ in this quotient. It follows that for a fixed $k\geq 4$, there are $(k-3)/2$ linearly independent terms for $k$ odd and $(k-2)/2$ linearly independent terms for $k$ even. For example, when $k=10$, the linearly independent terms we get are $$\theta_{10}(e_8,e_7),\theta_{10}(e_8,e_6),\theta_{10}(e_8,e_5), \theta_{10}(e_8,e_4).$$

\begin{Question}
    More generally, for a fixed $\alpha, \beta \geq 1$, $m>\alpha$ and $n> \beta$, one can ask if the barbells $\theta_{k_1}(e_{k_1-\alpha},e_{k_1-m})$ and $\theta_{k_2}(e_{k_2-\beta},e_{k_2-n})$ are linearly independent. 
    
    One possible approach is to set all terms $t_1^pt_3^q$ to zero except when $p-q$ equals to $m-\alpha$ or $n-\beta$. If we assume $n-\beta\neq \alpha-m$ or $k_1-\alpha$ or $k_1-m$, and $m-\alpha \neq k_2-\beta$ or $k_2-n$, then we get two polynomials $-(k_1-\alpha)t_1^{m-\alpha}t_3^{m-k_1}$ and $-(k_2-\beta)t_1^{n-\beta}t_3^{n-k_2}$ which are linearly dependent if and only if $m-\alpha=k_2-n$ and $m-k_1=\beta-n$. If $\alpha =\beta$, the above conditions imply that $k_1=k_2$ and $m+n=k_1+\alpha$. Similarly, if $k_1=k_2$, then the above conditions imply that $\alpha=\beta$ and $m+n=k_1+\alpha$. 
\end{Question}

\subsection{Computing $W_3$ for unknotted barbells in $S^1\times D^3$}

The following theorem indicates how the $(W_3)_1$ of general unknotted barbells in $S^1\times D^3$ is calculated.

\begin{theorem}
\label{barbellfactor}
    Let $\mathcal{B}$ be an unknotted barbell in $S^1\times D^3$ represented by an element $$ \nu_B^{i_1}t^{j_1}\nu_R^{k_1}\cdots \nu_B^{i_n}t^{j_n}\nu_R^{k_n}$$ 
    in $\langle \nu_R\rangle \backslash  F_3/\langle \nu_B\rangle$, where $n\in \mathbb{Z}$. Define $|\mathcal{B}|=j_1+\dots+j_n$, i.e.~it is the total signed number of times the bar of $\mathcal{B}$ loops around the circle factor of $S^1\times D^3$. Then $W_3(\mathcal{B})$ is equal to the sum
    $$\sum_{\alpha,\beta} c_{i_{\alpha} k_{\beta}}W_3(\mathcal{B}_{i_{\alpha} k_{\beta}})$$
    of $W_3$ invariants of a finite sequence of ``sub-barbells'' $\mathcal{B}_{i_{\alpha} k_{\beta}}$ with $\alpha, \beta \in \{1,2,\dots,n\}$, and coefficients $c_{i_{\alpha} k_{\beta}}$. Each $\mathcal{B}_{i_{\alpha} k_{\beta}}$ is an unknotted barbell which satisfies the condition that the bar links each of the two cuff spheres only positively once, and is represented by a word in the form
    $$t^{j_1}t^{j_2}\dots t^{j_{\alpha-1}}\nu_B^{1}t^{\alpha}\dots t^{\beta}\nu_R^{ 1}\dots t^{j_n}$$
    when $\alpha\leq \beta$, and in the form
    $$t^{j_1}t^{j_2}\dots t^{j_{\beta}}\nu_R^{1}t^{\beta+1}\dots t^{\alpha-1}\nu_B^{1}\dots t^{j_n}$$
    when $\alpha>\beta$. Also,  $|\mathcal{B}_{i_{\alpha}k_{\beta}}|=|\mathcal{B}|$ for all $\alpha,\beta$. Further, the coefficients satisfy
    $$c_{i_{\alpha} k_{\beta}}=i_{\alpha} k_{\beta}.$$
\end{theorem}

\begin{proof}
    The idea is the same as the proof of Lemma \ref{linear}. We first break the information of cohorizontal points from $\mathcal{B}$ into barbells in the form of 
    $$ \cdots t^{j_1}\cdot \nu_B^{i_{\alpha}}t^{j_{\alpha}}\cdots\nu_R^{k_{\beta}}\cdots t^{j_n},$$ 
    i.e.~barbells that keep only a single power of $\nu_B$ and a single power of $\nu_R$ but keep all powers of $t$. Each such barbell leads to triangle and tetrahedron diagrams that contain a piece of information of cohorizontal points of $\mathcal{B}$. Next, we further observe that the exponents $i_{\alpha}$ and $k_{\beta}$ have the effects of contributing to parallel copies of cohorizontal points of
    $$ \cdots t^{j_1}\cdot \nu_B^{1}t^{j_{\alpha}}\cdots\nu_R^{ 1}\cdots t^{j_n}$$
leading to coefficients $c_{i_{\alpha} k_{\beta}}=i_{\alpha}k_{\beta}$.\end{proof}

    \begin{figure}
    \centering
    \includegraphics[width=0.4\textwidth]{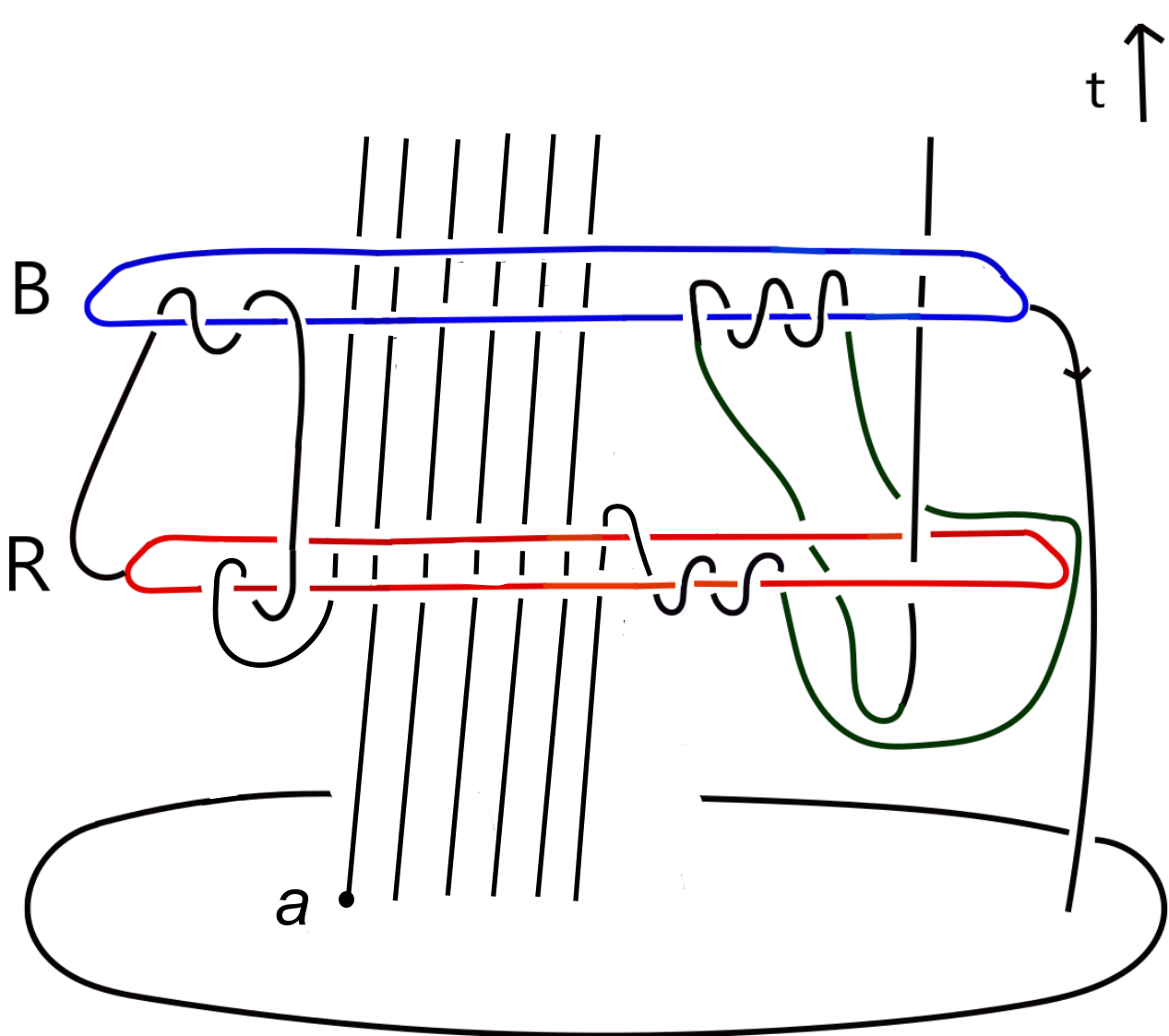}
    \caption{The barbell manifold $t^{-1}\nu_R {\nu_B}^3{\nu_R}^{-3}t^{-6}\nu_R{\nu_B}^2$.}
    \label{wordexample}
\end{figure}

\begin{figure}
    \centering
    \includegraphics[width=0.4\linewidth]{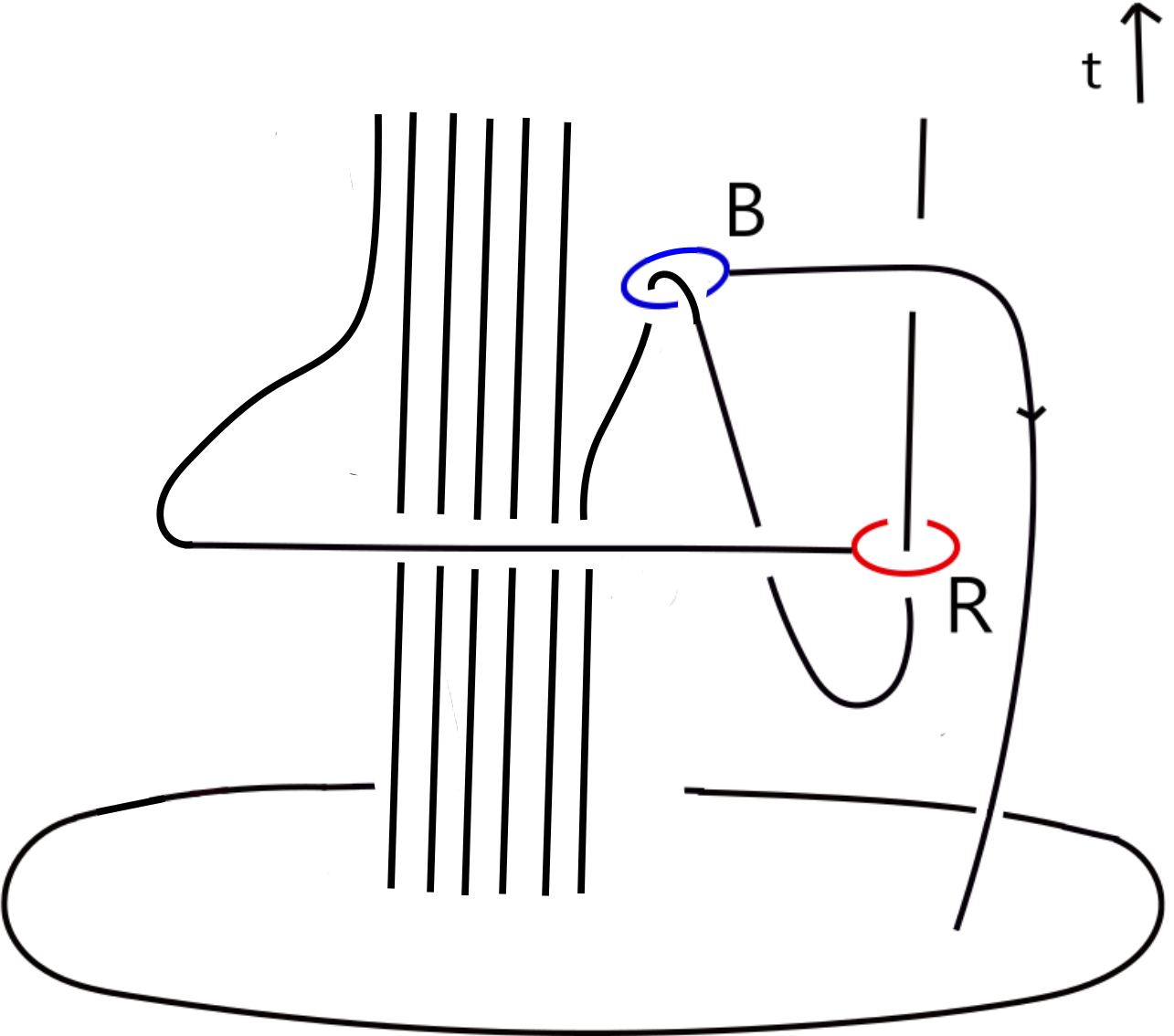}
    \caption{The barbell $\theta_k(e_7,-e_1)$.}
    \label{changebarbell1}
\end{figure}

\begin{theorem}
\label{barbellformulaem111}
   The $W_3$ of the barbell represented by the word 
$$ t^{j_1}\cdots \nu_R^{1}t^{j_{\alpha}}\cdots t^{j_{\beta}}\nu_B^{1}t^{j_{\beta+1}}\cdots t^{j_n}$$
is given by
\begin{align*}
&(|j_{\beta+1}|+\cdots+|j_{n}|)\mathrm{T}_1(j_1+\cdots+j_{\beta},j_{\alpha}+\cdots j_{n})+
(|j_{\alpha}|+\dots+|j_{\beta}|)\mathrm{T}_2(j_1+\cdots+j_{\beta},j_{\alpha}+\cdots j_{n})+\\
&(|j_1|+\dots+|j_{\alpha-1}|)\mathrm{T}_3(j_1+\cdots+j_{\beta},j_{\alpha}+\cdots j_{n})
\end{align*}
and the $W_3$ of the barbell represented by the word 
$$ t^{j_1}\cdots \nu_B^{1}t^{j_{\alpha}}\cdots t^{j_{\beta}}\nu_R^{1}t^{j_{\beta+1}}\cdots t^{j_n}$$
is given by
\begin{multline*}
(|j_1|+\dots+|j_{\alpha-1}|)\mathrm{T}_4(j_1+\cdots+j_{\alpha-1},j_{\beta+1}+\cdots j_n)+(|j_{\beta+1}|+\cdots+|j_n|)\mathrm{T}_6(j_1+\cdots+j_{\alpha-1},j_{\beta+1}+\cdots j_n).
\end{multline*}
Here $\mathrm{T}_1$, $\mathrm{T}_2$, $\mathrm{T}_3$, $\mathrm{T}_4$ and $\mathrm{T}_6$ are defined in Theorems \ref{mainformula} and \ref{newcasecomplete}.
\end{theorem}

\begin{proof}
    With out loss of generality, we can choose a scanning disk $\{pt\}\times D^3\subset S^1\times D^3$ such that there are exactly $|j_1|+\dots+|j_n|$ intersection points between the bar and the scanning disk. 
    
    For the former case, observe that there are $|j_1|+\dots+|j_{\alpha-1}|$ Type 3 intersection points, $|j_{\alpha}|+\dots+|j_{\beta}|$ Type 2 intersection points and $|j_{\beta+1}|+\cdots+|j_{n}|$ Type 1 intersection points. It follows that its $W_3$ can be calculated using the formula in Theorem \ref{mainformula} and plugging in $j_1+\cdots+j_{\beta}$ and $j_{\alpha}+\cdots j_{n}$.

    Similarly, for the latter case, observe that by performing an appropriate isotopy as in Lemma \ref{type4and6barbell} (for example by dragging the blue cuff along the bar in the direction of reducing the total unsigned number of times the bar loops around the circle factor by $t^{j_\alpha}\cdots t^{j_{\beta}}$ times), one obtains an equivalent barbell represented by the word $ t^{j_1}\cdots \nu_B^{1}\nu_R^{1}t^{j_{\beta+1}}\cdots t^{j_n}$. This equivalent barbell admits $|j_1|+\dots+|j_{\alpha-1}|$ Type 4 intersection points and $|j_{\beta+1}|+\cdots+|j_n|$ Type 6 intersection points and its $W_3$ is given by plugging in $j_1+\cdots+j_{\alpha-1}$ and $j_{\beta+1}+\cdots j_n)$.
\end{proof}

\begin{figure}
    \centering
    \includegraphics[width=0.37\linewidth]{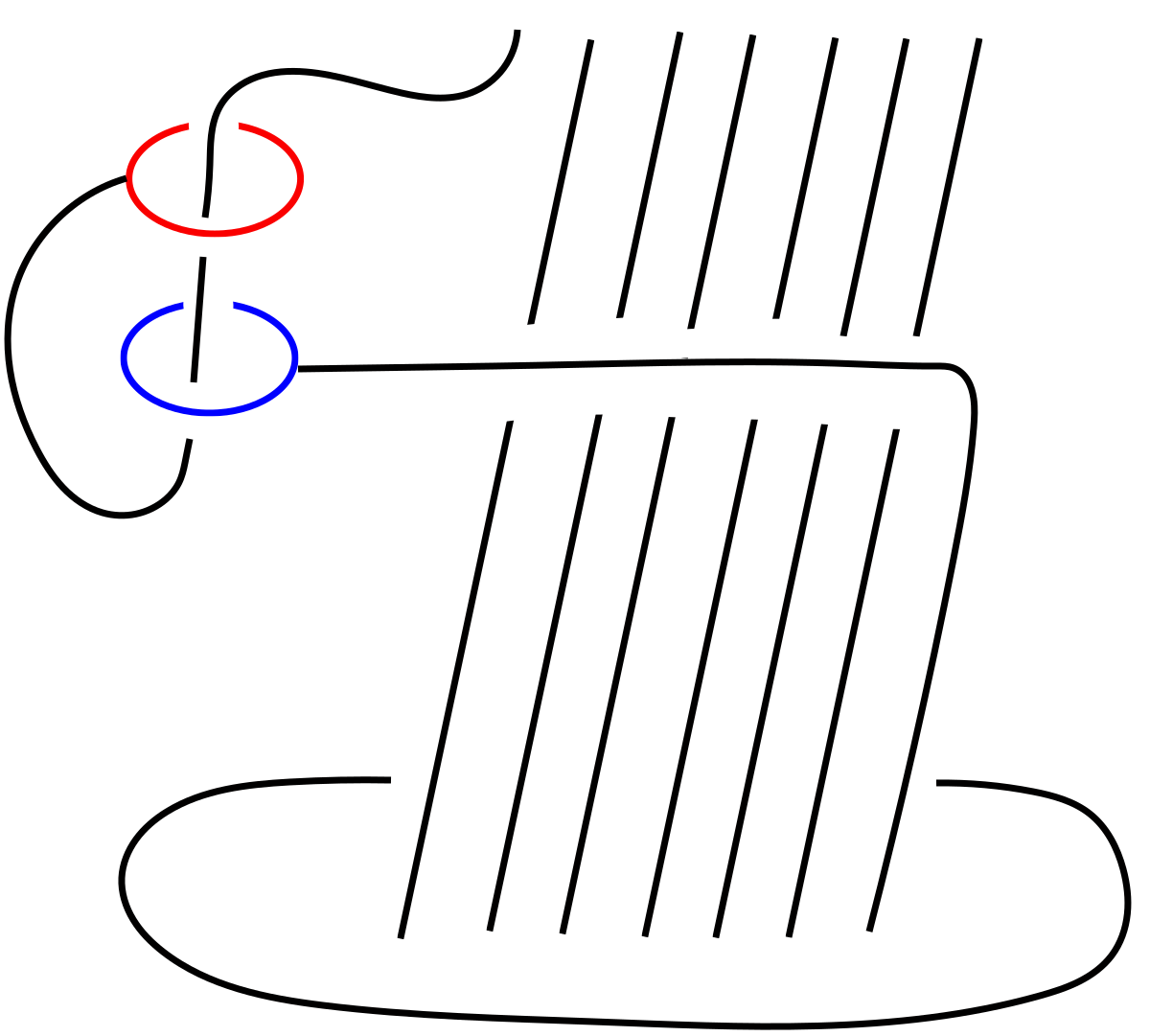}
    \caption{The barbell $t^{-7} \nu_{R}\nu_{B}$.}
    \label{08092024-2}
\end{figure}

\begin{example}
    Take the barbell represented by the word $t^{-1}\nu_R {\nu_B}^3{\nu_R}^{-3}t^{-6}\nu_R{\nu_B}^2$ as in Figure \ref{wordexample}. Note that we have $|t^{-1}\nu_R {\nu_B}^3{\nu_R}^{-3}t^6\nu_R{\nu_B}^2|=-7$. Scanning through a chosen $D^3$ as in Figure \ref{wordexample} gives rise to 7 intersection points that can be analyzed separately. For example, at the first intersection point $a$ (counting from left to right), we can similarly draw its intersection with the spanning disks of the cuff spheres as in Figure \ref{wordexample2}.

    \begin{figure}
    \centering
    \includegraphics[width=0.4\textwidth]{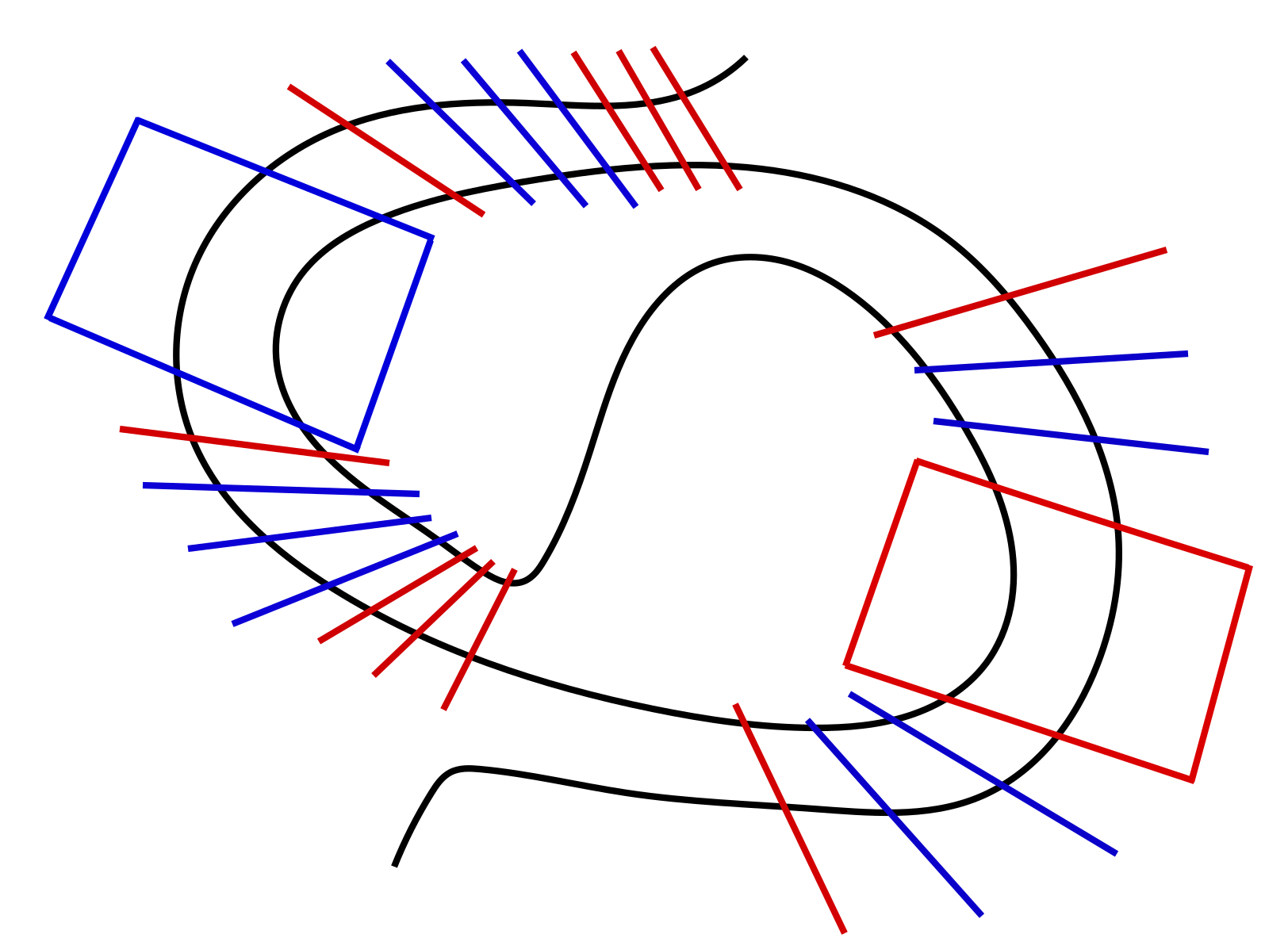}
    \caption{The intersection between the arc through point $a$ and the spanning disks of the cuff spheres.}
    \label{wordexample2}
\end{figure}

    In this case, one takes all possible pairs of two symmetric red arcs with two symmetric blue arcs, and every such pair leads to exactly the same pattern of linking number diagrams. The $W_3$ invariant of $t^{-1}\nu_R {\nu_B}^3{\nu_R}^{-3}t^{-6}\nu_R{\nu_B}^2$ is the sum of all these possible pairs at each intersection point of the bar with the scanning disk $D^3$. In particular, the above barbell can be factored into the following sub-barbells:
    \begin{center}
         \begin{tabular}{ c c c}
           $3  t^{-1}\nu_{R}\nu_{B}t^{-6}$ & $2  t^{-1}\nu_{R}t^{-6}\nu_{B}$ & $9 t^{-1}\nu_{B}\nu_{R}^{-1} t^{-6}$\\
         $6  t^{-1}\nu_{R}^{-1}t^{-6}\nu_{B}$&  $3  t^{-1} \nu_{B}t^{-6}\nu_{R}$ & $2  t^{-7} \nu_{R}\nu_{B}$.
    \end{tabular}
    \end{center}

For example, Figure \ref{changebarbell1} shows the barbell $t^{-1}\nu_{R}\nu_{B}t^{-6}$. In fact, this barbell is isotopic to the barbell $\theta_7(e_6,e_1)$ by dragging the blue cuff $B$ in Figure \ref{changebarbell1} downwards along the bar to the $t$ negative direction for one time (note that crossing changes of the bar are allowed as we work in dimension 4). Therefore, its $W_3$ invariant can be calculated using Theorem \ref{mainformula}. Further, one can verify, by drawing pictures, that except for
$  t^{-7} \nu_{R}\nu_{B},$
each of the above sub-barbells can be identified with some $\theta_k(\pm e_i,\pm e_j)$ for some $i,j\in \{1,\dots,7\}$, and $k\leq 8$. In particular, we have 
\begin{align*}
   & W_3(t^{-1}\nu_R {\nu_B}^3{\nu_R}^{-3}t^{-6}\nu_R{\nu_B}^2)=\\
   &3W_3(\theta_8(e_6,e_7))+2W_3(\theta_8(e_6,e_7))+9W_3(\theta_8(-e_7,e_1))+\\
   &6W_3(\theta_8(-e_6,e_7)+
   3W_3(\theta_8(e_1,e_1))+2W_3(t^{-7} \nu_{R}\nu_{B}),
\end{align*}
in which all terms can be easily calculated using Theorem \ref{mainformula} except $W_3(t^{-7} \nu_{R}\nu_{B})$ (See Figure \ref{08092024-2}).

The barbell $t^{-7} \nu_{R}\nu_{B}$ is not isotopic to a $\theta_k(v,w)$. However, we can calculate its $W_3$ invariant fairly easily by observing that its bar intersects $s_0\times D^3$ at $7$ points and all of these 7 points give rise to tetrahedron pictures with the same shape as a Type 3 intersection point shown in Figure \ref{type3qqq}. In particular, we can make use of Figure 17 of \cite{Budney-gabai2} together with the strategy we have been discussing throughout this section to obtain a slightly variant formula as in Theorem \ref{mainformula}. Namely, we take the formula for Type 3 points in Theorem \ref{mainformula}, and plug in $i=j=7$.
\end{example}

\section{Computing $(W_3)_m$ for $m\geq 2$}
In this section, we illustrate how one computes $(W_3)_m$ for $m\geq 2$ with a focus on $m=2$. We begin with a brief discussion of cubical complexes in 5.1, followed by introducing the notion of co-geodesic submanifolds as a generalization of collinear submanifolds (cf. Section 4.2) in 5.2. Finally, in 5.3, we compute $(W_3)_2$ for unknotted barbells in $\natural_2 S^1\times D^3$.
\subsection{CAT(0)-cubical complexes}
In this section, we review the notion of CAT(0)-cubical complexes that will be needed in later sections. Most of the material can be found in \cite{cat0}. 

We start with a brief motivation. Recall from Lemma \ref{collineardetect} that $W_3$ can be calculated using a special class of submanifolds of the universal cover of $C_3[S^1\times D^3]$. In this chapter, we will calculate a version of $W_3$ for $\pi_0\mathrm{Diff}(\natural_m S^1\times D^3)$, for $m\geq 2$, using a similar method and a special class of submanifolds of the universal cover of $C_3[\natural_m S^1\times D^3]$, called the \textit{co-geodesic submanifolds} (cf. Definition \ref{co-geodesicdef}) whose definitions depend on the uniqueness and existence of geodesics. In particular, we need the fact that the universal cover of $C_3[\natural_m S^1\times D^3]$ is a uniquely geodesic space. We shall achieve this by giving it a CAT(0) cubical complex structure (cf. Theorem \ref{cube} and Lemma \ref{universalcovercube}).

Let $C^n=[0,1]^n$ denote the $n$-dimensional cube. The (codimension-1) faces of $C^n$ are given by $F_{i,\epsilon}=\{x=(x_1,\dots,x_n)\in C^n\colon x_i=\epsilon\}$ for $\epsilon=1$ or $0$. Lower-dimensional faces are obtained by intersections of the codimension 1 faces. We equip cubes with the restricted Euclidean metric.

\begin{definition}
    Let $\mathcal{C}$ be a set of cubes and $\mathcal{S}$ a set of isometries between faces of the cubes called gluing maps such that no cube is glued to itself, and there are at most two distinct gluings between any two distinct cubes. Then such a pair defines a space $\bigcup_{C\in \mathcal{C}} C/\sim $ where $\sim$ is the equivalence relation induced by the gluing maps. This is called a cubical complex.

\end{definition}

\begin{definition}
    Let $X$ be a cubical complex and $x,y\in X$. A string $\Sigma$ from $x$ to $y$ is a sequence of points $x=x_0,x_1,\dots,x_m=y$ such that each pair of consecutive points $x_i$ and $x_{i+1}$ in the sequence is contained in a single cube $C_i$. The length $L(\Sigma)$ of a string $\Sigma$ is given by the sum $\sum_{i} d_{C_i}(x_i,x_{i+1})$ where $d_{C_i}$ is the Euclidean metric on $C_i$. If $X$ is string connected, meaning that any two points are connected by a string of intervals, then we equip it with the polyhedric metric:
    $$d(x,y)=\mathrm{inf}\{L(\Sigma)\colon \text{ $\Sigma$ is a string from $x$ to $y$}\}.$$
\end{definition}
The polyhedric metric can be equivalently defined as
$$d(x,y)=\mathrm{inf}\{l(\gamma)\colon \gamma \text{ is a rectifiable curve in $X$ from $x$ to $y$}\}.$$
where a rectifiable curve $\gamma\colon [0,1]\to X$ means a curve with finite length 
$$l(\gamma)\coloneqq\mathrm{sup}_{0=t_0\leq \cdots\leq t_n=1} \sum_{i=1}^{n-1} d(\gamma(t_i),\gamma(t_{i+1})).$$

\begin{definition}
    A cubical complex is \textit{finite dimensional} if there is an upper bound on the dimension of cubes.
\end{definition}

\begin{proposition}
\label{cube2}
    Finite-dimensional cubical complexes are complete geodesic spaces.
\end{proposition}

\begin{definition}
\label{links}
    The link of a vertex $v$ of a cubical complex $X$ is the induced simplicial complex obtained from $X$ by taking the space $S(v,\epsilon)=\{x\in X\colon d(x,v)=\epsilon\}$ for a small $\epsilon\geq 0$. A link is called \textit{flag} if there are no empty simplices, i.e.~whenever the 1-skeleton of a simplex exists, so does the entire simplex (i.e.~the higher dimensional faces all exist).
\end{definition}

The next theorem will be crucial for us:
\begin{theorem}
\label{cube}

    A cubical complex $X$ is CAT(0) if and only if it is simply-connected and all links are flag.
\end{theorem} 

This is called \textit{Gromov’s link condition} due to Gromov. In particular, this implies that such cubical complexes are uniquely geodesic spaces. Both Proposition \ref{cube2} and Theorem \ref{cube} can be found in Section 3 of \cite{cat0}.

We apply the above to the universal cover of $\natural_m S^1\times D^3$ (cf. Section 3.1 and Figure \ref{fundamental domain}).

\begin{figure}
    \centering
    \includegraphics[width=0.3\textwidth]{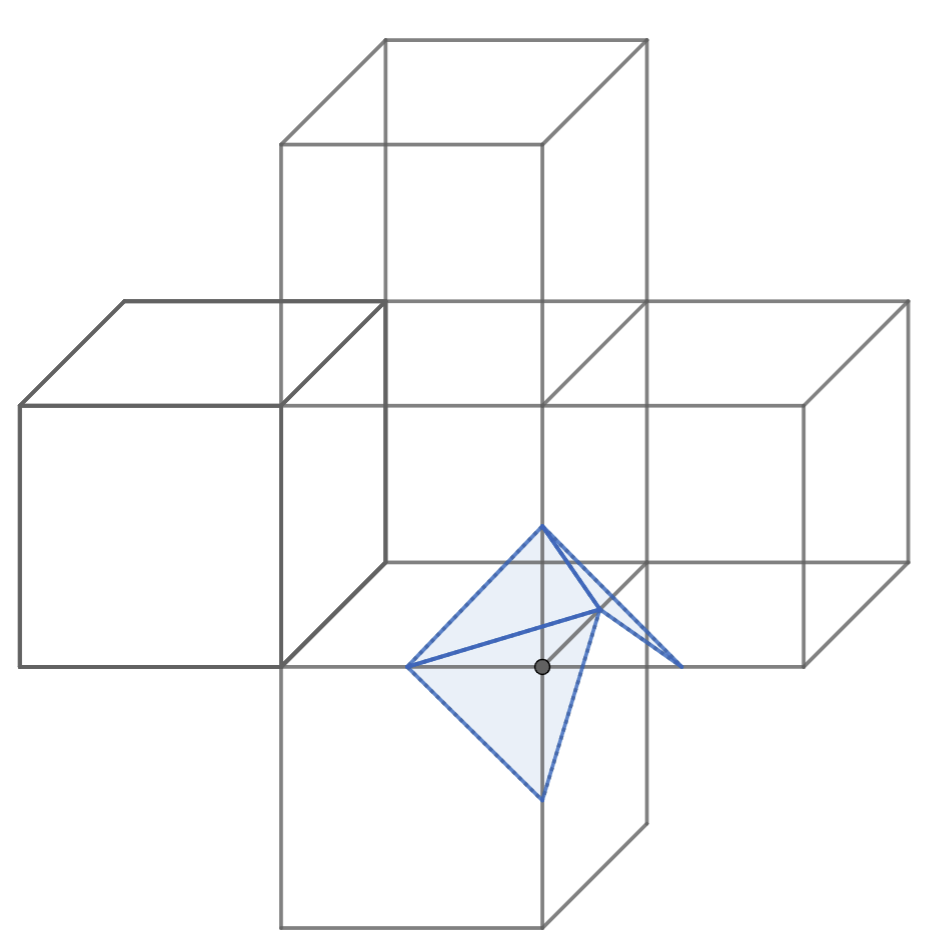}
    \caption{The links of the universal cover $T_2\times I\times I$.}
    \label{cat000}
\end{figure}

\begin{lemma}
    \label{universalcovercube}
    The space $T_2\times I\times I$ admits a natural cubical complex structure (cf. Section 5.1) which contains $5$ 4-dimensional cubes in each fundamental domain. Furthermore, all links of this cubical complex are flag. Thus, it is a uniquely geodesic metric space whose geodesics are concatenations of straight intervals. 
\end{lemma}

\begin{proof}
One can verify that each link of $T_2\times I\times I$ is a union of three ordered tetrahedrons, say $(\mathcal{T}_1,\mathcal{T}_2,\mathcal{T}_3)$ with $\mathcal{T}_i$ being glued to $\mathcal{T}_{i+1} $ along a common face, for $i=1,2$. This space satisfies the condition of Definition \ref{links}. See Figure \ref{cat000} for a 3-dimensional picture with a typical link drawn in blue. It is a union of three triangles with two consecutive triangles glued along a common edge. One can observe that the links of all the vertices in this picture are isomorphic.
  
The second paragraph follows from Proposition \ref{cube2} and Theorem \ref{cube}. \end{proof}

\begin{remark}
    In fact, one can go further and prove that Lemma \ref{universalcovercube} is true for all $m\geq 2$, meaning that the space $T_m\times I\times I$ admits a natural cubical complex structure (cf. Section 5.1), which contains $3m-1$ 4-dimensional cubes in each fundamental domain, and the links of all vertices are isomorphic to a union of $m+1$ tetrahedrons with consecutive ones being glued along a common face. These are all flag links.
\end{remark}
 
We can also view $T_m\times I\times I$ as a Riemannian manifold with corners. Locally, the Riemannian metric is given by the induced Euclidean metric in each fundamental domain.

\subsection{Co-geodesic submanifolds}

We end this section with a discussion of ideas for calculation of $(W_3)_m$. One refers to Section 3.4 for comparison to calculation of $W_3$. Fix a basepoint $(x_0,0)\in U_m\times \mathbb{R} = T_m\times I \times \mathbb{R}$ where $U_m$ is the universal tree. For a word $\nu$ in the free group generated by $(t_1)_i,(t_2)_i,(t_3)_i$ with $i=1,2\dots m$ (cf. the beginning of this section), we define the cohorizontal submanifold 
$$\nu.\mathrm{Co}_1^2=\{(p_1,p_2)\in {\overline{C_2(\natural_m S^1\times D^3)}}\colon \text{$p_2$ and $\nu.p_1$ coincide on $T_m$,} 
\text{and $p_2-\nu. p_1 =\lambda(x_0,0,1)$ for $\lambda>0$} \}.$$
Note that the subtraction makes sense since it only concerns the $I\times \mathbb{R}$ part. This submanifold $\nu.\mathrm{Co}_1^2$ detects the element $\nu.\omega_{12}$ as it intersects the image of $\nu.\omega_{12}$ at exactly one point, and does not intersect $\nu'.\omega_{12}$ if $\nu \neq \nu'$. 

To detect elements in $\pi_5 C_3[\natural_m S^1\times D^3]\otimes \mathbb{Q}$ as Whitehead products of elements in the form of $\nu.\omega_{12}$ and $\mu.\omega_{23}$ with $\nu$ and $\mu$ being words in the free group generated by the $3m$ generators $(t_1)_i$, $(t_2)_i$ and $(t_3)_i$ for $i=1,2,\dots, m$, we make use of the language of \textit{co-geodesic} submanifolds as a generalisation of Budney--Gabai's collinear submanifolds (cf. Section 3.4). In the universal cover $T_m\times I^2$ of $\natural_m S^1\times D^3$, straight lines are not natural choices anymore as we are not working in a convex space. However, since $T_m\times I^2$ is a simply-connected cubical complex with all links being flag, it is a uniquely geodesic space by Theorem \ref{cube}, with geodesics given by concatenations of straight arcs (polylines). This motivates the following definition.

\begin{definition}
\label{co-geodesicdef}
    The co-geodesic submanifolds are defined as
    $$\mathrm{Cog_{\nu,\mu}}^1=\{(p_1,p_2,p_3)\in T_m\times I^2\colon 
    \text{ $\nu.p_1$ lies on the geodesic determined by $\mu.p_3$ and $p_2$}\}$$
$$\mathrm{Cog_{\nu,\mu}}^3=\{(p_1,p_2,p_3)\in T_m\times I^2\colon
\text{ $\mu.p_3$ lies on the geodesic determined by $\nu.p_1$ and $p_2$}\}$$
where $\nu$ and $\mu$ are words in the free group generated by the $3m$ generators $(t_1)_i$, $(t_2)_i$ and $(t_3)_i$ for $i=1,2,\dots, m$.
\end{definition}

The manifold $\mathrm{Cog_{\nu,\mu}}^1$ intersects $\mu.\omega_{12}$ algebraically once, and $\mathrm{Cog_{\nu,\mu}}^3$ intersects algebraically $\nu.\omega_{23}$ once. Therefore, the pair $$(\mathrm{Cog_{\nu,\mu}}^1,\mathrm{Cog_{\nu,\mu}}^3)$$ detects the element $\mu \nu.[\omega_{12},\omega_{23}]$ in the sense that the preimage of this pair under $\mu \nu .[\omega_{12},\omega_{23}]$ is a standard linking pair in $S^5$ (cf. Section 3.4), and the preimage of the same pair under $\mu' \nu' [\omega_{12},\omega_{23}]$ is empty if $\mu \nu [\omega_{12},\omega_{23}]\neq \mu' \nu' [\omega_{12},\omega_{23}]$. For any element in $\pi_5 C_3[\natural_m S^1\times D^3]\otimes \mathbb{Q}$, we can write it as a linear combination of $\mu\nu.[\omega_{12},\omega_{23}]$ with coefficients given by words in the free group with $3m$ generators, and the linking number of the preimage of $(\mathrm{Cog_{\nu,\mu}}^1,\mathrm{Cog_{\nu,\mu}}^3)$ determines the coefficients. Thus computing $(W_3)_m$ boils down to a linking number calculation.

To calculate these linking numbers, we use a trick similar to Budney--Gabai's argument as in Lemma 3.4 of \cite{Budney-gabai2} (see also Lemma \ref{cohorizontal calculation}). Recall that $T_m$ is the universal cover of a regular neighbourhood of a wedge of $m$ circles in $\mathbb{R}^2$ (cf. Figure \ref{fundamental domain} and the discussion around it). Embed $T_m\times I^2$ naturally in $T_m\times I \times \mathbb{R}$ and consider the function $P_{\lambda}\colon T_m\times I \times \mathbb{R}\to T_m\times I \times \mathbb{R}$ defined by $$P_{\lambda}(p_1,p_2,p_3,p_4)=(p_1,p_2,p_3, p_4+\lambda d(p_1,p_2,p_3)^2)$$ where $\lambda >0$ is a positive parameter and $d(p_1,p_2,p_3)$ denotes the distance from the point $(p_1,p_2,p_3)$ to $x_0$ under the induced Euclidean metric. This 1-parameter family of diffeomorphisms sends polylines to parabolas that get steeper as $\lambda$ increases positively. As $\lambda$ goes to infinity, the polylines breaks into steep line segments. It follows that the linking number of the preimage of $(\mathrm{Cog_{\nu,\mu}}^1,\mathrm{Cog_{\nu,\mu}}^3)$ can be calculated by the pair $$(\mu\mathrm{Co}_2^1-\mu \nu^{-1} \mathrm{Co}_3^1,\nu \mu^{-1}\mathrm{Co}_1^3-\nu \mathrm{Co}_2^3)$$ (cf. Lemma \ref{cohorizontal calculation}). Here the cohorizontal submanifolds $\mathrm{Co}_i^j$ are defined in the same manner as before with point $i$ and $j$ being cohorizontal and the remaining coordinate is allowed to vary freely. This will enable us us to calculate the $(W_3)_m$ invariant for embedded barbells in $\natural_m S^1\times D^3$ for $m\geq 2$ using cohorizontal submanifolds in the next section.

\subsection{Unknotted barbells in $S^1\times D^3 \natural S^1\times D^3$}
From now on, we fix $m=2$ and study $\pi_0\mathrm{Diff}(\natural_2 S^1\times D^3)$.  In this section, we use the invariant 
$$(W_3)_2\colon \pi_0\mathrm{Diff}(S^1\times D^3 \natural S^1\times D^3,\partial)\to \Lambda_2$$ to study the mapping class group of $S^1\times D^3 \natural S^1\times D^3$ by detecting embedded barbells. In this case, we simplify our notation and just refer $(W_3)_2$ as $W_3$ when there is no ambiguity. We also denote the generators of $\pi_1C_3(S^1\times D^3 \natural S^1\times D^3)$ by $t_1,t_2,t_3,u_1,u_2,u_3.$ We discuss how to compute $W_3$ for unknotted barbells in $S^1\times D^3\natural S^1\times D^3$ that are specified by the quotient $$\mathrm{F}_4[\nu_R,\nu_B,t,u]\to \langle \nu_B\rangle \backslash  F_4/\langle \nu_R\rangle$$
where $t$, $u$ represent the circle factors and $\nu_R$, $\nu_B$ represent the meridians of $R$ and $B$ respectively. We view the boundary connected sum as attaching a handle diffeomorphic to $I\times D^3$. We fix a slice $\{\mathrm{pt}\}\times D^3\subset I\times D^3$ as the scanning disk $\Delta$ such that both $B$ and $R$ are disjoint from the scanning disk, and are contained in the same component (i.e.~on the same side). Let $\#w$ be the number of times the bar of $\mathcal{B}(w)$ passes through $\Delta$ that cannot be eliminated without moving $R$ and $B$ through $\Delta$, i.e.~we exclude pairs of intersections that can be eliminated by an inverse finger move. We also assume that we choose $\Delta$ so that $\#w$ is minimal. First, we note that the following relation analogous to Theorem \ref{barbellfactor} holds.
\begin{lemma}
\label{reallinearfactorisation}
    Let $w\in \langle \nu_R\rangle \backslash  F_4/\langle \nu_B\rangle$ and let $\mathcal{B}(w)$ denote the induced barbell from $w$. Then 
    $$ w=\nu_B^{i_1}\lambda_1\nu_R^{k_1}\cdots \nu_B^{i_n}\lambda_n\nu_R^{k_n}$$ 
    where $\lambda_1,\dots,\lambda_n\in \langle t,u \rangle$. Then $W_3(\mathcal{B}(w))$ is equal to the sum
    $$\sum_{\alpha,\beta} c_{i_{\alpha} k_{\beta}}W_3(\mathcal{B}_{i_{\alpha} k_{\beta}})$$
    of $W_3$ invariants of a finite sequence of ``sub-barbells'' $\mathcal{B}_{i_{\alpha} k_{\beta}}$ with $\alpha, \beta \in \{1,2,\dots,n\}$, and coefficients $c_{i_{\alpha} k_{\beta}}$. Each $\mathcal{B}_{i_{\alpha} k_{\beta}}$ is an unknotted barbell satisfies the condition that the bar links each of the two cuff spheres only positively once, and is represented by a word in the form
    $$\lambda_1 \lambda_2\dots \lambda_{\alpha-1}\nu_B^{1}\lambda_{\alpha}\dots \lambda_{\beta}\nu_R^{ 1}\dots \lambda_n$$
    when $\alpha\leq \beta$, and in the form
     $$\lambda_1 \lambda_2\dots \lambda_{\beta-1}\nu_R^{1}\lambda_{\beta}\dots \lambda_{\alpha}\nu_B^{ 1}\dots \lambda_n$$
    when $\alpha>\beta$. Further, the coefficients satisfies
    $$c_{i_{\alpha} k_{\beta}}=i_{\alpha} k_{\beta}.$$
\end{lemma}
\begin{proof}
    The proof is the same as the proof of Theorem \ref{barbellfactor} except that we replace the powers of $t$'s by monomials expressed as products of $\lambda$'s that still indicate the path one needs to travel (along the bar) along $t$- and $u$-directions of $S^1\times D^3\natural S^1\times D^3$ from one point of a cohorizontal pair to the other.
\end{proof}

Let $w\in \langle \nu_R\rangle \backslash  F_4/\langle \nu_B\rangle$ such that both $\nu_R$ and $\nu_B$ appear only once positively. We classify such barbells into two separate cases:
\begin{enumerate}
    \item $w=a\nu_R^1b\nu_Bc$ where $a,b,c\in \langle t,u \rangle.$
    \item $w=a\nu_B^1b\nu_Rc$ where $a,b,c\in \langle t,u \rangle.$
\end{enumerate}
For the latter case, by applying an appropriate isotopy if necessary (for example, dragging the blue cuff along the bar until it approaches the red cuff), we can change $w$ to the word
 $$b\nu_B^{1}\nu_R^{ 1}cb=a\nu_B^1\nu_R^1c'.$$
We make the important observation that the pattern of cohorizontal points of such barbells coincide with barbells in$S^1\times D^3$, and can be presented in exactly the same triangle and tetrahedron pictures as in the last chapter. For a word $k(t,u)\in  \langle t,u \rangle$, we denote the inverse of $k$ as $\bar{k}\coloneqq k^{-1}.$ Then we have:

\begin{lemma}
\label{type123formula}
    If $w=a\nu_Rb\nu_Bc$ where $a,b,c\in \langle t,u \rangle$, then we notice that $\# (a\nu_Rb\nu_Bc )= \#(\bar{c}\nu_R\bar{b}\nu_B\bar{a})$, and these intersection points are classified into Type 1, Type 2 and Type 3 intersection points as in Figures \ref{twisted arc}, \ref{type2qqq} and \ref{type3qqq} since they share the same pattern of triangle and tetrahedron pictures. In particular, we have 
     $$W_3(a\nu_Rb\nu_Bc)=k_1 T_1(\bar{b}\bar{a},\bar{c}\bar{b})+k_2 T_2(\bar{b}\bar{a},\bar{c}\bar{b})+k_3 T_3(\bar{b}\bar{a},\bar{c}\bar{b})$$
where $k_1+k_2+k_3=\#w$. Further, $k_1=l_3$, $k_3=l_1$ and $k_2=l_2$. Each constituent of a Type 1 intersection point contributes
     \begin{enumerate}
    \item$ \mathrm{lk}(\textcolor{red}{\mu\mathrm{Co}_2^1},\textcolor{blue}{\nu \mu^{-1}\mathrm{Co}_1^3})=0$ 

    \item $\mathrm{lk}(\textcolor{blue}{\mu\mathrm{Co}_2^1},\textcolor{red}{\nu \mu ^{-1}\mathrm{Co}_1^3}) =a_1b_1\bar{c}_3\bar{b}_3a_3 b_3$

    \item -$\mathrm{lk}(\textcolor{red}{\mu\mathrm{Co}_2^1},\textcolor{blue}{\nu\mathrm{Co}_2^3})=-\bar{c}_1\bar{b}_1\bar{b}_3\bar{a}_3$

    \item -$\mathrm{lk}(\textcolor{blue}{\mu\mathrm{Co}_2^1},\textcolor{red}{\nu\mathrm{Co}_2^3})=-\bar{b}_1\bar{a}_1\bar{c}_3\bar{b}_3$

        \item $\mathrm{lk}(\textcolor{red}{\mu \nu^{-1}\mathrm{Co}_3^1},\textcolor{blue}{\nu\mathrm{Co}_2^3})=\bar{c}_1\bar{b}_1a_1b_1a_3b_3$

    \item $\mathrm{lk}(\textcolor{blue}{\mu \nu^{-1}\mathrm{Co}_3^1},\textcolor{red}{\nu\mathrm{Co}_2^3})=0$
\end{enumerate}
leading to
$$T_1(\bar{b}\bar{a},\bar{c}\bar{b})=a_1b_1\bar{c}_3\bar{b}_3a_3 b_3+\bar{c}_1\bar{b}_1a_1b_1a_3b_3-\bar{c}_1\bar{b}_1\bar{b}_3\bar{a}_3-\bar{b}_1\bar{a}_1\bar{c}_3\bar{b}_3.$$
Note that the subscripts (1 or 3) mean that we viewing an element of $\langle t,u \rangle$ as an element in $\langle t_1,u_1 \rangle$ or $\langle t_3,u_3 \rangle$. Similarly, for a Type 2 intersection point we have 
     \begin{enumerate}
    \item$ \mathrm{lk}(\textcolor{red}{\mu\mathrm{Co}_2^1},\textcolor{blue}{\nu \mu^{-1}\mathrm{Co}_1^3})=-b_1c_1\bar{b}_3\bar{a}_3b_3c_3$

    \item $\mathrm{lk}(\textcolor{blue}{\mu\mathrm{Co}_2^1},\textcolor{red}{\nu \mu ^{-1}\mathrm{Co}_1^3}) =a_1b_1\bar{c}_3\bar{b}_3a_3b_3$

    \item -$\mathrm{lk}(\textcolor{red}{\mu\mathrm{Co}_2^1},\textcolor{blue}{\nu\mathrm{Co}_2^3})=0$

    \item -$\mathrm{lk}(\textcolor{blue}{\mu\mathrm{Co}_2^1},\textcolor{red}{\nu\mathrm{Co}_2^3})=0$

        \item $\mathrm{lk}(\textcolor{red}{\mu \nu^{-1}\mathrm{Co}_3^1},\textcolor{blue}{\nu\mathrm{Co}_2^3})=\bar{c}_1\bar{b}_1a_1b_1a_3b_3$

    \item $\mathrm{lk}(\textcolor{blue}{\mu \nu^{-1}\mathrm{Co}_3^1},\textcolor{red}{\nu\mathrm{Co}_2^3})=-\bar{b}_1\bar{a}_1b_1c_1b_3c_3$
    \end{enumerate}
    that leads to
    $$T_2(\bar{b}\bar{a},\bar{c}\bar{b})=-b_1c_1\bar{b}_3\bar{a}_3b_3c_3+a_1b_1\bar{c}_3\bar{b}_3a_3b_3+\bar{c}_1\bar{b}_1a_1b_1a_3b_3-\bar{b}_1\bar{a}_1b_1c_1b_3c_3$$
    and for a Type 3 intersection point
         \begin{enumerate}
    \item$ \mathrm{lk}(\textcolor{red}{\mu\mathrm{Co}_2^1},\textcolor{blue}{\nu \mu^{-1}\mathrm{Co}_1^3})=-b_1c_1\bar{b}_3\bar{a}_3b_3c_3+\bar{b}_1\bar{a}_1b_1c_1b_3c_3$

    \item $\mathrm{lk}(\textcolor{blue}{\mu\mathrm{Co}_2^1},\textcolor{red}{\nu \mu ^{-1}\mathrm{Co}_1^3}) =a_1b_1\bar{c}_3\bar{b}_3a_3b_3+\bar{b}_1\bar{a}_1b_3c_3\bar{b}_3\bar{a}_3$

    \item -$\mathrm{lk}(\textcolor{red}{\mu\mathrm{Co}_2^1},\textcolor{blue}{\nu\mathrm{Co}_2^3})=-b_1c_1a_3b_3$

    \item -$\mathrm{lk}(\textcolor{blue}{\mu\mathrm{Co}_2^1},\textcolor{red}{\nu\mathrm{Co}_2^3})=-a_1b_1b_3c_3$

        \item $\mathrm{lk}(\textcolor{red}{\mu \nu^{-1}\mathrm{Co}_3^1},\textcolor{blue}{\nu\mathrm{Co}_2^3})=b_1c_1\bar{b}_1\bar{a}_1\bar{b}_3\bar{a}_3+\bar{c}_1\bar{b}_1a_1b_1a_3b_3$

    \item $\mathrm{lk}(\textcolor{blue}{\mu \nu^{-1}\mathrm{Co}_3^1},\textcolor{red}{\nu\mathrm{Co}_2^3})=-a_1b_1\bar{c}_1\bar{b}_1\bar{c}_3\bar{b}_3$
\end{enumerate}
that leads to 
$$T_3(\bar{b}\bar{a},\bar{c}\bar{b})=-b_1c_1\bar{b}_3\bar{a}_3b_3c_3+a_1b_1\bar{c}_3\bar{b}_3a_3b_3+\bar{b}_1\bar{a}_1b_3c_3\bar{b}_3\bar{a}_3-b_1c_1a_3b_3-a_1b_1b_3c_3+$$
$$b_1c_1\bar{b}_1\bar{a}_1\bar{b}_3\bar{a}_3+\bar{c}_1\bar{b}_1a_1b_1a_3b_3-a_1b_1\bar{c}_1\bar{b}_1\bar{c}_3\bar{b}_3.$$
\end{lemma}
\begin{lemma}
\label{type46formula}
    If $w=a\nu_Bb\nu_Rc$ where $a,b,c\in \langle t,u \rangle.$, then as argued previously, we can assume $b=0$. Therefore, there are only Type 4 and Type 6 intersection points (as in Figure \ref{newtypes}) that are involved in its $W_3$, and
$$W_3(a\nu_B\nu_Rc)=k_1 T_4(\bar{a},\bar{c})+k_2 T_6(\bar{a},\bar{c})$$
with $k_1+k_2=\#w $ and $k_1=l_2$, $k_2=l_1$ For a Type 4 intersection point, we have  \begin{enumerate}
    \item$ \mathrm{lk}(\textcolor{red}{\mu\mathrm{Co}_2^1},\textcolor{blue}{\nu \mu^{-1}\mathrm{Co}_1^3})=0$

    \item $\mathrm{lk}(\textcolor{blue}{\mu\mathrm{Co}_2^1},\textcolor{red}{\nu \mu ^{-1}\mathrm{Co}_1^3}) =\bar{a}_1\bar{c}_3\bar{a}_3+a_1\bar{c}_3a_3$

    \item -$\mathrm{lk}(\textcolor{red}{\mu\mathrm{Co}_2^1},\textcolor{blue}{\nu\mathrm{Co}_2^3})=-\bar{c}_1\bar{a}_3-\bar{c}_1a_3$

    \item -$\mathrm{lk}(\textcolor{blue}{\mu\mathrm{Co}_2^1},\textcolor{red}{\nu\mathrm{Co}_2^3})=\bar{a}_1\bar{c}_3+a_1\bar{c}_3$

        \item $\mathrm{lk}(\textcolor{red}{\mu \nu^{-1}\mathrm{Co}_3^1},\textcolor{blue}{\nu\mathrm{Co}_2^3})=-\bar{c}_1\bar{a}_1\bar{a}_3-\bar{c}_1a_1a_3$

    \item $\mathrm{lk}(\textcolor{blue}{\mu \nu^{-1}\mathrm{Co}_3^1},\textcolor{red}{\nu\mathrm{Co}_2^3})=0$
\end{enumerate}
that leads to 
$$T_4(\bar{a},\bar{c})=\bar{a}_1\bar{c}_3\bar{a}_3+a_1\bar{c}_3a_3-\bar{c}_1\bar{a}_3-\bar{c}_1a_3+\bar{a}_1\bar{c}_3+a_1\bar{c}_3-\bar{c}_1\bar{a}_1\bar{a}_3-\bar{c}_1a_1a_3$$
and for a Type 6 intersection point, we have 
 \begin{enumerate}
    \item$ \mathrm{lk}(\textcolor{red}{\mu\mathrm{Co}_2^1},\textcolor{blue}{\nu \mu^{-1}\mathrm{Co}_1^3})=\bar{c}_1\bar{a}_3\bar{c}_3+c_1\bar{a}_3c_3$

    \item $\mathrm{lk}(\textcolor{blue}{\mu\mathrm{Co}_2^1},\textcolor{red}{\nu \mu ^{-1}\mathrm{Co}_1^3}) =-a_1c_3a_3-a_1\bar{c}_3a_3+\bar{a}_1c_3\bar{a}_3+\bar{a}_1\bar{c}_3\bar{a}_3$

    \item -$\mathrm{lk}(\textcolor{red}{\mu\mathrm{Co}_2^1},\textcolor{blue}{\nu\mathrm{Co}_2^3})=\bar{a}_1\bar{c}_3+\bar{a}_1c_3$

    \item -$\mathrm{lk}(\textcolor{blue}{\mu\mathrm{Co}_2^1},\textcolor{red}{\nu\mathrm{Co}_2^3})=-\bar{c}_1\bar{a}_3-c_1\bar{a}_3$

        \item $\mathrm{lk}(\textcolor{red}{\mu \nu^{-1}\mathrm{Co}_3^1},\textcolor{blue}{\nu\mathrm{Co}_2^3})=-\bar{c}_1\bar{a}_1\bar{a}_3-c_1\bar{a}_1\bar{a}_3+\bar{c}_1a_1a_3+c_1a_1a_3$

    \item $\mathrm{lk}(\textcolor{blue}{\mu \nu^{-1}\mathrm{Co}_3^1},\textcolor{red}{\nu\mathrm{Co}_2^3})=-\bar{a}_1\bar{c}_1\bar{c}_3-\bar{a}_1c_1c_3$
\end{enumerate}
that leads to 
$$T_6(\bar{a},\bar{c})=\bar{c}_1\bar{a}_3\bar{c}_3+c_1\bar{a}_3c_3-a_1c_3a_3-a_1\bar{c}_3a_3+\bar{a}_1c_3\bar{a}_3+\bar{a}_1\bar{c}_3\bar{a}_3+\bar{a}_1\bar{c}_3+\bar{a}_1c_3$$
$$-\bar{c}_1\bar{a}_3-c_1\bar{a}_3-\bar{c}_1\bar{a}_1\bar{a}_3-c_1\bar{a}_1\bar{a}_3+\bar{c}_1a_1a_3+c_1a_1a_3-\bar{a}_1\bar{c}_1\bar{c}_3-\bar{a}_1c_1c_3.$$
\end{lemma}

Combining Lemmas \ref{reallinearfactorisation}, \ref{type123formula} and \ref{type46formula} enables us to calculate $W_3$ for any random word $w\in \langle \nu_B\rangle \backslash  F_4/\langle \nu_R\rangle$.

\begin{figure}
    \centering
    \includegraphics[width=0.6\linewidth]{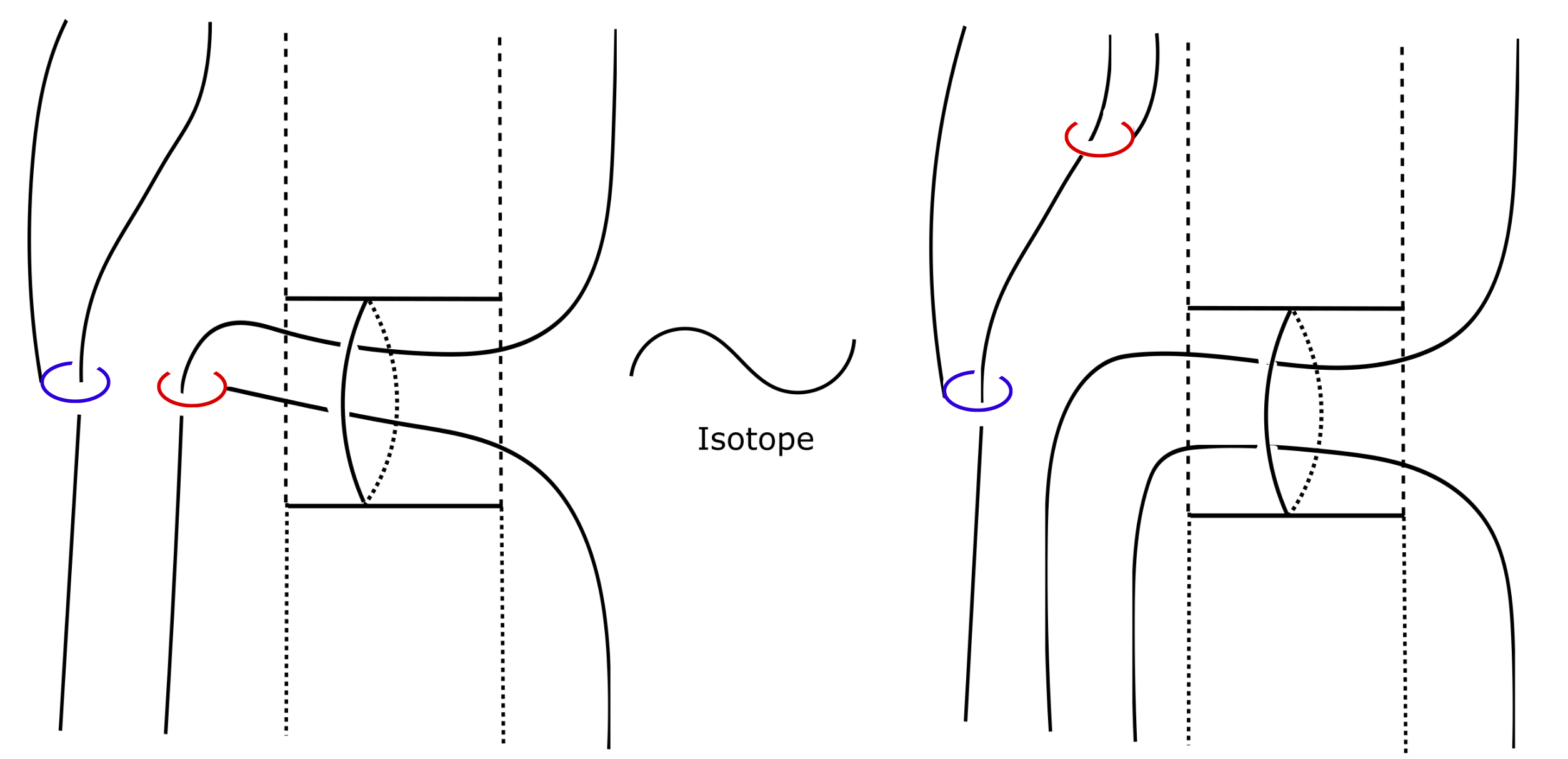}
    \caption{The barbells $\mathcal{B}(t\nu_B t\nu_R u)$ and $\mathcal{B}(t\nu_B\nu_R tut^{-1})$.}
    \label{extrafigureexample}
\end{figure}

Consider the example as in Figure \ref{extrafigureexample} that shows an unknotted barbell in $S^1\times D^3\natural S^1\times D^3$ denoted by $\mathcal{B}(t\nu_B t\nu_R u).$ The cylinder in the middle of Figure \ref{extrafigureexample} indicates the ``neck'' of the boundary connected-sum. Recall from the last section that to define $(W_3)_2$, we need to choose a scanning disk properly embedded in $S^1\times D^3 \natural S^1 \times D^3$. There is a natural choice $\Delta$ in the middle of the boundary connected sum neck, drawn in Figure \ref{extrafigureexample}. By assumption, the cuff spheres are disjoint from $\Delta$. The induced barbell diffeomorphism makes $\Delta=I\times I\times I$ into an element in $\pi_2 \mathrm{Emb}(I,S^1\times D^3 \natural S^1\times D^3)$ via the scanning map.

\begin{proposition}
\label{extraexampleprop}
   The barbell diffeomorphisms induced from $\mathcal{B}(t\nu_B t\nu_R u)$ (see Figure \ref{extrafigureexample}) is not isotopic to the identity in $S^1\times D^3\# S^1\times D^3$. Moreover, the barbells $\mathcal{B}(t\nu_B t\nu_R u^k)$ for $k=1,2,3,\dots$ induce linearly independent, non-isotopic barbell diffeomorphisms of $S^1\times D^3\# S^1\times D^3$.
\end{proposition}

\begin{proof}
First, we isotope $\mathcal{B}(t\nu_B t\nu_R u)$ to $\mathcal{B}(t\nu_B\nu_R tut^{-1})$. Scanning leads to two Type 4 intersection points and applying Lemma \ref{type46formula} we obtain:
   $$2(t_1^{-1}t_3u_3^{-1}t_3^{-2}+t_1t_3u_3^{-1}-t_1u_1^{-1}t_1^{-1}t_3^{-1}-t_1u_1^{-1}t_1^{-1}t_3+t_1^{-1}t_3u_3^{-1}t_3^{-1}+t_1t_3u_3^{-1}t_3^{-1}-t_1u_1^{-1}t_1^{-2}t_3^{-1}-t_1u_1^{-1}t_3)\neq 0.$$
Thus $\Phi_{ \mathcal{B}(t\nu_B t\nu_R u)}$ non-trivial in $\pi_0\mathrm{Diff}(S^1\times D^3\# S^1\times D^3,\partial)$. Similarly, applying the formulae to $\mathcal{B}(t\nu_B t\nu_R u^k)$ gives rise to 
   $$2(t_1^{-1}t_3u_3^{-k}t_3^{-2}+t_1t_3u_3^{-k}-t_1u_1^{-k}t_1^{-1}t_3^{-1}-t_1u_1^{-k}t_1^{-1}t_3+t_1^{-1}t_3u_3^{-k}t_3^{-1}+t_1t_3u_3^{-k}t_3^{-1}-t_1u_1^{-k}t_1^{-2}t_3^{-1}-t_1u_1^{-k}t_3)\neq 0.$$
Now, suppose the $W_3$ of the linear combination $$a_1 \mathcal{B}(t\nu_B t\nu_R u)+a_2 \mathcal{B}(t\nu_B t\nu_R u^2) +a_3 \mathcal{B}(t\nu_B t\nu_R u^3) +\cdots+a_n \mathcal{B}(t\nu_B t\nu_R u^n)=0.$$
Then setting all terms that contain powers of $u_1^{\pm 2}$ and $u_3^{\pm 2}$ being trivial implies that $a_1=0$. Repeating this argument with $u_1^{\pm k}$ and $u_3^{\pm k}$, $k=2,3,\dots,n-1$ implies that all coefficients are zero.
\end{proof}

Therefore, we have obtained an infinite sequence of knotted 3-balls $\Phi_{\mathcal{B}(t\nu_B t\nu_R u^k)}(\Delta)$ in $S^1\times D^3\natural S^1\times D^3$ for $k\geq 1$. This also shows that these elements do not come from $\pi_0 \mathrm{Diff}(S^1\times D^3,\partial)\times \pi_0 \mathrm{Diff}(S^1\times D^3,\partial)$ since it fixes $\Delta$. So, we have the following theorem.
\begin{theorem}
\label{chapter5main1}
The group $\pi_0\mathrm{Diff}(S^1\times D^3 \natural S^1\times D^3,\partial)/ \left( \pi_0 \mathrm{Diff}(S^1\times D^3,\partial)\right)^2$ has an infinitely generated subgroup. Moreover, $\Phi_{ \mathcal{B}(t\nu_B t\nu_R u^k)}(\Delta)$ gives infinitely many properly embedded separating 3-balls in $S^1\times D^3 \natural S^1\times D^3$ with common boundaries that are not isotopic relative to the boundary.
    \end{theorem}

\begin{remark}
Tatsuoka \cite{tatsuoka} shows that there is an infinite sequence of non-trivial diffeomorphisms $\tilde{\delta}_k$ of $S^1\times D^3\#S^1\times D^3$ that are supported in $S^1\times D^3\natural S^1\times D^3\subset S^1\times D^3\#S^1\times D^3$, and agree with $\delta_k$ in $S^1\times D^3$ (cf. Theorem \ref{deltakkkk}) away from a neighbourhood of the connected-sum sphere. Her argument goes as follows. Suppose that $\tilde{\delta}_k$ is isotopic to the identity, then gluing $S^2\times D^2$ to the right copy of $S^1\times D^3$ modifies $(\tilde{\delta}_k,S^1\times D^3\#S^1\times D^3)$ to $(\delta_k,S^1\times D^3\# S^4)$ and implies that $\delta_k$ is isotopic to the identity, which is a contradiction.

We observe that this argument does not say anything about the non-triviality of $\mathcal{B}(t\nu_B t\nu_R u^k)$. In particular, if we glue a copy of $S^2\times D^2$ to the right copy of $S^1\times D^3\# S^1\times D^3$,  $\mathcal{B}(t\nu_B t\nu_R u^k)$ becomes a trivial barbell in $S^1\times D^3$. 

\end{remark}

\clearpage

\bibliographystyle{abbrv} 
\bibliography{sample} 
\bigskip
\noindent\textsc{Yau Mathematical Sciences Center, Tsinghua University}\\
\textit{Email:} \texttt{weizheniu@mail.tsinghua.edu.cn}
\end{document}